\documentclass{cimart}

\usepackage{subcaption}

\usepackage{amsmath}
\usepackage{mathrsfs}    
\usepackage{amssymb}
\usepackage{amsthm}
\usepackage{graphicx,color}
\usepackage{hyperref}

\usepackage{verbatim}

\usepackage{pdfsync}

\definecolor{verde}{rgb}{0,0.35,0.1} 
\definecolor{rosso}{rgb}{0.7,0,0}
\definecolor{blue}{rgb}{0,0,1}
\definecolor{viola}{rgb}{0.6,0,0.4}

\newcommand{\mf}[1]{ \mathbf{#1}}

\usepackage{titlesec} 
\titleformat{\chapter}{\normalfont\huge}{\thechapter.}{20pt}{\huge\bf}  

\newtheoremstyle{definition2}{\topsep}{\topsep}%
     {}
     {}
     {\bfseries}
     {.}
     {.5em}
     {\thmnumber{(#2)}\thmname{ #1}\thmnote{ #3}}

\usepackage{indentfirst}

\def\N{\mathbb{N}}

\def\R{\mathbb{R}}

\def\hbar{\bar{h}}

\usepackage{bbm}

\def\dist{{\rm dist}}

\def\leq{\leqslant}
\def\geq{\geqslant}

\newcommand{\rst}[1]{\ensuremath{{\mathbin |}%
\raise-.5ex\hbox{$#1$}}}

\usepackage{mathtools}

 \numberwithin{equation}{section}

\makeatletter
\DeclareRobustCommand\widecheck[1]{{\mathpalette\@widecheck{#1}}}
\def\@widecheck#1#2{%
    \setbox\z@\hbox{\m@th$#1#2$}%
    \setbox\tw@\hbox{\m@th$#1%
       \widehat{%
          \vrule\@width\z@\@height\ht\z@
          \vrule\@height\z@\@width\wd\z@}$}%
    \dp\tw@-\ht\z@
    \@tempdima\ht\z@ \advance\@tempdima2\ht\tw@ \divide\@tempdima\thr@@
    \setbox\tw@\hbox{%
       \raise\@tempdima\hbox{\scalebox{1}[-1]{\lower\@tempdima\box
\tw@}}}%
    {\ooalign{\box\tw@ \cr \box\z@}}}
\makeatother

\title{%
    Topics in elliptic problems: from semilinear equations to shape optimization
    }

\author{%
    Hugo Tavares
    }

\authorinfo[%
    H. Tavares]{
Centro de An\'alise Matem\'atica, Geometria e Sistemas Din\^amicos, Departamento de Matem\'atica do Instituto Superior T\'ecnico,  Universidade de Lisboa, 1049-001 Lisboa, Portugal}{%
hugo.n.tavares@tecnico.ulisboa.pt
    }

\abstract{%
 In this paper, which corresponds to an updated version of the author's \emph{Habilitation lecture} in Mathematics, we do an overview of several topics in elliptic problems. We review some old and new results regarding the Lane-Emden equation, both under Dirichlet and Neumann boundary conditions, then focus on sign-changing solutions for Lane-Emden systems. We also survey some results regarding fully nontrivial solutions to gradient elliptic systems with mixed cooperative and competitive interactions. We conclude by exhibiting results on optimal partition problems, with cost functions either related with Dirichlet eigenvalues or to the Yamabe equation. Several open problems are referred along the text.
    }



\VOLUME{32}
\YEAR{2024}
\ISSUE{3}
\firstpage{171}
\DOI{https://doi.org/10.46298/cm.12568}

\begin{document}

\tableofcontents

\section{Introduction}

Elliptic partial differential equations are a very important class of equations with connections to applied sciences (e.g. physics or biology) as well as to other fields of Mathematics such as Differential Geometry, Functional Analysis and Calculus of Variations. Because of these facts, they are a fascinating topic and an increasingly active field of research.

The purpose of this text is to give an overview of most of my research from the last 10 years in this field, as well as to point out new directions and open problems. The text corresponds to a slightly updated version of the \emph{Habilitation lecture} in Mathematics I presented at Instituto Superior T\'ecnico, Universidade de Lisboa, on May 22--23, 2023.

\smallbreak

I have divided this lecture into four parts. 

Section \ref{sec2}, entitled \emph{Semilinear elliptic problems: old and new}, is mainly intended to give context to the sections. Therein, I review (for non-experts) some classical material related to variational methods and applications to the Lane-Emden equation,
\[
-\Delta u=|u|^{p-1} u\qquad  \text{ in some domain $\Omega\subset \mathbb{R}^N$},
\]
with $p>0$ and under homogeneous Dirichlet and Neumann boundary conditions. Even though most of the material is classical, the section finishes with some open questions and some recent contributions made by myself. This section contains (with permission) parts of my survey paper \cite{CIM}.

In the remaining sections I focus mostly on more recent problems and on my research. Section \ref{sec3} deals with \emph{Elliptic Hamiltonian systems}, in particular with Lane-Emden systems
\[
-\Delta u=|v|^{q-1}v,\qquad -\Delta v=|u|^{p-1}u \qquad \text{ in } \Omega.
\]
for $p,q>0$. Under Dirichlet boundary conditions $u=v=0$ on $\partial \Omega$, we focus on least energy nodal solutions in the subcritical case, while under Neumann boundary conditions (where all nontrivial solutions are necessarily sign-changing) we focus on least energy solutions, both in the subcritical and the critical cases. We are mostly interested in existence result, as well as symmetry or symmetry breaking phenomena in case the underlying domain is a bounded radial set.

Section \ref{sec4}, entitled \emph{Existence of fully nontrivial solutions to a class of gradient elliptic systems,} focus in systems of type
\begin{equation}\label{S-system}
\begin{cases}
-\Delta u_{i}+\lambda_{i}u_{i}=u_{i}|u_i|^{p-2}\sum_{j = 1}^{d}\beta_{ij}  |u_{j}|^p   ~\text{ in } \Omega,\\
u_{i}=0 \text{ on } \partial\Omega,  \quad i=1,...d,
\end{cases}
\end{equation}
with $d\geq 2$ equations, where $N\geq 3$, in a (Sobolev) critical or subcritical regime $0<p\leq2^*/2=N/(N-2)$,  and $\lambda_i \in \R$. We assume that the coupling terms are symmetric, that is $\beta_{ij}=\beta_{ji}$ for $i\neq j$, which provides a gradient structure to the problem. These systems admit semitrivial solutions, that is, solutions $(u_1,\ldots, u_d)$ with zero components (some $u_i$ might vanish identically), and we focus on reviewing conditions given in the literature that ensure the existence of fully nontrivial solutions (in particular of least energy). We conclude with a short subsection describing some results on normalized solution, a very active topic of research at the moment, dealing with the less explored case of bounded domains.

Finally, Section \ref{sec5} deals with \emph{Optimal partition problems}, which are problems of the following type:
\[
\inf\left\{\Phi(\omega_1,\ldots, \omega_m):\ \omega_i\in \mathcal{A},\ \omega_i\cap \omega_j=\emptyset\ \forall i\neq j\right\},
\]
where $\mathcal{A}$ is a class of admissible sets in a certain ambient space and $\Phi:\mathcal{A}^m\to \R$ is a cost function. We address the case of spectral partitions (where the cost function is related to Dirichlet eigenvalues) and of problems related with the Yamabe equation. We are interested in obtaing existence and regularity results, as well as in the relation with system \eqref{S-system}. Very recent results regarding optimal partition problems with volume constraints are also described.

\smallbreak

Part of my works did not make it into this overview for the sake of brevity and coherence of the material. I would like to emphasize, for instance, the work with D. Bonheure, J. F\"oldes, E. Moreira dos Santos and A. Salda\~na about general criteria for the uniqueness of critical points of functionals with a hidden convexity \cite{BFMST18}, the joint work with D. Cassani and J. Zhang on gradient systems with critical growth  in the sense of Moser in dimension two \cite{CassaniTavaresZhang}, the work with E. Moreira dos Santos, G. Nornberg and D. Schiera on the phenomenon of two principal half eigenvalues in the context of fully nonlinear Lane-Emden type systems with possibly unbounded coefficients and weights \cite{LEfullynonlinear}, and the work with F. Agostinho and S. Correia \cite{AgostinhoCorreiaTavares} on a rather complete study of the positive $H^1$-solutions of the equation $-u'' +\lambda u = |u|^{p-2}u$ on the $\mathcal{T}$-metric graph for $\lambda>0$ and $p>2$.

\section{Semilinear elliptic problems: old and new}\label{sec2}

\noindent This section contains parts of my survey paper \cite{CIM}.

\smallbreak

Many problems can be modelled with the aid of elliptic partial differential equations. One of the most well known examples is the classical Poisson equation: given a bounded regular domain $\Omega \subset \R^N$, take
\[
-\Delta u=f \text{ in } \Omega.
\]
Its solutions may represent  the shape of an elastic membrane in equilibrium subject to a vertical load  $f:\Omega \to \R$ ($u(x)$ corresponds to the vertical displacement at the point $x$); an electrostatic potential (for $f=\rho/\varepsilon$, where $\rho(x)$ is the volume charge density and $\varepsilon$ the permittivity of the medium), a gravitational potential (for $f=-4\pi G \rho$, where $\rho$ is the density of the object and $G$ the gravitational constant), or the stationary solutions for the heat equation (in this case, $u$ represents a temperature, and $f$ is a heat source).
 To obtain existence and uniqueness of solution, one couples the equation with boundary conditions: \emph{Dirichlet boundary conditions} ($u=g$ on $\partial \Omega$) or \emph{Neumann boundary conditions} ($u_\nu :=\nabla u\cdot \nu=g$ on $\partial \Omega$, where $\nu=\nu(x)$ is the outer  unit  vector at $x\in \partial \Omega$) are typical examples arising in applications. Linear problems are very well understood and can be found in classical textbooks (see for instance \cite{EvansPDE,SalsaPDE}), while current research aims at a good understanding of nonlinear problems. Among the wide class of possible nonlinearities, the simplest to treat (although already quite rich mathematically, as we will see), are semilinear problems, where $f:\R\to \R$, $f=f(u)$, is a nonlinear function, that is, the nonlinearity occurs at the level of the zero order terms.

Let $N\geq 3$ and  $p>0$. For simplicity, let us work from now on with the prototypical example of the Lane-Emden equation
\begin{equation}\label{eq:LE_example1}
-\Delta u=|u|^{p-1}u \text{ in } \Omega.
\end{equation}
Equation \eqref{eq:LE_example1} is not only a mathematical paradigm in nonlinear analysis of PDEs, but it has also a physical motivation, since, for $N=3$, normalized radial solutions of \eqref{eq:LE_example1} solve the Lane-Emden equation of index $p$, namely,
\[
-\frac{1}{\xi^{2}}(\xi^2 \theta')'=|\theta|^{p-1}\theta,\quad \theta(0)=1, \theta'(0)=0.
\]
The latter equation is used in astrophysics to model self-gravitating spheres of plasma, such as stars or self-consistent stellar systems in polytropic-convective equilibrium, where the pressure $P$ and the density $\rho\left(=k \theta^p\right)$ satisfy a nonlinear relationship $P=c \rho^{\frac{p+1}{p}}$, see \cite{C57}. In this setting, a positive solution $\theta$ is often called a polytrope, and it contains, up to constants, important physical information such as the radius of the star (the first zero $r_1$ of $\left.\theta\right)$, the total mass $\int_{B_{r_1}} \theta^p$, the pressure $\theta^{p+1}$ and for an ideal gas, the temperature is proportional to $\theta$.

The (homogeneous) Dirichlet case has been extensively studied:
\begin{equation}\label{eq:LE_example}
-\Delta u=|u|^{p-1}u \text{ in } \Omega,\quad u=0 \text{ on } \partial \Omega,
\end{equation}
and in the following four subsections we focus on this problem in the case where $\Omega$ is a bounded regular domain. We remark that the study of the Neumann case has been completed only recently with our contributions, see Section \ref{sec:Neumann_1eq} below.

\smallbreak

Clearly $u\equiv 0$ is always a solution, but we are interested in nontrivial ones. Recall that (formally) a weak solution is a function $u\in H^1_0(\Omega)$ such that
\[
\int_\Omega \nabla u\cdot \nabla v -\int_\Omega |u|^{p-1}u v=0\qquad  \forall v \in H^1_0(\Omega),
\]
and weak solutions correspond to critical points of the functional
\begin{equation}\label{eq:thefunctional}
\mathcal{I}:H^1_0(\Omega)\to \R,\qquad \mathcal{I}(u)=\frac{1}{2}\int_\Omega |\nabla u|^2 - \frac{1}{p+1}\int_\Omega |u|^{p+1}
\end{equation}
(observe that $(|t|^{p+1})'=(p+1)|t|^{p-1}t$ for every $t\in \R$, $p>0$). To make these statements precise and correct, we need a restriction on the exponent $p$, since the integral $\int_\Omega |u|^{p+1}$ is not always finite for $u\in H^1_0(\Omega)$. One needs to recall \emph{Sobolev inequalities}:  for
\[
 1\leq q\leq 2^*:=\frac{2N}{N-2},
 \]
 there exists $C_{N,q}>0$  such that
\begin{equation}\label{eq:Sobemb}
\left(\int_\Omega |u|^q\right)^{1/q} \leq C_{N,q} \left(\int_\Omega |\nabla u|^2\right)^{1/2}\ \forall u\in H^1_0(\Omega),\
\end{equation}
which amounts to saying that the embedding $H^1_0(\Omega)\hookrightarrow L^q(\Omega)$ is continuous. The number $2^*$ is the critical Sobolev exponent. Therefore, in conclusion, \eqref{eq:thefunctional} is well defined only for $p\leq 2^*-1$. In this case, in order to look for weak solutions of the problem \eqref{eq:LE_example}, one may try to find critical points of $\mathcal{I}$.

Now the question is: how can we find a critical point of $\mathcal{I}$? And how many of them are there? The answer depends on $p$: not only the geometry of $\mathcal{I}$ changes from $p<1$ to $p>1$, but also the situations $p+1<2^*$ and $p+1=2^*$ are very different:  the embedding of $H^1_0(\Omega)$ in $L^q(\Omega)$ is compact only for $1\leq q<2^*$. The discussion of the case $p>2^*-1$ is much harder and less is known.

Before moving on, we would like to point out that most of the results we describe below for $p\neq 1$ do not depend on the homogeneity of the map $t\mapsto |t|^{p-1}t$. Indeed, many results are true for more general nonlinearities $f(t)$.  We also recall that the main purpose of this section is to give some context to the forthcoming ones (even though I present some new results of my own here, see the last paragraph of Subsection \ref{subsec:sublinear} and Subsection \ref{sec:Neumann_1eq}). For this reason, I do not even dare to make a complete state of the art and many references are left out.

\subsection{Dirichlet boundary conditions: the linear case $p=1$.} Before going nonlinear, let us analyse what happens in the linear case $p=1$, that is: $-\Delta u=u$ in  $\Omega$, $u=0$ on $\partial \Omega$. This problem may or may not have a (nontrivial) solution; what we are asking, in other words, is if $\lambda=1$ is an eigenvalue of the operator $A:=-\Delta$ with Dirichlet boundary conditions.
In this context, indeed, $\lambda\in \R$ is called an eigenvalue whenever  $-\Delta u=\lambda u$  in $\Omega$, $u=0$  on $\partial \Omega$ admits a nontrivial (weak) solution. From the spectral theory of compact operators (using the compactness of the embedding $H^1_0(\Omega) \hookrightarrow L^2(\Omega)$), we deduce that the eigenvalues of $-\Delta$ (counting multiplicities) form a nondecreasing sequence
\[
0<\lambda_1(\Omega) <\lambda_2(\Omega)\leq \lambda_3(\Omega)\leq \ldots \to \infty
\]
and there is a Hilbert base of $H^1_0(\Omega)$ made of associated eigenfunctions $(v_n)_{n}$. Exactly as for eigenvalues of a matrix, the eigenvalues admit a variational formulation, namely the Courant-Fischer-Weyl min-max formulas
\begin{equation}\label{eq:eigenvalues_charact}
\lambda_1(\Omega) =\min\left\{ \mathcal{R}(u):\ u\in H^1_0(\Omega)\setminus \{0\}\right\}, \quad \lambda_k(\Omega)= \min_{V\subset H^1_0(\Omega)\atop  \textrm{dim} V=k} \max_{u\in V\setminus \{0\}}\mathcal{R}(u)\quad  (k\geq 2),
\end{equation}
where $\mathcal{R}(u)=\int_\Omega |\nabla u|^2 / \int_\Omega u^2$, for $u\not\equiv 0$, is called the Rayleigh quotient. The details can be found, for instance, in \cite[Chapter 6]{SalsaPDE}. Therefore, the question of whether problem \eqref{eq:LE_example} in the case $p=1$ admits a nontrivial solution or not depends on the domain: the answer is affirmative only for domains for which $1=\lambda_i(\Omega)$ for some $i$.

Below, in Section \ref{sec5}, I present my work on \emph{spectral partition problems}: optimal partition problems where the cost function depend on the eigenvalues of each set of the partition.

\subsection{Dirichlet boundary conditions: the sublinear case $0<p<1$.}\label{subsec:sublinear}
 If $0<p<1$, it is straightforward to see that $\mathcal{I}$ has a minimum in each direction: for a fixed $w\in H^1_0(\Omega)\setminus \{0\}$, this corresponds to studying the real function $f(t)=\mathcal{I}(t\omega)$, which has the form $at^2-b|t|^{p+1}$ for some $a,b>0$. Using Sobolev inequalities and the direct method of Calculus of Variations  \cite[Chapter 8.2]{EvansPDE}, one shows that $\mathcal{I}$ admits a global negative minimum in $H^1_0(\Omega)$: the level
\[
\inf\{\mathcal{I}(u):\ u\in H^1_0(\Omega)\}<0
\]
is achieved, providing a nontrivial solution (which is called a \emph{least energy solution}). We know a lot about minimizers. First of all, they are signed: either $u>0$ in $\Omega$ or $u<0$ in $\Omega$ (this is a consequence of the inequality $\mathcal{I}(|u|)\leq \mathcal{I}(u)$ and the strong maximum principle \cite[Chapter 2.2]{HanLin}). Positive solutions are unique \cite{Krasnoselskii}. This uniqueness also implies symmetry properties in symmetric domains: for instance, if the domain is radially symmetric (ball or annulus centered at the origin), the solution is radially symmetric (working in the space $H^1_{0,rad}:=\{u\in H^1_0(\Omega): u(x)=u(|x|)\ \forall x\in \Omega\}$ provides a positive solution). More generally, we can consider the situation of a domain $\Omega$ which is invariant under a subgroup $G$ of the orthogonal group $O(N)$.

In the previous paragraph we described properties of minimizers. Does $\mathcal{I}$ admit other critical points (\emph{i.e.}, solutions of the problem \eqref{eq:LE_example})? The answer is affirmative (see e.g. \cite{BartschWillem95}): there exists a sequence of critical points $(v_k)_k$ of $\mathcal{I}$, which satisfies
\[
\mathcal{I}(v_k)<0,\qquad \mathcal{I}(v_k)\to 0.
\]
This is a consequence of the $\mathbb{Z}_2$--symmetry of the problem (the functional is invariant under the map $u\mapsto -u$); solutions can be found as saddle points of $\mathcal{I}$, characterized via min-max methods in an analogous way to what happens for eigenvalues (recall \eqref{eq:eigenvalues_charact}). Observe that, since positive (and negative) solutions are unique, the previous multiplicity result yields the existence of infinitely many \emph{sign-changing} solutions. The next step is then to understand them as well as possible. The study of the zero-set of  sign-changing solutions (the \emph{free-boundary} set $\Gamma=\{x\in \Omega:\ u(x)=0\}$) is particularly challenging, as the function $f(t)=|t|^{p-1}p$ is not of class $C^1$ for $0<p<1$. The study of $\Gamma$ has been done recently: up to a subset with small Hausdorff dimension, $\Gamma$ is a regular hypersurface \cite{SoaveTerracini, SoaveWeth}. Moreover, one may also ask if, among all sign-changing solutions, there is one that minimizes the energy functional $\mathcal{I}$, that is, if the level
\[
c_{nod}=\inf\{\mathcal{I}(u):\ u \text{ is a sign-changing critical point of } \mathcal{I}\}
\]
 is achieved (solutions $u$ such that $\mathcal{I}(u)=c_{nod}$ are typically called \emph{least energy sign-changing solutions} or \emph{least energy nodal solutions}). The answer is affirmative, as shown recently in a joint work with D. Bonheure, E. Moreira dos Santos, E. Parini and T. Weth.
 \begin{theorem}[{{\cite{BSPTW}}}]
Let $0<p<1$. There exist $u\in H^1_0(\Omega)\cap L^\infty(\Omega)$ such that $\mathcal{I}(u)=c_{nod}$. Moreover, any function achieving the level $c_{nod}$ is a least energy nodal solution.
 \end{theorem}
In our  paper \cite{BSPTW} it is also shown that, quite remarkably, the type of critical point we find depends on the domain: there exist domains where the least energy nodal solution is a local minimizer of $\mathcal{I}$, and others where it is a saddle point (more precisely, of Mountain Pass type, see Theorem \ref{thm:MPT} below). A complete understanding of how the domain influences the type of critical point is an open problem. To conclude, we emphasize that the results in \cite{BSPTW} are valid for a class of nonlinearities $f(t)$  which also include the Allen-Cahn-type $f (t) = \lambda (t-|t|^{p-1}t)$, with $p > 1$ and $\lambda > \lambda_2(\Omega)$. For $C^1$ sublinear-type nonlinearities, we further prove that least energy nodal solutions on radial domains are not radially symmetric, but only foliated Schwarz symmetric, that is, there exists $p\in \partial B_1(0)$ such that the solution is  axially symmetric with respect to $p\R$, and strictly decreasing  in the polar angle $\theta=\arccos \left(\frac{x}{|x|}\cdot p\right)$.

\subsection{Dirichlet boundary conditions: the superlinear--subcritical case $1<p<2^*-1$.}

For the case $p>1$, by using Sobolev inequalities one can show that:
\begin{itemize}
\item[-] the origin $u=0$ is a strict local minimum of $\mathcal{I}$;
\item[-] $\mathcal{I}$ is unbounded from below and from above.
\end{itemize}
In this case, to obtain solutions we cannot simply minimize (nor maximize) the functional in the whole $H^1_0(\Omega)$. Based on the geometry of the functional, we can use the following version of the celebrated result by Ambrosetti and Rabinowitz \cite{AmbrosettiRabinowitz}.
\begin{theorem}[Mountain Pass Theorem]\label{thm:MPT}
Let $H$ be a Hilbert space and let $\mathcal{J}:H\to \R$ be a $\mathcal{C}^{1,1}$ functional satisfying
\begin{itemize}
\item $\mathcal{J}(0)=0$;
\item there exists $r>0$ such that
\[
\inf\{\mathcal{J}(u): \|u\|\leq r\}=0,\qquad \inf\{\mathcal{J}(u): \|u\|=r\}>0;
\]
\item there exists $v$ such that $\mathcal{J}(v)<0$.
\end{itemize}
Let $\Gamma:=\{\gamma\in C([0,1];H^1_0(\Omega)):\ \gamma(0)=0,\ \mathcal{J}(\gamma(1))<0\}$, and
\[
c:=\inf_{\gamma\in \Gamma} \sup_{u\in \gamma([0,1])} \mathcal{J}(u).
\]
Then there exists a sequence $(u_k)_k\subset H$ such that $\mathcal{J}(u_k)\to c$ and $\mathcal{J}'(u_k)\to 0$.
\end{theorem}
The proof of this result uses deformation lemmas and the study of steepest descending flows (a simple proof can be found in \cite[Chapter 8]{EvansPDE}). The existence of a sequence $(u_k)_k$ such that $\mathcal{J}(u_k)\to c$ and $\mathcal{J}'(u_k)\to 0$, by itself, does not imply the existence of a critical point (take the counterexample $H=\R$, $c=0$, $u_k=-k$ and $\mathcal{J}(x)=e^{x}$). A  new concept regarding compactness is needed:
\begin{quote}
A functional $\mathcal{J}\in C^1(H,\R)$ satisfies the \emph{Palais-Smale condition at $c$} if, whenever we have a sequence $(u_k)_k$ such that
$\mathcal{J}(u_k)\to c$ and  $\mathcal{J}'(u_k)\to 0$, then there exists a subsequence $(u_{k_j})_j$ of $(u_k)_k$ and $u\in H$ such that $u_{k_j}\to u$ in $H$. In particular, $\mathcal{J}'(u)=0$.
\end{quote}
Using the compactness of the Sobolev embeddings \eqref{eq:Sobemb} for $q<2^*$, one proves that $\mathcal{I}$ defined in \eqref{eq:thefunctional} satisfies this condition, and the Mountain Pass Theorem provides the existence of a critical point of $\mathcal{I}$, hence a solution of \eqref{eq:LE_example}. What can we now say about this solution? Another possible variational characterization is via a Nehari manifold:
\[
c=\inf_{\mathcal{N}}\mathcal{I},
\]
where
\[
\mathcal{N}=\{u\in H^1_0(\Omega)\setminus \{0\}:\ \mathcal{I}'(u)u=0\}=\{u\in H^1_0(\Omega)\setminus \{0\}:\ \int_\Omega |\nabla u|^2=\int_\Omega |u|^{p+1}\}.
\]
Observe that $\mathcal{I}$ is bounded from below on $\mathcal{N}$, and that the condition $I'(u)u=0$ is a free constraint (in the sense that the associated Lagrange multiplier is zero).  The solution achieving $c$ is also a \emph{least energy solution}, in the sense that
\[
c=\inf\{\mathcal{I}(u):\ u\in H^1_0(\Omega) \setminus \{0\},\ \mathcal{I}'(u)=0\}.
\]
Exactly as in the sublinear case, least energy solutions can be shown to be signed: they are either strictly positive or strictly negative in $\Omega$. However, uniqueness of positive solutions does not hold in general, as an effect of the topology of the domain (there are multiplicity results in annular domains) or of the geometry (dumbbell shaped domains). There is a long standing conjecture \cite{Dancer1988,Kawohl_uniq} that, if the domain is convex, then there is uniqueness of positive solution of \eqref{eq:LE_example} for $1<p<2^*-1$. A good review of the state-of-the-art regarding this subject can be found in the introduction of \cite{Grossi:2021aa}.
What about the symmetry in radial domains? When the domain is a ball, positive solutions are radially symmetric (consequence of the so called \emph{moving plane method}, which uses many types of maximum principles, see \cite{GidasNiNirenberg} or \cite[Chapter 2.6]{HanLin}). However, if $\Omega$ is an annulus, the solutions (at least for large $p$) lose one axis of symmetry, being foliated Schwarz symmetric \cite{BartschWethWillem} (axially symmetric and strictly decreasing with respect to the polar angle from the symmetry axis). As we can see, there are some key changes between the cases $p<1$ and $p>1$.

Regarding the multiplicity of solutions, again by the $\mathbb{Z}_2$--invariance of the functional, there exists infinitely many (sign-changing) solutions; however, unlike the sublinear case, this time we can find a sequence of solutions $(u_k)_k$ such that $\mathcal{I}(u_k)\to \infty$. A long standing open question is whether the symmetry of the functional is necessary to obtain multiplicity results; see the introduction of \cite{MollePassaseo,RamosTavaresZou} for a good overview. A least energy nodal solution, on the other hand, can be characterized by
\begin{equation}\label{eq:cnod_Dirichletsuper}
c_{nod}=\inf_{\mathcal{N}_{nod}} \mathcal{I},\quad \text{ where }\quad  \mathcal{N}_{nod}=\{u\in H^1_0(\Omega)\setminus \{0\},\ \mathcal{I}'(u)u^+=\mathcal{I}'(u)u^-=0\},
\end{equation}
see \cite{BartschWethWillem,CastroCossioNeuberger} (observe this set is not a $C^1$-manifold). On bounded radial domains, the associated solutions are not radial \cite{AftalionPacella}, but only foliated Schwarz symmetric \cite{BartschWethWillem}.

The study of the regularity of the zero-set of sign changing solutions is actually simpler in the superlinear case $p>1$ than in the sublinear one $p<1$  (although, in any case, is not at all simple); this is as a consequence of the map $f(t)=|t|^{p-1}t$ being of class $C^1$ for $p>1$ \cite{Hardtetal,Lin}.

\subsection{Dirichlet boundary conditions: the critical case $p=2^*-1=(N+2)/(N-2)$}\label{eq:critical_1eq}

 In this case, we are dealing with
\[
-\Delta u=|u|^{2^*-2}u \text{ in } \Omega,\qquad u=0 \text{ on } \partial \Omega,
\]
and the associated functional $\mathcal{I}$ does not satisfy the Palais-Smale condition for all levels $c$. The question of whether there are (nontrivial) solutions or not for $p=2^*-1$ or $p>2^*-1$ depends strongly on the domain.  When $\Omega$ is star-shaped, for instance, there are no solutions (by the Pohozaev identity, see for instance \cite[Theorem 3.4.26]{BadialeSerra}); however, there are examples of contractible domains where solutions do exist.  This shows that the topology of the domain is not enough to characterize the situation, although it has some influence: if, for some positive $d$, the homotopy group of $\Omega$ with $\mathbb{Z}_2$ coefficients is non trivial, $\mathcal{H}_d(\Omega,\mathbb{Z}_2)\neq \{0\}$, then we have a positive solution \cite{BahriCoron}. Multiplicity results are much more recent (and challenging); recent contributions are, for instance, \cite{ClappMussoPistoia,ClappWeth2004, MussoPistoiaJMPA2006, MussoPistoia2008}.

In order to emphasize how delicate the situation is in the critical case $p=2^*-1$, we make two remarks:
\begin{enumerate}
\item If the domain is not bounded but instead the whole $\R^N$, then we have (explicit!) solutions, the so called \emph{bubbles}:
\begin{equation}\label{eq:bubble}
U_{\delta,\xi}=(N(N-2))^{(N-2)/4}\frac{\delta^\frac{N-2}{2}}{(\delta^2+|x-\xi|^2)^\frac{N-2}{2}},\qquad \text{ for } \delta>0,\ \xi\in \R^N.
\end{equation}
\item If we consider a linear perturbation of the problem, namely:
\[
-\Delta u=\lambda u+|u|^{2^*-2}u \text{ in } \Omega,\qquad u=0 \text{ on } \partial \Omega
\]
 the situation changes. This problem has positive solutions for $\lambda\in (0,\lambda_1(\Omega))$ and $N\geq 4$ (the problem is commonly known as the Brezis-Nirenberg problem \cite{BrezisNirenberg}), and for $\lambda\in (\lambda^*(\Omega),\lambda_1(\Omega))$ in $N=3$, for some $\lambda^*(\Omega)>0$. The topology of the domain, in this situation, also influences multiplicity results: there exist at least $\text{cat}_\Omega(\Omega)$ solutions, where the (Lyusternik-Schnirelmann)  category of $\Omega$ is the least integer $d$ such that there exists a covering of $\Omega$ by $d$ closed contractible sets. As $\lambda\to 0$, the solutions tends to concentrate and blowup at certain points which depend on geometric properties of $\Omega$ \cite{HanAIHP,Rey1}.
\end{enumerate}
We recommend the survey \cite{PistoiaSurvey} for more results in the critical case. Therein, the reader can also find a nice and simple general explanation of the use of the Lyapunov-Schmidt reduction method as a powerful and useful technique to build solutions to semilinear elliptic problems.

\subsection{Lane-Emden equations with Neumann boundary conditions}\label{sec:Neumann_1eq}

Quite surprisingly, for the Neumann problem
\begin{equation}\label{eq:Neumann}
-\Delta u=|u|^{p-1}u \text{ in } \Omega,\qquad \partial_\nu u=0 \text{ on } \partial \Omega,
\end{equation}
very little was known before our work. Observe that solutions satisfy the compatibility condition
\[
\int_\Omega |u|^{p-1}u=0,
\]
hence all nontrivial solutions necessarily change sign. Therefore, \emph{least energy solutions} are actually least energy \emph{nodal} solutions.\footnote{The situation changes drastically if we consider, instead, the problem $-\Delta u+\lambda u=|u|^{p-1}u$ with $\lambda>0$, which allows positive solutions. This has been extensively studied since the celebrated papers \cite{AdimurthiMancini,LinNiTakagi,Wang1991}. Since the results are different in nature, we do not make a literature review of this case.}

For the subcritical case $p<1$, the existence of least energy (nodal) solutions was established in \cite{PariniWeth2015}. When $\Omega$ is a ball, the authors proved that these solutions are \emph{not} radial but only foliated Schwarz symmetric (axially symmetric and decreasing as a function of the polar angle).

For the critical exponent case $p=2^*-1$, differently to the Dirichlet case (recall Subsection \ref{eq:critical_1eq}), in the Neumann one there are solutions, as was shown in \cite{CK91} using  a dual variational formulation.

Combining this with two of my works, one  with A. Salda\~na for the subcritical case \cite{SaldanaTavares}, and the other with A. Pistoia and D. Schiera for the critical one \cite{PistoiaSchieraTavares}, we have the following:

\begin{theorem}[Combination of \cite{PariniWeth2015,CK91} with my recent papers
\cite{PistoiaSchieraTavares,SaldanaTavares}]\label{thm:singleeq_neumann}\ \newline
Let $\Omega$ be a bounded smooth domain in $\R^N$, $N\geq 1$. Assume that  $p<2^*-1$ and $p\neq 1$ (with the convention that $2^*=\infty$ for $N=1,2$), or $p=2^*-1$ and $N\geq 4$. Then there exist least energy solutions, and they are all classical solutions. Moreover:
\begin{enumerate}
\item (Case $N=1$ and $\Omega=(-1,1)$): every least energy solution is {strictly monotone} in $\Omega$.
\item (Case $N\geq 2$ and  $\Omega$ a ball or an annulus):
\begin{enumerate}
\item every least energy solution is {foliated Schwarz symmetric} and it is \underline{not} radially symmetric.
\item there exist least energy radial solutions; they are classical and  strictly monotone in the radial variable.
\end{enumerate}
\end{enumerate}
\end{theorem}

The symmetry breaking is done by contradiction, and this is why it is so important to prove the existence and monotonicity of least energy \emph{radial} solutions. These are defined, when $\Omega$ is symmetric, as solutions of the problem \eqref{eq:Neumann} which achieve the level
\[
c_{rad}:=\inf\{\mathcal{I}(u):\ u\in H^1(\Omega) \text{ is a nontrivial radial solution of } \eqref{eq:Neumann}\}.
\]A key element in the proof of the monotonicity  is the use of a dual variational formulation combined with a new $L^t$-norm-preserving transformation introduced in \cite{SaldanaTavares}, which combines a suitable flipping with a decreasing rearrangement. This combination allows us to treat annular domains, sign-changing functions, and Neumann problems, which are non-standard settings to use rearrangements and symmetrizations. Both \cite{PistoiaSchieraTavares} and \cite{SaldanaTavares} prove the results of Theorem \ref{thm:singleeq_neumann} for the more general context of Lane-Emden systems, and the single equation case follows as corollary. We will discuss this in more detail in Section \ref{sec3} below.

\smallbreak

We point out two other things:
\begin{itemize}
\item in a recent work with A. Salda\~na \cite{SaldanaTavaresNoDEA} we showed the convergence of least energy nodal solutions in terms of $p$; in particular, the limit as  $p\to 1$ depends on the domain.
\item jointly with M. Grossi and A. Salda\~na \cite{GrossiSaldanaTavares}, we deduced the blowup behavior as $p\nearrow 2^*-1$ of all radial solutions  of \eqref{eq:Neumann}; incidently, in order to prove it we had to prove at the same time the behaviour of all radial Dirichlet solutions \eqref{eq:LE_example}, generalizing the work \cite{HanAIHP}.
\end{itemize}

Another interesting open question is whether or not one has a solution for all $p>2^*-1$. In a joint work with A. Pistoia and A. Salda\~na \cite{PistoiaSaldanaTavares}, using the Lyapunov-Schmidt reduction method, we proved the existence of solutions in the slightly supercritical case, when the domain has some symmetries. For instance, if $\Omega$ is the ball, the main result therein is the following:
\begin{theorem}[{{\cite{PistoiaSaldanaTavares}}}]\label{thm:main:ball} Take $N\geq 4$, let $\Omega\subset \R^N$  be the unit ball centered at the origin, and $p=2^*-1+\varepsilon$.  There exists $\varepsilon_0>0$ such that, for $\varepsilon\in (0,\varepsilon_0)$, the problem \eqref{eq:Neumann} has a solution $u_\varepsilon$ which is even in $x_1,\ldots, x_{N-1}$ and odd in $x_N$; this solution ``looks like'' the difference of two bubbles \eqref{eq:bubble}, concentrating at antipodal points as $\varepsilon\to 0$.
\end{theorem}
This result has been recently extended to Lane-Emden systems, see \cite{GuoPeng}.

\section{Elliptic Hamiltonian systems}\label{sec3}

A natural extension of studying the single equation  $-\Delta u=|u|^{p-1}u$ is to deal with the following particular example of an \emph{elliptic  Hamiltonian system}\footnote{The name has its origins in the following: \eqref{eq:Hamiltonian} has the shape $-\Delta u=-H_v(u,v)$, $-\Delta v=H_u(u,v)$ for $H(u,v)=|u|^{p+1}/(p+1)-|v|^{q+1}/(q+1)$.}, also known in the literature as \emph{Lane-Emden system}:
\begin{equation}\label{eq:Hamiltonian}
-\Delta u=|v|^{q-1}v,\qquad -\Delta v=|u|^{p-1}u \qquad \text{ in } \Omega.
\end{equation}
One considers $p,q>0$ in either the \emph{sublinear} ($pq<1$) or the \emph{superlinear} ($pq>1$) cases. The correct notion of criticality correspond to $(p,q)$ being on the \emph{critical hyperbola}
\begin{equation}\label{eq:critical}
\frac{1}{p+1}+\frac{1}{q+1}=\frac{N-2}{N}.
\end{equation}
(see \cite{ClementdeFigueiredoMitidieri, PeletiervanderVorst}). In this setting, the more general notion of ``linearity'' is $pq=1$ or, equivalently, $1/(p+1)+1/(q+1)=1$, see \cite{ClementvanderVorst}; being subcritical means that:
\begin{align}\label{sc:intro}
p,q>0,\qquad pq\neq 1,\qquad \text{ and }\qquad \frac{1}{p+1}+\frac{1}{q+1}>\frac{N-2}{N},
\end{align}
and this condition is trivially satisfied if $N=1,2$ or if $pq<1$.

A simple way of motivating these notions is to formally write $v:=-|\Delta u|^{\frac{1}{q}-1}\Delta u$, so that \eqref{eq:Hamiltonian} reduces to the higher order problem
\begin{equation}\label{eq:higher_order}
\Delta \left(|\Delta u|^{\frac{1}{q}-1}\Delta u\right)=|u|^{p-1}u \qquad \text{ in } \Omega.
\end{equation}
This equation has the following associated Euler-Lagrange functional
\[
u\mapsto \frac{q}{q+1}\int_\Omega |\Delta u|^\frac{q+1}{q}-\frac{1}{p+1}\int_\Omega |u|^{p+1}.
\]
So, the notion of criticality \eqref{eq:critical} is related with the validity of the embedding $W^{2,\frac{q+1}{q}}(\Omega) \hookrightarrow L^{p+1}(\Omega)$, while ``linear'' (with an abuse of language) is related to having the same homogeneity on the equation \eqref{eq:higher_order}. This variational formulation can be made precise, see for instance \cite[Section 4]{BonheureSantosTavares} for the case of Dirichlet boundary conditions.
\begin{figure}
    \centering
    \includegraphics[scale=0.35]{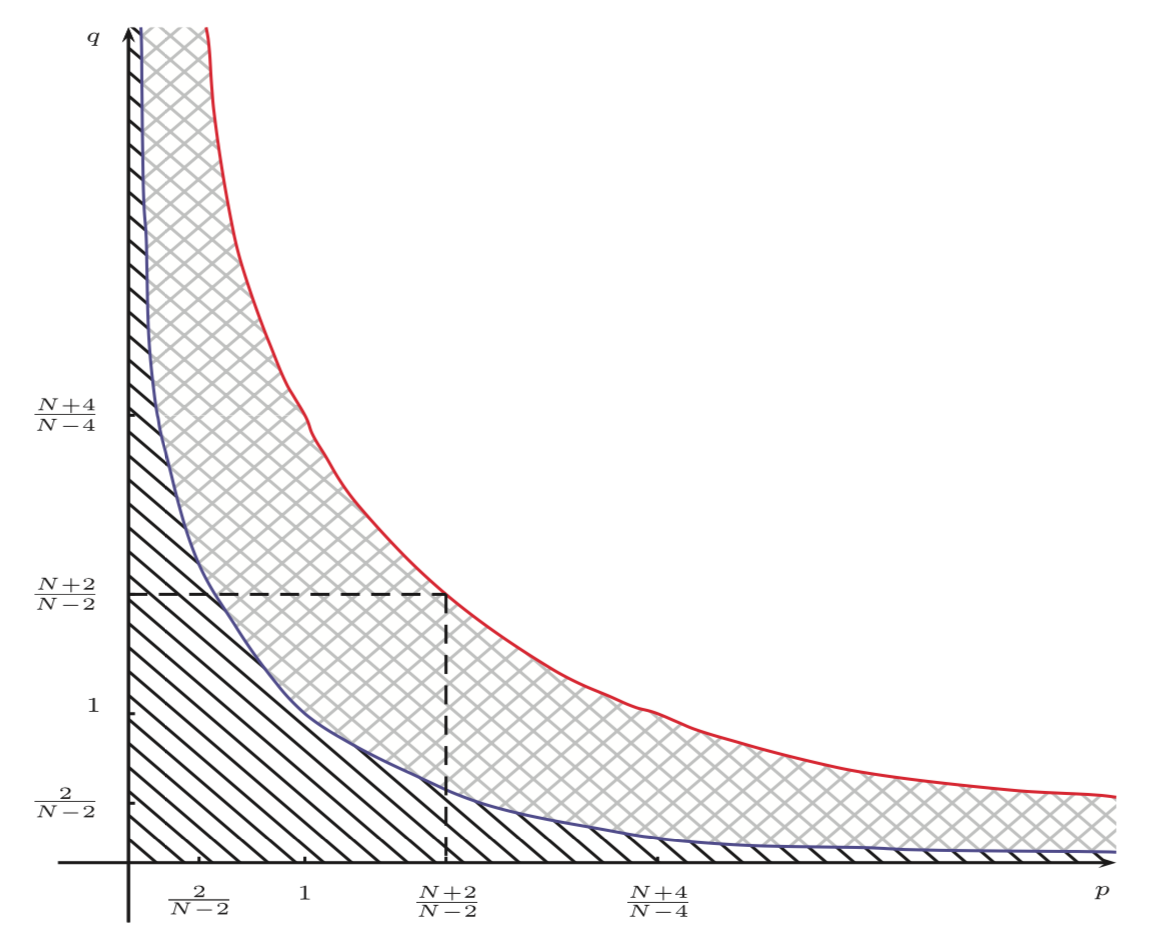}
    \caption{The curve on top is the \emph{critical hyperbola}. The hyperbola below corresponds to $pq=1$. Image taken from \cite{BonheureSantosTavares}.}
\end{figure}

Systems with a Hamiltonian structure such as \eqref{eq:Hamiltonian} have been extensively studied in the past 25 years, and many results are known regarding existence, multiplicity, concentration phenomena, positivity, symmetry, Liouville theorems, \emph{etc}.  We refer to the surveys \cite{BonheureSantosTavares, deFigueiredo, Ruf} for an overview of the topic. In the Dirichlet case, for $(p,q)$ below the critical hyperbola and not belonging to the linear one, existence, multiplicity and symmetry results for least energy solutions have been shown using several variational approaches, each one with its pros and cons. Together with D. Bonheure and E. Moreira dos Santos, I have written a  long and detailed survey \cite{BonheureSantosTavares} exploring the advantages and disadvantages of each approach, defining rigorously the notion of least energy solution, and providing several proofs.  As for the critical or supercritical case (still with Dirichlet boundary conditions), a Pohozaev-type identity also rules out the existence of nontrivial solutions in star shaped domains for systems. When $(p,q)$ approaches a point on the critical hyperbola, we are only aware of concentration and blowup results in the paper \cite{Guerra} where, however, there is the technical restriction of considering either $p$ or $q$ as being fixed. See also \cite{Guerra2} for the case when $(p,q)$ approaches asymptotically at infinity the critical hyperbola. Some relations to an 1-biharmonic equation have been shown in \cite{AbatangeloSaldanaTavares}, when either $p$ or $q$ go to infinity.

Recalling what is known for the single equation case (recall Section \ref{sec2} above), one is tempted to ask about least energy nodal solutions in the Dirichlet case, and what happens in the case of Neumann boundary conditions. These have been, indeed, my contributions to the field, and I describe them below.

\subsection{Dirichlet boundary conditions: least energy nodal solutions}

Condition \eqref{eq:critical} or \eqref{sc:intro}  (resp. the critical and subcritical cases) together with Sobolev embeddings and the Rellich--Kondrachov theorem imply the following embeddings
\begin{align}\label{embed}
W^{2,\frac{p+1}{p}}(\Omega)\hookrightarrow L^{q+1}(\Omega)\  \text{ and }\  W^{2,\frac{q+1}{q}}(\Omega)\hookrightarrow L^{p+1}(\Omega),
\end{align}
which are compact in the subcritical case. A strong solution of
\begin{equation}\label{NHS_Dirichlet}
-\Delta u=|v|^{q-1}v,\qquad -\Delta v=|u|^{p-1}u \text{ in } \Omega \qquad u=v=0 \text{ on } \partial \Omega,
\end{equation}
 is defined as a pair $(u,v)\in (W^{2,\frac{q+1}{q}}(\Omega)\cap W^{1,\frac{q+1}{q}}_0(\Omega)) \times (W^{2,\frac{p+1}{p}}(\Omega) \cap W^{1,\frac{p+1}{p}}_0(\Omega))$ satisfying the equations a.e. in $\Omega$, and the boundary conditions in the trace sense. The system is strongly coupled, in the sense that $u\equiv 0$ if, and only if, $v\equiv 0$, or $u$ is sign-changing if, and only if, $v$ is sign-changing.
Problem \eqref{NHS_Dirichlet} has a variational structure, and \eqref{NHS_Dirichlet} are the Euler-Lagrange equations of the energy functional
\begin{equation}\label{eq:functional_I}
(u,v)\mapsto I(u,v)= \int_\Omega \nabla u\cdot \nabla v-\frac{|u|^{p+1}}{p+1}-\frac{|v|^{q+1}}{q+1}\, dx.
\end{equation}
We define a \emph{least energy solution} as a nontrivial strong solution of \eqref{NHS_Dirichlet} achieving the level
\begin{align}\label{c:level}
c:=\inf\left\{ I(u,v):\ (u,v)\not\equiv (0,0),\ (u,v) \text{ is a strong solution of \eqref{NHS_Dirichlet}}\right\}.
\end{align}
In view of \eqref{sc:intro} and \eqref{embed}, the functional $I$ is well defined at strong solutions. As we said before, existence of least energy solutions is established via several different approaches in the subcritical case \eqref{sc:intro} (see \cite{BonheureSantosTavares}). Here we are interested in the least energy nodal solutions, \emph{i.e.}, strong solution of \eqref{NHS_Dirichlet} achieving the level
\begin{align}\label{c:level}
c_{nod}:=\inf\left\{ I(u,v):\ \ (u,v) \text{ is a strong solution of \eqref{NHS_Dirichlet}},\ u^\pm \not \equiv 0,\ v^\pm \not \equiv 0\right\}.
\end{align}
Actually we proved, together with D. Bonheure, E. Moreira dos Santos and M. Ramos \cite{BMRT15}, the existence and partial symmetry of least energy nodal solutions for the more general problem of H\'enon--type:
\begin{equation}\label{Henon}
-\Delta u =|x|^\beta |v|^{q-1}v,\quad -\Delta v=|x|^\alpha |u|^{p-1}u \text{ in } \Omega,\qquad u=v=0 \text{ on } \partial \Omega.
\end{equation}
in the superlinear--subcritical case.  It is not obvious that the level $c_{nod}$ is achieved, since this no longer follows from a simple minimization argument. Indeed, even if we have enough compactness to extract a converging subsequence, the limit could be a critical point $(u,v)$ such that both $u$ and $v$ are positive (or negative). The existence of a least energy nodal solution for the scalar Lane-Emden equation \cite{BartschWethWillem,CastroCossioNeuberger} follows from the minimization of the functional over a nodal Nehari set (recall \eqref{eq:cnod_Dirichletsuper}).  However, it is not clear at all how such a nodal Nehari set associated with the energy functional
\[
(u,v)\mapsto \int_\Omega \nabla u\cdot \nabla v-|x|^\alpha \frac{|u|^{p+1}}{p+1}-|x|^\beta\frac{|v|^{q+1}}{q+1}\, dx.
\] could be defined. We follow a dual variational framework and polarization techniques, proving the following.

\begin{theorem}[{{\cite{BMRT15}}}]\label{thm:main1}
Let $N\geq 1$, $\alpha \geq 0$, $\beta\geq 0$ and suppose that $(p,q)$ is superlinear and subcritical. Then there exists a least energy nodal solution of \eqref{Henon}.

Moreover, when $N \geq 2$ and  $\Omega$ is either a ball or an annulus centered at the origin, then every  least energy nodal solution $(u,v)$ is foliated Schwarz symmetric with to respect to some direction $e\in \partial B_1(0)$.
\end{theorem}

For the Lane-Emden equation, symmetry breaking in radial domains is proved via a Morse index argument \cite{AftalionPacella}. For the H\'enon-Lane-Emden system \eqref{Henon}, it is not clear how to compute (or even define) the Morse index of the solutions. Nevertheless, using a perturbation argument, in \cite{BMRT15} we prove symmetry breaking (\emph{i.e.}, least energy nodal solutions are not radial) for some ranges of the parameters, namely when $\alpha\sim \beta$ and $p\sim q$. We also observe that our results contain, as a particular case, the following one for the biharmonic problem, complementing some results from \cite{WethTMNA2006}.

\begin{corollary}[{{\cite{BMRT15}}}]\label{th:biharmonic}
Let $\Omega \subset \R^N$, $N \geq 1$ and assume that $\frac{1}{2}> \frac{1}{p+1}>\frac{N-4}{2N}$. Then the fourth order problem
\[
\Delta^2 u= |x|^{\alpha}|u|^{p-1}u \quad \text{ in } \Omega, \qquad u=\Delta u=0\quad \text{ on } \partial \Omega
\]
admits a least energy nodal solution. Moreover, if $\Omega$  is either a ball or an annulus centered at the origin, $N\geq 2$, then any least energy nodal solution is such that $u$ and $-\Delta u$ are foliated Schwarz symmetric with respect to the the same unit vector $e\in \R^N$.
\end{corollary}

\subsection{Neumann boundary conditions: least energy solutions in the subcritical and critical cases}

The papers mentioned before in this section work with Dirichlet boundary conditions and, up to our knowledge, the few papers addressing Neumann problems are
\cite{AvilaYang, PistoiaRamosNeumann,RamosYang,  Zeng}, where existence of positive solutions and concentration phenomena are studied, and \cite{BonheureSerraTilli}, which focuses on existence of positive radial solutions. However, these papers focus on a \emph{different} operator of the form $Lw=-\Delta w+V(x) w$, with $V$ positive and bounded.  In comparison with problem
\begin{equation}\label{NHS}
-\Delta u=|v|^{q-1}v,\qquad -\Delta v=|u|^{p-1}u \text{ in } \Omega \qquad u_\nu=v_\nu=0 \text{ on } \partial \Omega,
\end{equation}
the shape of solutions changes drastically; for instance, the operator $L$ with Neumann boundary conditions induces the $H^1$-norm
\[
w\mapsto \int_\Omega (|\nabla w|^2+V(x)w^2),
\] and this allows the existence of positive solutions, while all nontrivial solutions of \eqref{NHS} are \emph{sign-changing}. Indeed, if $(u,v)$ is a classical solution of \eqref{NHS}, then by the Neumann boundary conditions and the divergence theorem,
\begin{align}\label{comp}
\int_\Omega |u|^{p-1}u=\int_\Omega |v|^{q-1}v=0.
\end{align}
Since $u\equiv 0 $ if and only if $v\equiv 0$, \eqref{comp} is only satisfied if $(u,v)$ is trivial or if both components are sign-changing (recall also what happens in the single equation case, Section \ref{sec:Neumann_1eq})
  As far as we know, we were the first to study problem \eqref{NHS}.

\paragraph{The subcritical case.}   The study of the Neumann problem was initiated recently in a joint paper with A. Salda\~na \cite{SaldanaTavares}, where we prove that, in the subcritical case, least energy (nodal) solutions exist and, whenever $\Omega$ is a radial domain, they are not symmetric but only foliated Schwarz symmetric.

For the Neumann problem, a strong solution of \eqref{NHS} is defined as a pair $(u,v)\in W^{2,\frac{q+1}{q}}(\Omega)\times W^{2,\frac{p+1}{p}}(\Omega)$ satisfying the equations a.e. in $\Omega$, and the boundary conditions in the trace sense.
Least energy solutions can be defined exactly as in the previous section.

\begin{theorem}[{{\cite{SaldanaTavares}}}]\label{thm1:SaldanaTavares}
Consider $(p,q)$ in the subcritical regime, \emph{i.e.},
\[
\frac{1}{p+1}+\frac{1}{q+1}>\frac{N-2}{N},\qquad pq\neq 1.
\]
The least energy level is achieved and, if $(u,v)$ is a least energy solution, then:
\begin{itemize}
\item it is a classical solution.
\item (Monotonicity) If $N=1$ and $\Omega=(-1,1)$, then $u'v'>0$ in $\Omega$;  in particular, $u$ and $v$ are both strictly monotone increasing or both strictly monotone decreasing in $\Omega$.
\item (Partial symmetry \& symmetry breaking) If $N\geq 2$ and $\Omega$ is either a ball or an annulus, then $u$ and $v$ are foliated Schwarz symmetric with respect to the same vector. Moreover, $u$ and $v$ are \emph{not}  radially symmetric.
\end{itemize}
\end{theorem}

In particular, the case $p=q$ leads to the results in the subcritical case in Theorem \ref{thm:singleeq_neumann}. The approach to show this theorem is based on a variant of the \emph{dual method} \cite{AlvesSoares, ClementvanderVorst}. Before describing rigorously the dual framework, we formally observe that
\[
-\Delta u=|v|^{q-1}v,\  -\Delta v=|u|^{p-1}u \iff u=(-\Delta )^{-1}(|v|^{q-1}v),\   v=(-\Delta )^{-1}(|u|^{p-1}u);
\]
by introducing the new dual variables $f=|u|^{p-1}u$, $g=|v|^{q-1}v$, we obtain
\[
|f|^{\frac{1}{p}-1}f=(-\Delta )^{-1}g,\ |g|^{\frac{1}{q}-1}f=(-\Delta )^{-1}f.
\]
In the definition of $(-\Delta)^{-1}$, and recalling that we are dealing with Neumann boundary conditions, one needs to take into account the normalization \eqref{comp}. To rigorously perform these steps, we introduce some notation. Let $p$ and $q$ satisfy \eqref{sc:intro} and, for $s>1$, let
\begin{align}
X^s=\Big\{f\in L^{s}(\Omega): \ \int_\Omega f =0\Big\},\qquad
 \text{and}\qquad
X:=X^\frac{p+1}{p}\times X^\frac{q+1}{q},
\end{align}
endowed with the norm $\|(f,g)\|_X=\|f\|_\frac{p+1}{p}+\|g\|_\frac{q+1}{q}$. Let $K$ denote the inverse (Neumann) Laplace operator with zero average, that is, if $h\in X^s(\Omega)$, then $u := Kh\in W^{2,s}(\Omega)$ is the unique strong solution of $-\Delta u = h$ in $\Omega$ satisfying $\partial_\nu u=0$ on $\partial \Omega$ and $\int_\Omega u = 0$.  In this setting, the (dual) energy functional $\phi:X\to \R$ is given by
\begin{align}\label{phi:def}
 \phi(f,g):=\frac{p}{p+1}\int_\Omega |f|^{\frac{p+1}{p}}+ \frac{q}{q+1}\int_\Omega |g|^{\frac{q+1}{q} }-\int_\Omega g\, K f,\qquad (f,g)\in X.
\end{align}
Since the conditions $\int_\Omega f=\int_\Omega g=0$ are included in $X$, in order to associate an equation to a critical point of $\phi$, we require a suitable translation of $K$ (which is related with \eqref{comp}). For $t>0$, let $K_t:X^\frac{t+1}{t}\to W^{2,\frac{t+1}{t}}(\Omega)$ be given by
\begin{align}\label{Ks:def}
K_t h:=Kh +c_t(h)\qquad \text{ for some }\quad c_t(h)\in\R\quad \text{ such that }\begin{cases}-\Delta (K_t h)=h \text{ in } \Omega\\
\partial_\nu (K_t h)=0 \text{ on } \partial \Omega,\\
\int_\Omega |K_t h|^{t-1}K_t h=0.
\end{cases}
\end{align}
Then, a critical point $(f,g)$ of $\phi$ solves the dual system $K_{q} f=|g|^{\frac{1}{q}-1}g$ and $K_{p} g= |f|^{\frac{1}{p}-1}f$ in $\Omega$. This is the starting point for the existence part in Theorem \ref{thm1:SaldanaTavares}. The symmetry breaking result, on the other hand, is based on a contradiction argument which, in turn, follows once we deduce the monotonicity of least energy radial solutions. The  proof of  the latter is based on a new $L^t$--norm-preservation transformation introduced in \cite{SaldanaTavares}, which we now recall.

For $\Omega=B_{R}(0)\setminus \overline B_{r}(0)$ an annulus or $\Omega=B_R(0)$ a ball (in which case we define $r:=0$), let
\begin{align*}
&\mathcal{I}:L_{rad}^\infty(\Omega)\to C_{rad}(\overline{\Omega}),\qquad {\cal I}h(x):=\int_{\{r \leq|y|{\leq}|x|\}} h(y)\ dy
=N\omega_N\int_{r}^{|x|}h(\rho)\rho^{N-1}\ d\rho\\
&\mathfrak{F}: C_{rad}(\overline{\Omega})\to L_{rad}^\infty(\Omega),\qquad \mathfrak{F} h:=(\chi_{\{{\cal I}h>0\}}-\chi_{\{{\cal I}h\leq 0\}})\, h.
\end{align*}
\begin{definition}[{{\cite{SaldanaTavares}}}]
For $h\in C_{rad}(\overline{\Omega})$, the $\divideontimes$-transformation is given by
\begin{align*}
h^{\divideontimes}\in L^\infty_{{rad}}(\Omega),\qquad h^{\divideontimes}(x) := ({\mathfrak F} h)^\#(\omega_N |x|^N-\omega_N {r}^N),
\end{align*}
where $\omega_N= |B_1|$ is the volume of the unitary ball in $\R^N$ and $\#$ is the decreasing rearrangement given by
\begin{align*}
h^\#:[0,|\Omega|]\to\R,\qquad h^{\#}(0):=\text{ess sup}_{\Omega} h\quad h^\#(s):=\inf\{t\in\R\::\: |\{h>t\}|<s\},\ s>0.
\end{align*}
\end{definition}

\begin{figure}[h]
\centering
\begin{subfigure}{.45\textwidth}
  \begin{picture}(150,170)
    \put(0,5){\includegraphics[width=.90\textwidth]{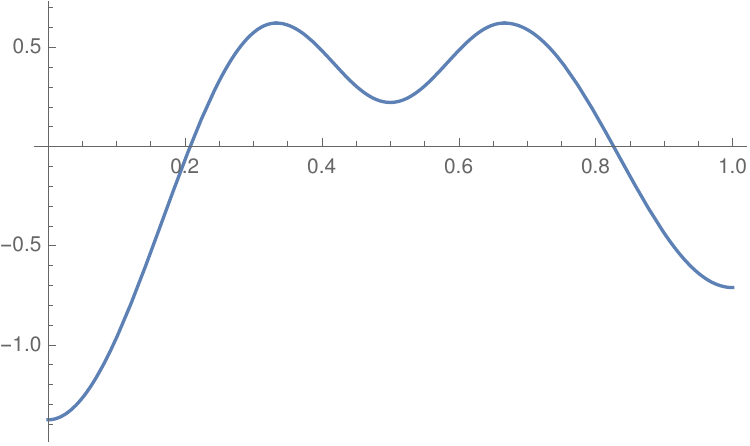}}
    \put(15,115){$h(r)$}
    \put(185,85){$r$}
  \end{picture}
\end{subfigure}%
\begin{subfigure}{.45\textwidth}
\begin{picture}(150,170)
    \put(0,5){\includegraphics[width=.90\textwidth]{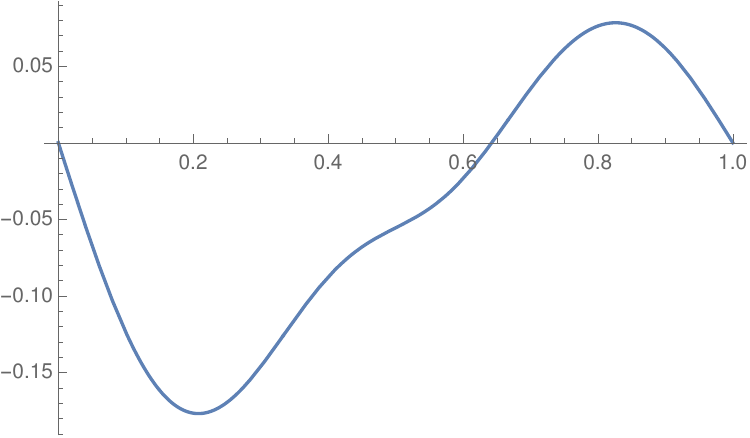}}
    \put(15,115){$\mathcal{I} h(r)$}
    \put(185,85){$r$}
  \end{picture}
\end{subfigure}
\begin{subfigure}{.45\textwidth}
\begin{picture}(150,170)
    \put(0,5){\includegraphics[width=.90\textwidth]{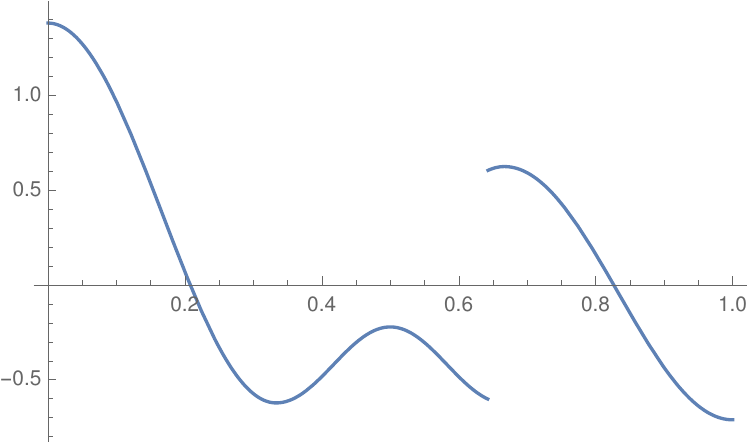}}
    \put(15,115){$\mathcal{F} h(r)$}
    \put(185,50){$r$}
  \end{picture}
\end{subfigure}%
\begin{subfigure}{.45\textwidth}
\begin{picture}(150,170)
  \put(0,5){\includegraphics[width=.90\textwidth]{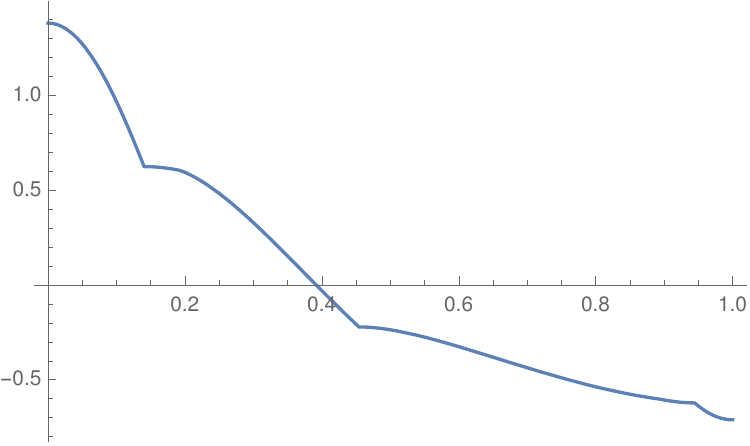}}
    \put(15,115){$h^\divideontimes(r)$}
    \put(185,50){$r$}
  \end{picture}
\end{subfigure}
\caption{Examples of the functions $\mathcal{I} h$, $\mathfrak{F} h$ and  $h^\divideontimes$  for a particular radial function $h\in C(\overline{B_1(0)})$. Images taken from \cite{SaldanaTavares}.
}
\label{fig}
\end{figure}
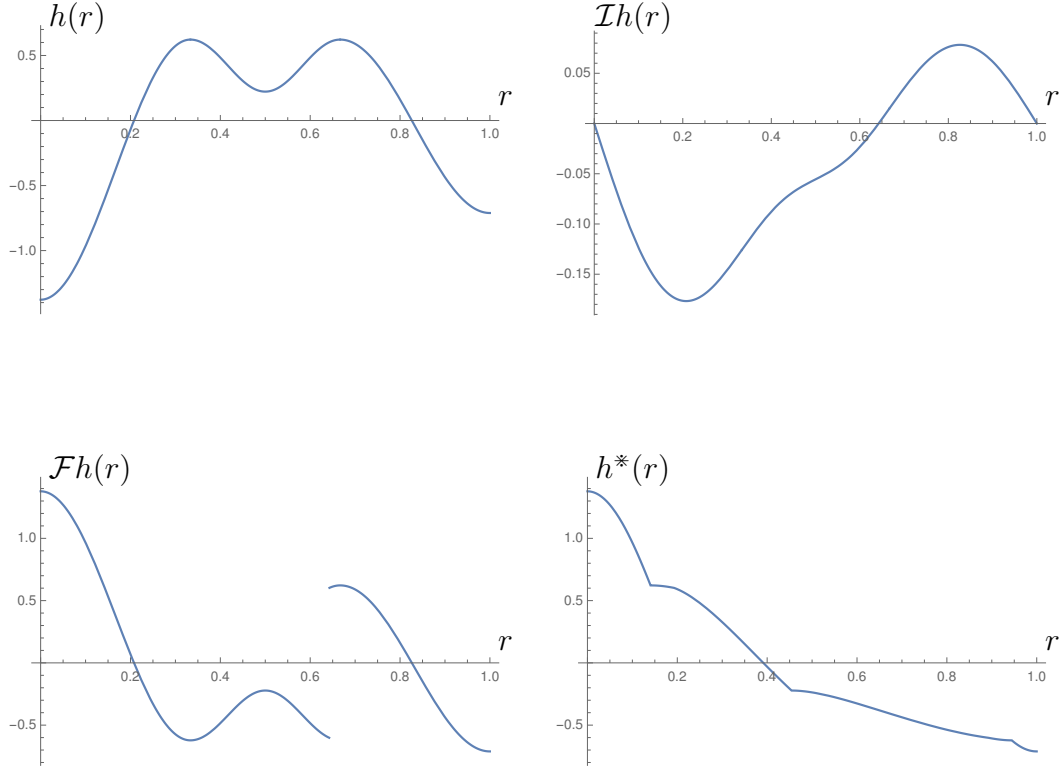

This transformation, in practice, is applied to radial functions $h$ with zero average: if $\Omega=B_{R}(0)\setminus \overline B_{r}(0)$, then $\mathcal{I}h(R)=\int_\Omega h=0$. Loosely speaking, this transformation does the following: in many situations, the domain of $h$ may be split in $r=:r_0<r_1<\ldots<r_N=R$, where $\int_{r_i<|x|<r_{i+1}}h=0$; however, it may not be true that $\mathcal{I}h(x)$ is nonnegative for every $x$. We flip the graph of $h$ in the annuli $[r_i,r_{i+1}]$ where we do not have this property, obtaining at the end $\mathcal{I}(\mathfrak{F} h)(x)\geq 0$ for every $x\in \Omega$. We then apply a decreasing rearrangement, finally placing the result back in the original domain using the transformation $x\mapsto \omega_N|x|^N-\omega_N r^N$. For more details, insights, examples and comments regarding the definition of the  \emph{flip-\&-rearrange transformation} $\divideontimes$ we refer to \cite[Section 3.2]{SaldanaTavares}, see also Figure \ref{fig} for an example. The following is a combination of Theorem 1.3 and Proposition 3.4 from \cite{SaldanaTavares}.

\begin{theorem}\label{thm:trans} Let $p,q>0$ and take $\Omega$ to be a ball or an annulus centered at the origin. Take  $f,g:\overline{\Omega}\to \R$ be continuous and radially symmetric functions with $\int_\Omega f = \int_\Omega g= 0$. Then $(f^\divideontimes,g^\divideontimes)\in X_{rad}$,
\begin{align}\label{in}
\|f^\divideontimes\|_\frac{p+1}{p}=\|f\|_\frac{p+1}{p}\quad \|g^\divideontimes\|_\frac{q+1}{q}=\|g\|_\frac{q+1}{q} \text{ and }\quad  \int_\Omega fK g\leq \int_\Omega f^\divideontimes K g^\divideontimes.
\end{align}
Furthermore, if $f,g$ are nontrivial and the last statement in \eqref{in} holds with equality, then $f,g$ are monotone in the radial variable. Moreover, $(Kf,Kg)$ is radially symmetric and $(Kf)_r(Kg)_r>0$.
\end{theorem}

Observe that annuli, sign-changing functions, and Neumann boundary data are non-standard conditions to work with rearrangements. We take advantage of working with \emph{radial} functions and of the fact that we use \emph{Lebesgue spaces} (and not Sobolev ones) within a dual framework; this gives more flexibility in the construction of our transformation although, on the other hand, it becomes harder to control the nonlocal operator $K$. In this sense, instead of rearrange $u,v$ directly, we are transforming the \emph{dual variables} $f,g$ to obtain, together with variational techniques,  information about the monotonicity of solutions.

\paragraph{The critical case.} In a joint work with A. Pistoia and D. Schiera \cite{PistoiaSchieraTavares}, we extended some of the results of the previous paragraph to the critical case. In this direction, the main result of our paper is the following:
\begin{theorem}[{{\cite{PistoiaSchieraTavares}}}]\label{main thm}
Let $p, q$ satisfy \eqref{eq:critical}, and moreover
\[
N \ge 6 \text{ and }
p,q > \frac{N+2}{2(N-2)}, \quad
\text{ or }\quad
N=5\text{ and } p, q> \frac{17}{13},\quad \text{ or }\quad
N=4\text{ and }p, q > \frac 7 3.
\]
 Then there exists a least energy (nodal) solution of \eqref{NHS}, which is a classical solution.
\end{theorem}

Since the problem is critical, the embeddings \eqref{embed} are not compact, and in general the dual functional does not satisfy the Palais-Smale condition. We prove, however, a compactness condition, which is based on a new class of Cherrier--type inequalities: for every $\varepsilon>0$ there exists $C(\varepsilon)>0$ such that
\begin{equation*}
 \|u\|_{\frac{N\eta}{N-2\eta}} \le  \left (\frac{2^{\frac2N}}{S} +\varepsilon\right ) \|\Delta u\|_{\eta} + C(\varepsilon) \|u\|_{W^{1,\eta}}, \quad \forall u \in W^{2, \eta}_\nu(\Omega),
 \end{equation*}
where $W^{2, \eta}_\nu(\Omega):=\{u\in W^{2, \eta}(\Omega):\ \partial_\nu u=0 \text{ on } \partial \Omega\}$. Here, we are inspired by \cite{BonheureCheikhNascimento}, in which the case $\eta=2$ is shown. We would like to observe that, exactly as in \cite{CK91, SaldanaTavares}, we use the dual method.

\smallbreak

If we take $N \ge 5$, $p=1$, $q = \frac{N+4}{N-4}$, or equivalently, $q=1$, $p = \frac{N+4}{N-4}$, system \eqref{NHS} reduces to the fourth order problem
\begin{equation}\label{biharmonic}
\Delta^2 u = |u|^{\frac{8}{N-4}} u \text{ in } \Omega,\qquad
u_\nu=(\Delta u)_\nu= 0  \text{ on } \partial \Omega.
 \end{equation}
If $N > 6$, then
\[  \frac{N+2}{2(N-2)} < \frac{N+4}{N-4}, \quad \text{ and } \quad \frac{N+2}{2(N-2)} < 1, \]
hence the study of \eqref{biharmonic}
is contained in Theorem \ref{main thm} if $N >6$.
However, the case $N=5, 6$ and $p=1$ is not included. However, we prove directly that, for $N\geq 5$, there exists a least energy (nodal) solution to problem \eqref{biharmonic}.
As a consequence, via a perturbation argument we show the following for systems.
\begin{theorem}[{{\cite{PistoiaSchieraTavares}}}]\label{main thm1.2}
Let $N=5,6$. There exists $\varepsilon=\varepsilon(N,\Omega)$ such that, if $p,q$  satisfy \eqref{eq:critical} and either
\[
|p-1|+\left|q-\frac{N+4}{N-4}\right|<\varepsilon \quad \text{ or } \quad \left|p-\frac{N+4}{N-4}\right|+|q-1|<\varepsilon,
\]
then there exists a least energy (nodal) solution of \eqref{NHS}, which is a classical solution.
\end{theorem}

We point out that the fact that least energy solutions are classical solutions is a consequence of the following result, which is new.

\begin{proposition}[{{\cite[Proposition 1.3]{PistoiaSchieraTavares}}}] Let $(u, v)$ be a strong solution to \eqref{NHS}, where $p, q$ satisfy \eqref{eq:critical}. Then $(u, v) \in C^{2, \zeta}(\bar{\Omega}) \times C^{2, \eta}(\bar{\Omega})$, with: $\zeta<q$ if $0<q<1$, and $\zeta \in(0,1)$ if $q \geq 1 ; \eta<p$ if $0<p<1$, and $\eta \in(0,1)$ if $p \geq 1$.
\end{proposition}

We observe that, unlike what happens in the single equation case, a Brezis-Kato type argument does not seem to work to prove this regularity result. Instead, we rely on the bound $|G(x,y)|\leq C/|x-y|^{N-2}$ for the (Neumann Laplacian's) Green function, together with the Hardy-Littlewood-Sobolev inequality.

In our paper, for $p,q$ satisfying the conditions of Theorem \ref{main thm} or Theorem \ref{main thm1.2}, and when $\Omega$ is a ball or an annulus, we also prove that least energy solutions are foliated Schwarz symmetric with respect to the same vector and are not radial. Moreover, when $\Omega$ is an annulus, least energy radial solutions exist; when $\Omega$ is a ball this is in general an open problem, unless $p=q=2^*$, where the answer is negative - see \cite[Remark 6.5]{PistoiaSchieraTavares} for more details.
\section{Existence of fully nontrivial solutions to a class of gradient elliptic systems}\label{sec4}
Consider the following system with $d\geq 2$ equations
\begin{equation}\label{S-system}
\begin{cases}
-\Delta u_{i}+\lambda_{i}u_{i}=u_{i}|u_i|^{p-2}\sum_{j = 1}^{d}\beta_{ij}  |u_{j}|^p   ~\text{ in } \Omega,\\
u_{i}=0 \text{ on } \partial\Omega,  \quad i=1,...d,
\end{cases}
\end{equation}
where $\Omega$ is a domain of $\R^N$, $N\geq 1$, $\lambda_i \in \R$, in a (Sobolev) critical or subcritical regime $0<p\leq2^*/2=N/(N-2)$ if $N\geq 3$, or $0<p<+\infty$ for $N=1,2$. We assume from now on that $\beta_{ij}=\beta_{ji}$ for $i\neq j$ and so, from a mathematical point of view, this is an example of a weakly coupled elliptic system with gradient terms\footnote{Indeed, it has the form $-\Delta u_i+\lambda_i u_i=H_{u_i}(u_1,\ldots, u_m)$, with $H(u_1,\ldots, u_m)=\frac{1}{2p}\sum_{i, j=1}^m \beta_{i j}|u_i|^p|u_j|^p$.}. From a physical point of view, this arises naturally when looking for standing wave solutions ($\Phi_i(x,t)=e^{\imath \lambda_i t} u_i(x)$)) of the following system of Gross-Pitaevskii/nonlinear Schr\"odinger equations:
\begin{equation}\label{S-system-2}
\begin{cases}
\imath \partial_t \Phi_i + \Delta \Phi_i + \Phi_i |\Phi_i|^{p-2}\sum_{j=1}^d\beta_{ij}  |\Phi_j|^p=0, \\
\Phi_i=\Phi_i(t,x),\quad i=1,...d,
\end{cases}
\end{equation}
where $\imath$ is the imaginary unit. These equations model important phenomena in Nonlinear Optics \cite{AkAn} and Bose-Einstein condensation \cite{Rogel_Salazar_2013, Timmermans}. In the models, the solutions are the corresponding condensate amplitudes, $\beta_{ii}$ represent self-interactions within the same component, while $\beta_{ij}$ $(i\neq j)$ express the strength and the type of interaction between different components $i$ and $j$. When $\beta_{ij}>0$ this represents cooperation, while $\beta_{ij}<0$ represents  competition.  Both cases $\Omega = \R^N$ and
$\Omega$ bounded are of interest \cite{FibichMerle2001,Fukuizumi2012}, the latter appearing also as a limiting case of the system in $\R^N$ with
(confining) trapping potential.

Weak solutions of \eqref{S-system} correspond to critical points of the functional $J:H^1_0(\Omega; \R^d)\to \R$ defined by
\[
J(\mf{u})=J\left(u_1, \ldots, u_d\right):=\frac{1}{2} \sum_{i=1}^d\int_\Omega (|\nabla u_i|^2+\lambda_i u_i^2)-\frac{1}{2 p} \int_{\Omega} \sum_{i, j=1}^d \beta_{i j}\left|u_i\right|^p\left|u_j\right|^p.
\]
One is typically interested in least energy solutions, that is, solutions of the following problem
\[
\inf\{J(\mf{u}): \mf{u}\neq \mf{0},\ \mf{u} \text{ solution of \eqref{S-system}}\}.
\]
To prove existence of least energy solutions for $\lambda_i>0$ does not require, in general, different methods from the ones used for the single equation case (recall Section \ref{sec2}). A more challenging question, instead, is the following.

\begin{quote}
Q1: Do least energy solutions $(u_1,\ldots, u_d)$ have nontrivial components, that is, is it true that $u_i\not\equiv 0$ for every $i$?
\end{quote}
If the answer to the above question is negative, then:
\begin{quote}
Q2: Are there solutions (at higher energy levels) satisfying such property? Is there, in particular, a least energy positive solution, as defined below?
\end{quote}
Having this in mind, we make the following definition (see for instance \cite{CorreiaOliveiraTavares,TavaresYou}).
\begin{definition}
A vector $\mathbf{u}=\left(u_1, \cdots, u_d\right)\in H^1_0(\Omega;\R^d)$ is called a \emph{fully nontrivial} solution of \eqref{S-system} if it is a weak solution of the system with $u_i \not \equiv 0$ for every $i=1, \ldots, d$; if this is not the case but nevertheless $\mf{u}\neq \mf{0}$, then we call it \emph{semitrivial}. The vector $\mathbf{u}$ is called a positive solution of  \eqref{S-system} if $\mathbf{u}$ is a solution and $u_i>0$ for every $i=1, \ldots, d$. A positive solution $\mathbf{u}$ is called a least energy positive solution if $J(\mathbf{u}) \leq J(\mathbf{v})$ for any positive solution $\mathbf{v}$ of  \eqref{S-system}. 
\end{definition}

The literature around the subject exploded since the seminal paper \cite{LinWei}, and it would be impossible to cite all contributions. Moreover, some papers deal with the case $\Omega$ bounded, others with the case of the whole space $\Omega=\R^N$; some deal with general $p$, others with particular choices of $p$. In order to highlight my contributions to the topic and to give the reader a coherent and general picture, from now on I will be always referring to the case $\Omega$ bounded and smooth, and treat the case of a general $p$. It should be remarked, however, that I will be mentioning papers which deal with the case $p=2$ only (but the method, in my opinion, works for general $p$), or papers that deal with the case $\Omega=\R^N$ and radial functions, but the methods also work  for $\Omega$ bounded.

We divide our discussion from now on between the subcritical case $2p<2^*$ and the critical one $2p=2^*$.
\subsection{Subcritical case}\label{sec:subcritical}

Consider in this subsection the subcritical case $2p<2^*$. As mentioned before, the existence of least energy solutions is not an issue for $\lambda_i>0$: one can consider the energy level
\[
c:=\inf_{\mathcal{N}} J, \qquad \text{ where } \mathcal{N}=\left\{\mf{u} \not\equiv \mathbf{0}: \ J'(\mf{u})[\mf{u}]=0\right\} \text{ is the Nehari manifold.}
\]
We start our discussion with the case of systems with $d=2$ equations:
\begin{equation}\label{S2-system}
\begin{cases}
-\Delta u_{1}+\lambda_{1}u_{1}=\beta_{11}u_{1}|u_1|^{2p-2}+ \beta_{12}u_{1}|u_1|^{p-2}|u_{2}|^p   ~\text{ in } \Omega,\\
-\Delta u_{2}+\lambda_{2}u_{2}=\beta_{22}u_{2}|u_2|^{2p-2}+ \beta_{12}u_{2}|u_2|^{p-2} |u_{1}|^p   ~\text{ in } \Omega,\\
u_{1}=u_2=0 \text{ on } \partial\Omega.
\end{cases}
\end{equation}
 Regarding question Q1, after preliminary results by \cite{AmbrosettiColorado,MaiaMontefuscoPellacci,Sirakov2007}, a full answer was obtained in \cite{Mandel}.

\begin{theorem}[{{\cite[Theorem 1]{Mandel}}}] There exists $\bar{\beta}=\bar{\beta}(\lambda_2/\lambda_1,\beta_{11},\beta_{22})>0$ such that
\begin{itemize}
\item for $\beta_{12}<\bar{\beta}$, all least energy solutions of \eqref{S2-system} are semitrivial.
\item for $\beta_{12}>\bar{\beta}$, all least energy solutions of \eqref{S2-system} are fully nontrivial.
\end{itemize}
\end{theorem}
Therefore, it makes sense to ask question Q2 for $\beta_{12}<\bar \beta$. In this situation, one needs to consider a Nehari-type set of a different kind:
\[
d_2:=\inf_{\mathcal{N}_2} J, \qquad \text{ where } \mathcal{N}_2=\left\{\mf{u}: u_i\not\equiv 0,\ \ \partial_i E(\mf{u})u_i=0\ \forall i\right\}.
\]
\begin{theorem}[{{\cite{Mandel,Sirakov2007}}}]  There exists $\underline{\beta}=\underline{\beta}(\lambda_2/\lambda_1,\beta_{11},\beta_{22})\leq \bar \beta$ such that a least energy  positive solution exists for $\beta_{12}<\underline{\beta}$.
\end{theorem}
 For the expression of the optimal values $\underline{\beta}$ and $\overline{\beta}$, we refer to \cite{ChenZou, Mandel, Sirakov}. However, in the case $1<p<2$ it is known that $\bar \beta=\underline{\beta}=0$, see \cite[Lemma 1 combined with Lemma 2]{Mandel}.

 I would also like to highlight my joint work with T. Weth \cite{TavaresWeth}, where we study the symmetry of least energy positive solutions in radial domains in the competitive case:
\begin{theorem}[{{\cite[Theorem 1.4]{TavaresWeth}}}]\label{TavaresWeth_symm} Let $\Omega$ be a radial bounded domain, $(u_1,u_2)$ a least energy positive solution of \eqref{S2-system} and $\beta_{12}<0$. Then $u_1,u_2$ are foliated Schwarz symmetric with respect to antipodal points.
\end{theorem}

In general, solutions are not radial, see for instance \cite[Remark 5.4]{TavaresWeth}. On the other hand, if $\Omega$ is a ball and $\beta_{12}>0$, then by Schwarz symmetrization it is easy to see that least energy solutions are radial. Also for $\beta_{12}>0$, the paper \cite{WangWillem} treats the case of the annulus, where symmetry breaking may also occur. To understand the possible symmetries in the $d\geq 3$-equations case remains a challenging open problem.

\medbreak

The question now is what happens when we increase the number of equations. We observe that there is an increase in complexity mainly due to the several possible  combinations of $\beta_{ij}$'s. However, surprisingly enough, we will see that in some situations even the parameters $\lambda_i$ play an important role.

From now on, we focus on the case $4\leq 2p<2^*$ (where the functional is of class $C^2$). Following the introduction in \cite{SoaveTavares}, to summarize the main results which were known before our work, we split our discussion into several cases. We first focus on three situations where we have the \emph{same} type of interaction terms, providing some references in which existence of least energy positive solutions is proved.
\begin{itemize}
\item \emph{Strong cooperation}: $\lambda_1=\dots=\lambda_d=\lambda>0$, $\beta_{ii}>0$, and $\beta_{ij}=\beta$ (for every $i \neq j$) larger than a positive constant depending on $\beta_{ii}$ and $\lambda$ (see Corollary 2.3 and Theorem 2.1 in \cite{LiuWang}; see also Theorem 1.6 and Remark 3 in \cite{Soave}). Other sufficient conditions in a purely cooperative setting have been given in  \cite[Section 4]{Sirakov}, \cite[Theorem 2.1]{LiuWang} and \cite{Chang}.
\item \emph{Weak cooperation}: $\lambda_i >0$, $\beta_{ii}>0$, $0<\beta_{ij}\leq \Lambda$ for some small $\Lambda$ depending on $\lambda_i$ and $\beta_{ii}$, and the matrix $(\beta_{ij})$ is positive definite (see Theorem 2 in \cite{LinWeiErratum});
\item \emph{Competition}: if $\lambda_i >0$, $\beta_{ii}>0$, and $\beta_{ij} \le 0$ for every $i \neq j$, then there exists a least energy positive solution (see Theorem 1.1. plus Remark 1.5 in \cite{LinWei2}; we refer also to Theorem 3.1 in \cite{LiuWang}, and to Corollary 1.4 plus Proposition 1.5 in \cite{Soave}).
\end{itemize}
It is natural to assume that $\beta_{ij}$ is either large, or small, with respect to $\beta_{ii}$ and $\beta_{jj}$. Indeed, if, for instance $\beta_{ii} \le \beta_{ij} \le \beta_{jj}$ and $\lambda_i>\lambda_j$, then a positive solution of \eqref{S-system} does not exist, see Theorem 1-($ii$) in \cite{Sirakov} or Theorem 0.2 in \cite{BartschWang}.

As far as the possible occurrence of simultaneous cooperation and competition is concerned, in \cite{SatoWang}  a $d=3$ components system is considered, showing that a least energy positive solution of \eqref{S-system} does exist if $\beta_{13},\beta_{23} \le 0$, and $\beta_{12} \gg 1$ is very large (depending on $\beta_{13}$ and $\beta_{23}$ fixed \emph{a priori}, which is a technical downsize that can be removed, see the upcoming paragraph \textbf{Mixed coefficients case}). In \cite[Theorems 1.6, 1.7 and 1.9]{Soave} the author considered an arbitrary $d$-component system, proving the existence of least energy positive solutions whenever the $d$ components are divided into $m$ groups, with $m \le d$, and
\begin{itemize}
\item the relation between components of the same group is purely cooperative, with coupling parameters greater than an explicit positive constant;
\item the relation between components of different groups is competitive, and the competition is very strong.
\end{itemize}
When restricted to a $3$-component system, this leads for instance to the existence of a least energy solution if $\beta_{12}> \overline{\beta}>0$, and $\beta_{13},\beta_{23} \ll -1$ (depending on $\beta_{12}$, which again is a downsize). 

After giving this context, we now summarize our main contributions to the field.

\paragraph{Cooperative case} ($\beta_{ij}>0$ for every $i,j$) Even though in the $d=2$ equations case a result holds for arbitrary $\lambda_1<\lambda_2$, all results cited in the strongly cooperative case for $d\geq 3$ impose  $\lambda_i\sim \lambda$ and $\beta_{ij}\sim \beta$. In the case $2<2p<4$, with F. Oliveira \cite{OliveiraTavares} we observed (using an induction argument in the number of equations) that the situation is the same as in the two equation case:
\begin{theorem}[{{\cite{OliveiraTavares}}}]
 Let $N\geq 1$, $\lambda_i>0$, $\beta_{ij}>0$  for every $i,j=1,\ldots, d$, and $\beta_{ij}=\beta_{ji}$ for $i\neq j$. For $1<p<2$, all possible least energy solutions of \eqref{S-system} are fully nontrivial.
\end{theorem}
As for the case $4\leq 2p<2^*$,  it was proven in \cite{CorreiaNA2016} that, when $\lambda_1=\dots=\lambda_d$, these questions may be reduced to a maximization problem in $\R^d$ and to the solution of a linear system. This reduction allowed the construction of examples (see Section 6 in \cite{CorreiaNA2016}) which gave evidence, for the first time, of the increase in complexity when one passes from $d=2$ to $d\ge3$ equations. Indeed, we stated qualitatively in \cite{CorreiaOliveiraTavares} (joint with S. Correia and F. Oliveira) what kind of combinations on the parameters give rise either to semitrivial or to fully nontrivial least energy solutions. In particular, it became evident from our analysis that the different families of parameters play distinct roles: while the choice of the $\beta_{ii}$ coefficients can be somehow arbitrary, only some combinations between different $\lambda_i$, and also between different $\beta_{ij}$'s allow for fully nontrivial least energy solutions to arise.

\begin{theorem}[{{\cite{CorreiaOliveiraTavares}}}]\label{thm:COT_main}
Let $d\geq 3$, $0<\lambda_1\le \lambda_2\le \cdots \le \lambda_d$ and $\beta_{ij}\equiv \beta$.
\begin{enumerate}
\item (Existence result) There exists  $\alpha=\alpha(\lambda_1/\lambda_2,d,N)$ such that, if
$$
\lambda_k \leq \alpha \lambda_2\qquad \text{ for every } k\neq 2,
$$
then there exists a constant $B=B(\beta_{ii})>0$, such that, for $\beta>B$, all least energy solutions of \eqref{S-system} are fully nontrivial.
\item (Nonexistence result) There exists a constant $\Lambda=\Lambda(\lambda_1/\lambda_2)$ such that, if
\[
\lambda_2 \Lambda \leq \lambda_i \text{ for some } i\geq 3, \text{ and } \beta>\max\{\beta_{11},\ldots, \beta_{dd}\},
\]
then every least energy solution of \eqref{S-system} is semitrivial (more precisely, $u_j\equiv 0$ for every $i\geq j$).
\end{enumerate}
\end{theorem}
Therefore, for $d\geq 3$, in a way, only perturbations of the 2--equation case (in terms of the parameters) allow for least energy solutions which are fully nontrivial.

In \cite{CorreiaOliveiraTavares}, we have similar results regarding the parameters $\beta_{ij}$: least energy solutions are fully nontrivial if these coefficients are large and ``close'' to each other; otherwise, they necessarily become semitrivial.  The proofs are based on classification results, comparison of energies between the main systems and appropriate subsystems, and \emph{a priori} bounds.

A natural open question is whether there are positive solutions under the range of parameters where least energy solutions are semitrivial, and how to characterize them variationally.

\paragraph{Mixed coefficients case}




Having in mind the idea of organizing the components of a solution to the system into $m\leq d$ groups, we follow \cite{Soave,SoaveTavares}.

\begin{definition}\label{def:groups}
 Given an arbitrary $1\leq m \leq d$, we say that a vector $\mathbf{a}=(a_{0},...,a_{m})\in \mathbb{N}^{m+1}$ is an $m$-decomposition of $d$ if
  \begin{equation*}
  0=a_{0}<a_{1}<\cdot \cdot \cdot<a_{m-1}<a_{m}=d.
  \end{equation*}
Given an $m$-decomposition $\mathbf{a}$ of $d$, for $h=1,...,m$ we define
\begin{align*}
&  I_{h}:=\left\{i\in \{1,...,d\}:a_{h-1}<i\leq a_{h}\right\},
\end{align*}
and
\begin{align*}
& \mathcal{K}_{1}:=\left\{(i,j)\in I_{h}^{2}  \text{ for some } h=1,...,m, \text{ with } i\neq j\right\},\\
& \mathcal{K}_{2}:=\left\{(i,j)\in I_{h}\times I_{k} \text{ with } h\neq k\right\}.
  \end{align*}
 In this way, we  say that $u_i$ and $u_j$ belong to the same group if $(i,j)\in \mathcal{K}_1$ and to a different group if $(i,j)\in \mathcal{K}_2$.
\end{definition}
As we will see below, the general idea is that we obtain existence results wherever the interaction between elements of the same group is strongly cooperative, while there is either weak cooperation or competition between elements of different groups.

The following is our main result, from a project with N. Soave.

\begin{theorem}[{{\cite{SoaveTavares}}}]\label{thmST}\
\begin{enumerate}
\item There exists $K=K(\lambda_i,\beta_{ii})>0$ such that, if
\[
-\infty< \beta_{ij} < K \qquad \text{for every $i \neq j$},
\]
then the system \eqref{S-system} admits a least energy positive solution.
\item Consider a decomposition of $\{1,\ldots, d\}=I_1\cup\ldots\cup I_m$ and assume the following.
\begin{enumerate}
\item[$i$)] Inside each group $I_h$:
\begin{align*}
&\beta_{ij}\equiv \beta_h > \max\{\beta_{ii}: i \in I_h\}\text{ for every } (i,j) \in I_h^2\text{ with } i \neq j;\\
&\lambda_i \equiv \lambda_h \text{ for every } i \in I_h
\end{align*}
\item[$ii$)] Between different groups: there exists $K=K(\lambda_i,\beta_{ii})>0$:
\[
\beta_{ij}=\beta < K \text{ for every } (i,j)\in \mathcal{K}_2;
\]
\end{enumerate}
Then the system \eqref{S-system} admits a least energy positive solution.
\end{enumerate}
\end{theorem}

In summary, if we particularize the discussion to the $d=3$ equations case, our results combined with what was known allows for a good understanding of the bigger picture: there exists $0<\underline{\beta}\leq \overline{\beta}$ such that the system  \eqref{S-system}  admits a least energy positive solution when one of the following conditions is verified:
\[
\begin{split}
& \beta_{12}\sim \beta_{13}\sim \beta_{23}>\overline{\beta};\quad  \lambda_1< \lambda_2\sim \lambda_3\\[5pt]
& \beta_{12}>\overline{\beta}\text{ and }-\infty<\beta_{13}=\beta_{23}<\underline{\beta};\qquad \lambda_1=\lambda_2 \\[5pt]
& -\infty<\beta_{12},\beta_{13},\beta_{23}<\underline{\beta}.
\end{split}
\]
In particular, we improve the dependences of some $\beta_{ij}$ which were present in the aforementioned \cite{SatoWang,Soave}, in the sense that $\underline{\beta}$ and $\bar \beta$ only depend on $\lambda_i$ and $\beta_{ii}$ for $i=1,\ldots, d$. Observe that having cooperative parameters too far apart, or too different $\lambda_i$, may lead to semitrivial solutions (recall Theorem \ref{thm:COT_main} and the paragraph that follows it).

A partial symmetry result is also proved in \cite{SoaveTavares} in the case of $m=2$ groups of components (in the line of Theorem \ref{TavaresWeth_symm}); however, the symmetry of the general case is, up to our knowledge, an open problem.

The variational formulation associated with these solutions would be too technical to explain in this document; here I just give the idea that they are Nehari--type sets with $m$ equations; it is quite straightforward to prove that the associated critical points have at least $m$ nontrivial components; to prove that all components are nontrivial, we make use of the $C^2$ regularity of the functional $J$ in the case $4\leq 2p$.

We conclude our section by mentioning that  other related and recent results in the subcritical case can be found in \cite{ClappPistoia2022, ClappSzulkin2019, CorreiaJDE2016,CorreiaNA2016,OliveiraTavares,PengWangWang2019,WeiWu}.
\subsection{Critical case}\label{sec:criticalcase_gradsystem}

For the \emph{critical} case $2p=2^*$, when $d=1$, system \eqref{S-system} is reduced to the classical Br\'ezis-Nirenberg problem \cite{BrezisNirenberg} (recall also Subsection \ref{eq:critical_1eq}):
\[
-\Delta u+\lambda_1 u=\beta_{11} u_1 |u_1|^{2p-2}, \quad u=0 \text{ on } \partial \Omega,
\] where the existence of a positive ground state is shown for $-\lambda_1(\Omega)<\lambda_1<0$ when $N\geq 4$, where $\lambda_1(\Omega)$ is the first Dirichlet eigenvalue. For the $d=2$ equation case \eqref{S2-system}, in \cite{ChenZouARMA2012} it is shown that there exist $0<\beta_{1}<\beta_{2}$ (depending on $\lambda_{i}$ and $\beta_{ii}$) such that
\begin{multline}\label{Range-1}
\text{the system } \eqref{S2-system}\text{ has a least energy positive solution if }\\
  \beta_{12}\in (-\infty, \beta_{1}) \cup (\beta_{2}, \infty)\text{ when } N=4,\ p=2.
\end{multline}
We mention that, when $p=2$ and $\beta\in [\min\{\beta_{11},\beta_{22}\},\max\{\beta_{11},\beta_{22}\}]$, system \eqref{S2-system} does not have a least energy positive solution.
Still for $d=2$ equations but in the higher dimensional case, the same authors in \cite{ChenZou2} proved that
\begin{multline}\label{Range-2}
\text{the  system }\eqref{S2-system} ~\text{ has a least energy positive solution for any }\\ \beta_{12}\neq 0 ~\text{ when } N\geq 5,\ 2p=2^*.
\end{multline}
In particular, from \eqref{Range-1} and \eqref{Range-2}, one deduces that the structure of least energy positive solutions in the critical case changes significantly from $N=4$ to $N\geq 5$. Once again, the reason behind this change is the fact that $p\in (1,2)$ whenever $N\geq 5$, while $p=2$ for $N=4$ (we recall that, in the subcritical case, the importance of this fact has been implicitly pointed out by Mandel in \cite{Mandel}, see also \cite{OliveiraTavares}). For more results regarding the general critical case with $d=2$ equations, see \cite{ChenLinZou,ClappPistoia2018,PengPengWang2016}.

For three or more equations ($d\geq 3$), in the critical case $2p=2^*$, before our work only the  purely competitive case \cite{ClappSzulkin2019,Wu} and the purely cooperative case \cite{YinZou} had been studied, and conditions for the existence of least energy positive solutions had been provided.

Together with S. You  \cite{TavaresYou}, working with $N=4$, we considered for the first time the critical case with simultaneous cooperation and competition; the higher dimensional case $N\geq 5$ was treated later in a collaboration with S. You and W. Zou \cite{TavaresYouZou}. We make the following assumptions:
\begin{equation}\label{eq:coefficients}
-\lambda_1(\Omega)<\lambda_1,\ldots, \lambda_d<0,\qquad \Omega \text{ is a bounded smooth domain of } \mathbb{R}^N,
\end{equation}
and
\begin{equation}\label{eq:beta_ij}
\beta_{ii}>0 \quad \forall i=1,\ldots, d,\qquad \beta_{ij}=\beta_{ji}\quad \forall i, j=1,\ldots, d,\ i\neq j.
\end{equation}

%
%

\begin{theorem}[{{\cite{TavaresYou}}}, critical case $N=4$]\label{thmSY} Fix an  $m$-decomposition $\mathbf{a}$ of $d$, for some $1< m< d$. There exists a least energy positive solution under \eqref{eq:coefficients}--\eqref{eq:beta_ij} with  $N=4$  ($p=2$) in each one of the following situations:
\begin{itemize}
\item $-\infty<\beta_{ij}<\Lambda \quad\forall i\neq j$, for some $\Lambda>0$ depending only on $\beta_{ii}, \lambda_{i}$;
\item $\lambda_{i}=\lambda_{h}$ for every $i\in I_{h}, h=1,\ldots, m$;\newline $\beta_{ij}=\beta_{h}>\max\{\beta_{ii}: i\in I_{h}\}$ for every $(i,j)\in I_{h}^{2}$ with $i\neq j, h=1,\ldots, m$;\newline $\beta_{ij}=b<\Lambda$ for every $(i,j)\in \mathcal{K}_{2}$;
\item $\lambda_{i}=\lambda_{h}$ for every $i\in I_{h}, h=1,\ldots, m$;\newline  $\beta_{ij}=\beta_{h}>\frac{\alpha}{\alpha-1}\max_{i\in I_h}\{\beta_{ii}\}$ for every $(i,j)\in I_{h}^{2}$ with $i\neq j, h=1,\ldots, m$;\newline  $ |\beta_{ij}|\leq \frac{\Lambda}{\alpha d^{2}}$ for every $(i,j)\in \mathcal{K}_{2}$;
\end{itemize}
Here $\Lambda$ is a precise constant and $\alpha > 1$ is arbitrary.
\end{theorem}

\begin{theorem}[{{\cite{TavaresYouZou}}}, critical case $N\geq 5$]\label{thm:leps} Assume that \eqref{eq:coefficients} and \eqref{eq:beta_ij} hold, $N\geq 5$. Then the system admits a least energy positive solution in each one of the following situations:
\begin{enumerate}
\item  $\beta_{ij}> 0 \ \forall i,j=1,\ldots, d$, $i\neq j$.
\item  $\beta_{ij}\leq 0 \ \forall i,j=1,\ldots, d$, $i\neq j$.
\item  $\mathbf{a}$ is an $m$-decomposition of $d$ for some $1< m< d$, and
\[
\beta_{ij}\geq 0  \ \forall (i,j)\in \mathcal{K}_{1},\qquad -\varepsilon\leq\beta_{ij}< 0\ \forall (i,j)\in \mathcal{K}_{2},
\]
 for some $\varepsilon=\varepsilon(\lambda_i,\beta_{ii},(\beta_{ij})_{(i,j)\in \mathcal{K}_1})>0$;
\item  $\mathbf{a}$ is an $m$-decomposition of $d$ for some $1< m< d$, we have $\beta_{ij}\geq 0  \ \forall (i,j)\in \mathcal{K}_{1}$,
and for every $M>1$ there exists $b=b(\lambda_i,\beta_{ii},(\beta_{ij})_{(i,j)\in \mathcal{K}_1},M)>0$ such that
\[
\frac{1}{M}\leq \left|\frac{\beta_{i_1j_1}}{\beta_{i_2j_2}}\right|\leq M, ~~\forall (i_1,j_1), (i_2,j_2)\in \mathcal{K}_{2}\quad \text{ and }\quad  \beta_{ij}\leq -b ~~\forall (i,j)\in \mathcal{K}_2.
\]
\end{enumerate}
Moreover, under case 1, the solution is a least energy solution.
\end{theorem}
The results in Theorem \ref{thmSY} are similar to the ones in Theorem \ref{thmST}; in both, the fact that the functional $J$ is of class $C^2$  is explored. For $N\geq 5$ this is not the case and this is why the results of Theorem \ref{thm:leps} (and the techniques used to prove it) are different.  The proofs of cases 3 and 4 in Theorem \ref{thm:leps} are based on an asymptotic study as $\beta_{ij}^n\to 0^-$ and $\beta_{ij}^n\to -\infty$ for $(i,j)\in \mathcal{K}_2$, respectively. This study also allows to answer open questions in the literature, see the forthcoming paragraph ``A few works about the strongly competing case''.
Based upon Theorem \ref{thm:leps}, it is natural to ask what happens when some interactions between elements of different groups are neither too strong nor too weak. For simplicity and to avoid too technical conditions, we present here the following result for the case of only two groups ($m=2$), for the partition $I_1=\{1,2\}$, $I_2=\{3\}$.

\begin{theorem}\label{existence-positive-two groups2}
Assume that $m=2$, $d=3$ and let $0<\sigma_0<\sigma_1$. Then there exists $\widehat{\beta}=\widehat{\beta}((\sigma_i)_i,(\beta_{ii})_i,(\lambda_i)_i)>0$ such that, if
\[
\beta_{13},\beta_{23} \in [-\sigma_1,-\sigma_0],\qquad \beta_{12}>\widehat{\beta},
\] then the system \eqref{S-system} has a least energy positive solution.
\end{theorem}

We leave a few open problems for the case $N\geq 5$.  Does a least energy positive solution exist under the general conditions $\beta_{12}> 0$ and $\beta_{13},\beta_{23}<0$? This seems a natural generalization of \eqref{Range-2}. And what happens when  more equations are present? Is it true that, in general, a least energy solution exists when $\beta_{ij}>0$ for $(i,j)\in \mathcal{K}_1$, and  $\beta_{ij}<0$ for $(i,j)\in \mathcal{K}_2$? Is it possible to obtain optimal thresholds, if this is not the case?

\medbreak

 For other topics related to critical systems (e.g. blowing up solutions as $\lambda_i\to 0^-$ or the case $\lambda_i=0$ in the whole space), see \cite{ChenLin, ClappPistoia2022, DovettaPistoia, GuoLuoZou,HeYang,PistoiaSoave,PistoiaTavares,PistoiaSoaveTavares}. In particular, I highlight a joint work with Angela Pistoia \cite{PistoiaTavares}, where the study of blowing up solutions as $\lambda_i\to 0^-$ is treated in the spirit of \cite{HanAIHP,MussoPistoiaIndiana2002,Rey1}. I also mention the works \cite{PistoiaSoaveTavares,PistoiaTavares} where, for the first time, a Coron-type problem for these systems was studied (the second is a collaboration between myself, A. Pistoia and N. Soave). These three works use the Lyapunov-Schmidt reduction method, and we recall Subsection \ref{eq:critical_1eq} for the references in the 1--equation case.

\paragraph{A few works about the strongly competing case} The asymptotic study of \eqref{S-system} as $\beta_{ij}\to -\infty$ for $(i,j)\in \mathcal{K}_2$ has been performed in a joint publication with N. Soave, S. Terracini and A. Zilio \cite{STTZsurvey}. Therein, it is showed that uniform bounds in $L^\infty$--norm imply uniform bounds in H\"older spaces, which then allow to pass to a strong limit $u_i$ as $\beta_{ij}\to -\infty$. Moreover, the regularity of the common nodal set  $\Gamma:=\{x\in \Omega:\ u_i(x)=0\ \forall i\}$ is studied. This follows previous work \cite{NTTV1, TavaresTerracini2}.  All least energy positive solutions satisfy these uniform $L^\infty$--bounds, therefore the results in \cite{STTZsurvey} may be used. While in the subcritical case $2p<2^*$ it is straightforward to check that all components do not vanish in the limit \cite{Soave}, this is not as easy to check in the critical case $2p=2^*$, due to the lack of compactness in some Sobolev embeddings. We proved it in \cite{TavaresYouZou}, and present here the actual statement in the case $m=d$:

\begin{corollary}[Combination of \cite{STTZsurvey} with \cite{TavaresYouZou}]\label{Phase Separation2}
Assume that $N\geq 4, d\geq2, -\lambda_1(\Omega)<\lambda_i<0$. Let $\beta_{ij}^n<0, n\in \mathbb{N}$ and $\beta_{ij}^n\rightarrow -\infty$ when $i\neq j$, and let $(u_1^n,\ldots, u_d^n)$ be a least energy positive solution with $\beta_{ij}=\beta_{ij}^n$. Then, passing to a subsequence, we have
\begin{equation*}
u_i^n\rightarrow u_i^\infty \text{ strongly in } H^1_0\cap C^{0,\alpha}(\Omega), ~~i=1,\ldots,d, \ \alpha\in (0,1),
\end{equation*}
and
\begin{equation*}
 \lim_{n\rightarrow\infty}\int_{\Omega}\beta_{ij}^n|u_i^n|^p |u_j^n|^p= 0,\quad \text{ and } \quad  u_i^\infty\cdot u_j^\infty\equiv 0 ~\text{ for every } i\neq j,
\end{equation*}
where $u_i^\infty\in C^{0,1}(\overline{\Omega})$ and $\overline{\Omega}=\bigcup_{i=1}^d\overline{\{u_i^\infty>0\}}$.
Moreover, for every $i=1,\ldots,d$, $\{u_i^\infty>0\}$ is a connected domain, and $u_i^\infty$ is a least energy positive solution of
\begin{equation*}
-\Delta u+\lambda_iu=\beta_{ii}|u|^{2^*-2}u, \quad u\in H^1_0(\{u_i^\infty>0\}).
\end{equation*}
Finally, the set $\Gamma:=\{x\in \Omega:\ u_i^\infty(x)=0\ \forall i\}$ is, up to a subset of Hausdorff dimension at most $N-2$, a collection of regular hypersurfaces.
\end{corollary}

This result answers a question left open in \cite[Remark 1.4(i)]{ChenZou2}, namely, whether or not the limiting configuration was fully nontrivial. So far this was known to be the case only for  $N\geq 9$ and $\beta_{ij}\equiv \beta$ for every $i\neq j$ (see \cite[Theorem 1.3]{Wu}). We have shown the answer is always positive in  general.

 As an application of Corollary \ref{Phase Separation2}, we can obtain the existence of a least energy nodal solution to the Br\'ezis-Nirenberg problem also in some lower dimensions:
\begin{equation}\label{B-N problem}
-\Delta u+ \lambda u=\mu |u|^{2^*-2}u, \quad u\in H^1_0(\Omega),
\end{equation}
where $\mu>0, -\lambda_1(\Omega)<\lambda<0$ and $N\geq 4$.
\begin{theorem}[{{\cite{TavaresYouZou}}}]\label{Sign-changing solutions}
Assume that $N\geq 4$, $d=2$, $\lambda_1=\lambda_2=\lambda\in (-\lambda_1(\Omega),0)$, and $\beta_{11}=\beta_{22}=\mu>0$. Let $(u_1^\infty, u_2^\infty)$ be as in Corollary \ref{Phase Separation2} for $d=2$. Then $u_1^\infty - u_2^\infty$ is a least energy nodal solution of \eqref{B-N problem} and has two nodal domains.
\end{theorem}

The existence of sign-changing solutions (not necessarily of least energy) to the Br\'ezis-Nirenberg problem \eqref{B-N problem} has been studied in \cite{CeramiSoliminiStruwe,ChenZou2,HeHeZhang,HeZhangWuLiang,SchechterZou} with $N\geq 4$ in a general domain. In some symmetric domains, see \cite{AtkinsonBrezisPeletier,CastroClapp,ClappWeth}. We mention that in \cite{CeramiSoliminiStruwe, ChenZou2} the authors proved the existence of least energy nodal solutions for $N\geq 6$. However, there are few results considering the  lower-dimensional situations $(N=4,5)$ in the literature. Recently, the authors in \cite{RoselliWillem} proved that \eqref{B-N problem} has a least energy sign-changing solution for $N=5$ with $\lambda\in (-\lambda_1(\Omega),-\overline{\lambda})$, for some $\overline{\lambda}\in(0,\lambda_1(\Omega))$. Here, we improve and extend this result to the case $N\geq 4$.

\subsection{A quick detour on normalized solutions for NLS in bounded domains}

In all previous subsections, the coefficients $\lambda_i$ are fixed \emph{a priori}. Here we briefly explore a different perspective. To simplify the presentation (since the focus here is not on the number of equations), we deal with the  $d=2$ equations case only, where the system \eqref{S-system-2} becomes (with $\beta=\beta_{12}=\beta_{21}$, $\mu_1=\beta_{11}$ and $\mu_2=\beta_{22}$)

\begin{equation}\label{eq:system_schro}
\begin{cases}
\imath\partial_t\Psi_1 + \Delta\Psi_1 + \Psi_1( \mu_1 |\Psi_1|^{2p-2} +\beta  |\Psi_1|^{p-2}|\Psi_2|^{p} )=0\\
\imath\partial_t\Psi_2 + \Delta\Psi_2 + \Psi_2( \mu_2 |\Psi_2|^{2p-2} +\beta  |\Psi_2|^{p-2}|\Psi_1|^{p} )=0.
\end{cases}
\end{equation}

The flow generated by solutions to the system \eqref{eq:system_schro} preserves, at least formally, the \emph{masses}
\[
\mathcal{Q}(\Psi_1(t))=\int_\Omega |\Psi_1(t)|^2, \quad \mathcal{Q}(\Psi_2(t))=\int_\Omega |\Psi_2(t)|^2.
\]
We look for standing wave solutions $(\Psi_1(t,x),\Psi_2(t,x))=(e^{\imath \lambda_1 t}u_1(x),e^{\imath \lambda_2 t}u_2(x))$ of \eqref{eq:system_schro} such that $(u_1,u_2) \in H^1_0(\Omega;\R^2)$ and
\begin{equation}\label{eq:mass_constraint}
\mathcal{Q}(\Psi_1(t))\equiv \mathcal{Q}(u_1)=\rho_1, \quad \mathcal{Q}(\Psi_2(t))\equiv \mathcal{Q}(u_2)=\rho_2,
\end{equation}
for some prescribed $\rho_1,\rho_2\ge0$. Therefore, unlike fixing $\lambda_1,\lambda_2$ \emph{a priori} as in the previous subsections and simply looking for solutions of \eqref{S2-system}, here we look for \emph{normalized solutions} of \eqref{S2-system}; namely, we ask if, given $\rho_1,\rho_2\geq 0$, there exists  $\lambda_1,\lambda_2\in \R$ and $u_1,u_2 \in H^1_0(\Omega)$ so that
\begin{equation}\label{eq:system_elliptic}
\begin{cases}
-\Delta u_1+ \lambda_1 u_1=\mu_1 u_1|u_1|^{2p-2}+\beta u_1|u_1|^{p-2} |u_2|^{p}\\
-\Delta u_2+ \lambda_2 u_2=\mu_2 u_2|u_2|^{2p-2}+\beta u_2 |u_2|^{p-2} |u_1|^{p}\\
\int_\Omega u_i^2=\rho_i, \quad i=1,2,\\
u_1,u_2 \in H^1_0(\Omega).
\end{cases}
\end{equation}
Throughout this section we are interested in positive solutions.

Solutions of \eqref{eq:system_elliptic} can be seen as critical points of the energy functional
\[
\mathcal{E}(\Psi_1,\Psi_2) := \frac12\int_{\Omega}|\nabla \Psi_1|^2+|\nabla \Psi_2|^2 -
\frac{1}{2p}\int_\Omega \mu_1 |\Psi_1|^{2p} + 2\beta |\Psi_1|^{p} |\Psi_2|^{p} + \mu_2 |\Psi_2|^{2p}.
\]
(another quantity formally preserved by the flow generated by the solutions of system \eqref{eq:system_schro})  constrained to the manifold
\begin{equation}\label{eq:defM}
\mathcal{M}_{\rho_1,\rho_2} := \left\{(u_1,u_2)\in H^1_0(\Omega;\R^2):\int_\Omega u_1^2=\rho_1, \ \int_\Omega u_2^2=\rho_2\right\}.
\end{equation}
From this point of view, the main aim is to
provide conditions on $p$ and $\rho_1,\rho_2$ (and also on $\mu_1,\mu_2,\beta$) so that $\left.\mathcal{E}\right|_{\mathcal{M}_{\rho_1,\rho_2}}$ has critical points or, more specifically, if it
 admits minima, either global or local.  We call such solutions \emph{energy ground states} (in the literature, the least energy solutions studied in the previous subsections are also called \emph{action ground states}). In this context, the \emph{unknowns} $\lambda_1,\lambda_2$ appear as  Lagrange multipliers. As a second aim, one considers the stability properties of such ground states with respect
to the evolution system \eqref{eq:system_schro}.

The simplest case one can face is that of a single Nonlinear Schr\"odinger (NLS) equation in $\R^N$, with a pure power nonlinearity:
\begin{equation}\label{fixedmass_RN}
-\Delta u_1+\lambda_1 u_1=\mu_1 u_1|u_1|^{2p-2} \text{ in } \R^N,\quad \int_{\R^N}u_1^2=\rho_1,\quad  u_1\in H^1(\R^N).
\end{equation} In such case, the problem can be completely solved by
simple scaling arguments. Indeed, it is known that such problem, up to translation, has a unique positive solution; if we denote by $Z$ the unique radial (decreasing) solution for $\lambda_1=1$, see \cite{Kwong}, then $u(x)=hZ(h^{p-1}x)$, $h>0$, solves \eqref{fixedmass_RN} with $\lambda_1=h^{2p-1}$ and $\rho_1=h^{2+N(1-p)}\|Z\|_{L^2(\R^N)}^2$. Therefore, for $p\neq 1+2/N$, the problem \eqref{fixedmass_RN} admits a positive solution for every value of the mass, while for $p=1+2/N$ it admits a  positive solution only for the mass $\rho_1=\|Z\|^2_{L^2(\R^N)}$. This, among other things, leads to the classification of the exponent $p$ in \eqref{eq:system_elliptic} according to the following four cases:
\begin{itemize}
\item[(H1)] superlinear, $L^2$-subcritical: $1<p<1+2/N$;
\item[(H2)] $L^2$-critical: $p=1+{2/N}$;
\item[(H3)] $L^2$-supercritical, Sobolev--subcritical: $1+2/N<p<2^*/2$;
\item[(H4)] Sobolev--critical: $p=2^*/2$, for $N\geq 3$.
\end{itemize}
Moreover, observe that the solutions found in the $L^2$-subcritical case  (H1) are associated with orbitally stable solitary waves of the corresponding evolution equation, while in the remaining cases there is instability \cite{CazenaveLions1982, Cazenave2003}.

However, whenever  one considers a system, as well as non-homogeneous nonlinearities, bounded domains or confining potentials, the situation cannot be solved by such simple scaling arguments. Apart from when global minimization
can be applied, see  \cite{RoseWeinstein88}, as far as we know the first
result in the literature is due to Jeanjean \cite{Jeanjean}, for the superlinear, Sobolev-subcritical
NLS single equation on $\R^N$ with a non-homogeneous nonlinearity. In recent years, other papers appeared, dealing with the NLS equation or system, always in the Sobolev subcritical regime,  either on $\R^N$ \cite{BartschJeanjean,BartschJeanjeanSoave,BartschSoave2,Bellazzinietal,GuoJeanjean,GuoJeanjean2,BartschSoave}
or on a bounded domain \cite{CirantVerzini,MR2928850,ntvAnPDE,NTV2,MR3689156}. In this short subsection, we focus mainly in our contributions to the bounded domain case.

These two settings ($\Omega$ bounded and $\Omega=\R^N$) are rather different in nature: each one requires a specific approach, and the results are in general not comparable.
A key difference is that $\R^N$ is invariant under translations and dilations, which has advantages and disadvantages: on the one hand, translations are responsible for a loss of compactness;
on the other hand, in the Sobolev subcritical case, dilations can be used to produce variations and eventually construct natural constraints such as the so-called Pohozaev manifold.
This tool is not available when working in bounded domains, and also the gain of compactness is lost when we face the Sobolev critical case. However, a common key tool in the study of normalized solution is the Gagliardo-Nirenberg inequality, which can be used to estimate the non-quadratic part in $\mathcal{E}$
in terms of the quadratic one, which also leads naturally to the threshold $p=1+2/N$ appearing in the classification (H1)-(H4).

In the first three cases (H1)-(H3), the study of the single equation
\begin{equation}\label{eq:singleeq}
-\Delta u_1 + \lambda_1 u=\mu_1 u_1|u_1|^{2p-2} \text{ in } \Omega,\qquad \int_\Omega u_1^2=\rho_1,\quad u_1\in H^1_0(\Omega),
\end{equation}
in a \emph{bounded} domain has been carried out in \cite{ntvAnPDE,MR3689156}, where the first is a joint work with B. Noris and G. Verzini. Notice that \eqref{eq:singleeq} is a particular case of \eqref{eq:system_elliptic}, when $\rho_2=0$, with associated energy
$u_1\mapsto \mathcal{E}(u_1,0)=\mathcal{E}(u_1)$. We also denote $\mathcal{M}_{\rho_1}:=\mathcal{M}_{\rho_1,0}$. Summarizing, it is known that
\begin{itemize}
 \item (H1) implies that \eqref{eq:singleeq} has a solution for every $\rho_1$, which is a global minimizer of $\mathcal{E}|_{\mathcal{M}_{\rho_1}}$;
  \item (H2) implies that \eqref{eq:singleeq} has a solution for $0\le\rho_1<\rho_*(\Omega,N,p,\mu_1)<+\infty$, which is a global minimizer of $\mathcal{E}_{\mathcal{M}_{\rho_1}}$;
 \item (H3) implies that \eqref{eq:singleeq} has at least two  solutions for $0\le\rho_1<\rho^*(\Omega,N,p,\mu_1)<+\infty$, and one of these is  a local minimizer of $\mathcal{E}|_{\mathcal{M}_{\rho_1}}$.
\end{itemize}
More precise results are given if $\Omega=B_1(0)$. Moreover, all the minimizers above are associated with orbitally stable solitary waves of the corresponding evolution equation. This shows, in particular, that the boundary has a stabilizing effect;  in the $L^2$-critical and
$L^2$-supercritical cases there exist standing waves which are orbitally stable (which, we recall, is not the case in the whole $\R^N$).

Up to our knowledge, the first paper dealing with the NLS system \eqref{eq:system_elliptic} (with both $\rho_i>0$) is another joint work with B. Noris and G. Verzini \cite{NTV2}. Among other things, in this paper we deal with the
$L^2$-supercritical, Sobolev--subcritical case (H3), obtaining the existence of orbitally stable solitary waves, in case both $\rho_1$, $\rho_2$ are sufficiently small and $\rho_1/\rho_2$ is uniformly bounded away from $0$ and $+\infty$. This result is perturbative in nature. The existence results follow by a multi-parametric extension of a Ambrosetti-Prodi-type reduction \cite{AmbrosettiProdiBook}, while the stability follows from  the
Grillakis-Shatah-Strauss stability theory \cite{GrillakisShatahStrauss}.

In a third paper with B. Noris and G. Verzini \cite{NTV3}, on the one hand, in the cases (H1)-(H2)-(H3) we extend to systems defined in a bounded domain the above described results
\cite{ntvAnPDE,MR3689156} for the single equation; on the other hand, we treat for the first time the Sobolev critical case (H4), obtaining results which are
new also in the case of a single equation.

In conclusion, our works remain important as they showed that the presence of the boundary has a stabilizing effect, complementing earlier observations by
\cite{FibichMerle2001,Fukuizumi2012}. In these papers, it is proved that also in the $L^2$-critical and
$L^2$-supercritical cases there exist standing waves which are orbitally stable.

The topic of normalized solutions has been very active in the past few years, mainly in the case of $\Omega=\R^N$, and it would be impossible and out of the scope of this document to mention all contributions and do a state of the art. Therefore, we  conclude by simply mentioning the following recent literature regarding the case of the whole space \cite{L2_5,L2_6,L2_2,L2_1,L2_3,L2_7,L2_8,L2_9,L2_4}.

\section{Optimal partition problems}\label{sec5}

Shape optimization problems are a class of problems where the general goal is to minimize (or maximize) a certain cost functional among a class  of  shapes, which are typically subsets of Euclidean spaces, manifolds or even metric graphs.  Two of the most famous examples (stated here in a slightly informal way) are the \emph{isoperimetric problem}:
\[
\min \{\text{per}(\omega):\ \omega\subset \R^N, \ |\omega|=a \}
\]
or the problem of finding the drum of a fixed $N$-volume that has the lowest fundamental frequency:
\[
\min\{\lambda_1(\omega):\  \omega\subset \R^N,\ |\omega|=a\},
\]
for a fixed $a>0$. Here, $per(\omega)$, $|\omega|$ and $\lambda_1(\omega)$ denote respectively the perimeter, the measure and the first Dirichlet eigenvalue of a set $\omega$.
In both situations the solution is a ball, as a consequence of the isoperimetric and the Rayleigh-Faber-Krahn inequalities, respectively.

In this section we focus on a subclass of shape optimization problems, namely on  the so called \emph{optimal partition problems}. Generally speaking, the aim is to study
\begin{equation}\label{eq:OPP}
\inf\left\{\Phi(\omega_1,\ldots, \omega_m):\ \omega_i\in \mathcal{A},\ \omega_i\cap \omega_j=\emptyset\ \forall i\neq j\right\},
\end{equation}
where $\mathcal{A}$ is a class of admissible sets in a certain ambient space and $\Phi:\mathcal{A}^m\to \R$ is a cost function. Observe that, here, the term \emph{partition} simply means that the shapes are disjoint; the condition that their union exhausts (in some sense) the whole domain is usually a consequence (\emph{a posteriori}) of the minimizing property of an optimal partition.

The problem of finding a partition that minimizes a certain cost function depending on disjoint shapes, despite its clear mathematical interest, appears quite naturally both in physics (e.g.  in liquid crystals or Cahn-Hilliard fluids \cite{AmbrosioBraides}), engineering (in situations where it is necessary to minimize the cost of a structure made of several materials) or image processing \cite{vanGoethem}. They are also important to characterize the limiting behavior of solutions to competing systems such as \eqref{S-system}, and play a fundamental role in the study of the nodal sets of eigenfunctions of Schr\"odinger operators \cite{BBHU,Berkolaiko,BHV,helffer,HHT,HHOT2,HHOT3,ShapeOptimizationBook}, as well as in the proof of monotonicity formulae \cite{acf,ContiTerraciniVerziniOPP,SoaveZilio}.

 In general, these kind of problems may only have a solution in a relaxed sense \cite{ButtazzoDalMaso, ButtazzoTimofte}, except when one imposes certain geometric constraints on the admissible domains, or some monotonicity properties on the cost function (we refer the reader to the book by Bucur and Buttazzo \cite{BucurButtazzoBook} for a good survey on these issues).

In this chapter we focus on our contributions to the field in four different classes of problems.

\subsection{Spectral optimal partitions}\label{sec:spectral}

Let $\Omega$ be a smooth bounded domain of $\R^N$, and let $k\geq 1$, $m \geq 2$ be integers. Consider the class of $m$ open partitions of $\Omega$:
\[
\mathcal{P}(\Omega)=\left\{  (\omega_1,\ldots, \omega_m)\left| \  \omega_i \subset \Omega \text{ is a nonempty open set for all $i$} ,\ \omega_i \cap \omega_j = \emptyset\  \forall  i\neq j\right.\right\}.
\]
We wish to solve
\begin{equation}\label{eq:OPPmain}
c_0=\inf \left\{ \sum_{i=1}^m \lambda_{k}(\omega_i):\  (\omega_1,\ldots, \omega_m)\in \mathcal{P}(\Omega) \right\}
\end{equation}
where $\lambda_k(\omega_i)$ is the $k$--th eigenvalue of $(-\Delta,H^1_0(\omega_i))$, counting multiplicities.\footnote{For simplicity, here we are considering the sum of eigenvalues, but other combinations are possible, see for instance \cite[p. 365]{RamosTavaresTerracini}.} The cost function
\[
\Phi(\omega_1,\ldots, \omega_m)=\sum_{i=1}^m \lambda_{k}(\omega_i)
\]
has good properties: it is monotone decreasing with respect to set inclusion, and it is lower semicontinuous for the $\gamma$--convergence. Therefore, the general abstract result of \cite{BucurButtazzoHenrot} implies the existence of a solution for \eqref{eq:OPPmain} in the class of \emph{quasi-open sets}. Going from quasi-open to open sets is, however, not an easy task. By using a penalization technique with partition of unity functions, Bourdin, Bucur and Oudet \cite{BourdinBucurOudet} gave a different proof for the existence of quasi-open solutions, while proving the existence of open solutions for the two-dimensional case $N=2$ (by using a compactness result \cite{Sverak} which only holds in dimension two).

The general goals are to determine the existence of a solution to this problem,  the optimal regularity of the associated eigenfunctions, and the regularity of the interfaces. Related with the latter, we make the following definition ($\mathscr{H}_{\text {dim }}(\cdot)$ denotes the Hausdorff dimension of a set):

\begin{definition}[{{\cite{RamosTavaresTerracini}}}]\label{def:regular_partitions} An open partition $\left(\omega_1, \ldots, \omega_m\right) \in \mathcal{P}(\Omega)$ is called regular if:
\begin{enumerate}
\item denoting $\Gamma=\Omega \backslash \bigcup_{i=1}^m \omega_i$, it holds $\mathscr{H}_{\text {dim }}(\Gamma) \leq N-1$;
\item there exists a relatively open  subset $\mathcal{R} \subseteq \Gamma$, such that
\begin{itemize}
\item $\mathscr{H}_{\operatorname{dim}}(\Gamma \backslash \mathcal{R}) \leqq N-2$;
\item  $\mathcal{R}$ is a collection of hypersurfaces of class $C^{1, \alpha}$ (for some $0<\alpha<1$ ), each one separating two different elements of the partition.
\end{itemize}
\end{enumerate}
\end{definition}

\begin{figure}
    \centering
    \includegraphics[scale=0.35]{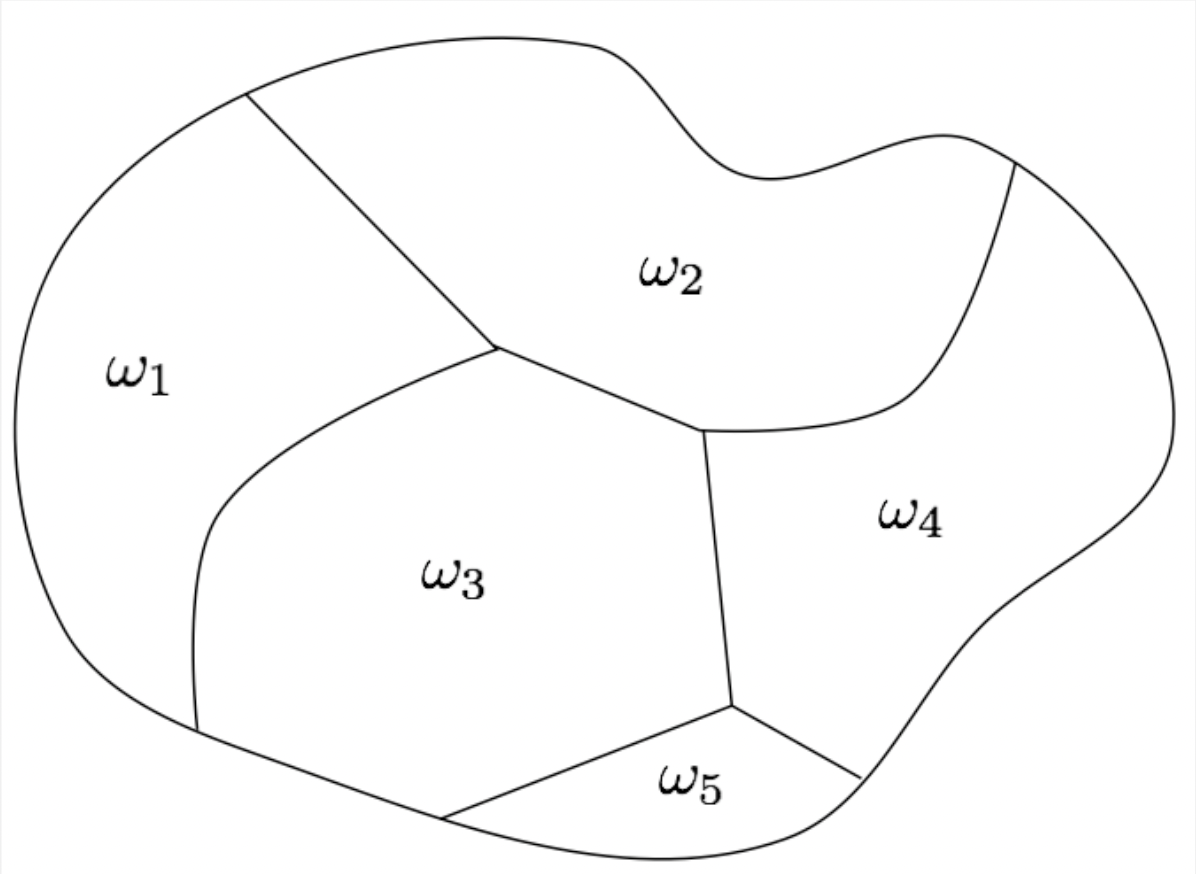}
    \caption{Example of a regular 5-partition in dimension $N=2$.}
\end{figure}

The study of \eqref{eq:OPPmain} in the case $k=1$ is the simplest, as it can be characterized by an absolute minimization of an energy functional in a singular space, namely
\begin{equation}\label{eq:equivalent_characterization_of_sum_lambda1}
c_0=\tilde c_0:=\inf \left\{\sum_{i=1}^m \int_\Omega |\nabla u_i|^2\, dx:\begin{array}{c}\displaystyle u_i\in H^1_0(\Omega)\text{ and } \int_\Omega u_i^2\, dx=1 \ \forall i,\\[8pt]
 u_i\cdot u_j\equiv 0 \ \forall i\neq j \end{array}\right\}.
\end{equation}

This minimization problem was the main object of study in \cite{CaffarelliLin2,ContiTerraciniVerziniOPP}, where it is proved that a nonnegative Lipschitz continuous solution $(u_1,\ldots, u_m)$ to \eqref{eq:equivalent_characterization_of_sum_lambda1} exists, and that the partition $(\{u_1>0\},\ldots, \{u_m>0\})$ is a regular element of $\mathcal{P}(\Omega)$, achieving $\inf_{\mathcal{P}(\Omega)} \sum_{i=1}^m \lambda_1(\omega_i)$ (\emph{i.e.}, it is an optimal partition). The existence result was first proved in \cite{ContiTerraciniVerziniOPP} together with the regularity for $N=2$. The regularity result in any space dimension was then stated in \cite{CaffarelliLin2}, see also \cite[Section 8]{TavaresTerracini1} for a detailed proof (the latter is a joint work with S. Terracini). In my opinion,
 the simplest approach to this situation nowadays is to consider limiting profiles of solutions to the singularly perturbed problem
%
\begin{equation}\label{eq:BEC}
-\Delta u_i=\lambda_{i,\beta} u_i +\beta \mathop{\sum_{j=1}^m}_{j\neq i} u_iu_j^2,\quad u_i\in H^1_0(\Omega), \qquad i=1,\ldots, m,
\end{equation}
under the constraints
\[
\int_\Omega u_i^2\, dx=1,\qquad i=1,\ldots, m,
\]
so that the parameters $\lambda_{i,\beta}$ appear as Lagrange multipliers.

These systems (similar to the ones that have appeared before in Section \ref{sec4}) have been the object of an intensive study in the last fifteen years, in particular in the case of competitive interaction $\beta<0$ and the study of the singular limit $\beta\to -\infty$. Their relation with optimal partition problems has also been addressed, for instance in \cite{BTWW,CaffarelliLin2,clll,ContiTerraciniVerziniOPP,HHT,TavaresTerracini2}. We have shown in \cite{NTTV1, TavaresTerracini2} that, in some situations, phase separation occurs between different components as the competition parameter increases, \emph{i.e.}, $\beta\to -\infty$. In particular it is shown that, by taking an $L^\infty$ bounded family of solutions $(u_\beta)_\beta$, and corresponding bounded coefficients $(\lambda_{i,\beta})_\beta$, then there exists a limiting profile $u_i:=\lim_{\beta\to +\infty} u_{i\beta}$ such that $(\{u_1\neq 0\},\ldots, \{u_m\neq 0\})\in \mathcal{P}(\Omega)$, and
\[
-\Delta \tilde u_i=\lambda_i \tilde u_i \qquad \text{ in } \{u_i\neq 0\}.
\]
This clearly illustrates the relation between optimal partitions involving eigenvalues and the system of Schr\"odinger equations \eqref{eq:BEC}. In particular, it is known that \eqref{eq:equivalent_characterization_of_sum_lambda1} can be well approximated (as $\beta\to -\infty$) by the ground state (least energy) levels of \eqref{eq:BEC}, namely:
\[
\inf \left\{\int_\Omega \sum_{i=1}^m |\nabla u_i|^2 - \beta \sum_{i<j} u_i^2 u_j^2:\ u_i\in H^1_0(\Omega) \text{ and } \int_\Omega u_i^2\, dx=1\ \forall i \right\}.
\]
Thus, using this approach and the results from a joint paper with S. Terracini \cite{TavaresTerracini1}, one proves once again the existence of a regular partition to the problem of summing first eigenvalues. However, passing to higher eigenvalues is not an easy task, as one needs to construct suitable minimax characterizations at higher energy levels of \eqref{eq:BEC}. In another paper with S. Terracini \cite{TavaresTerracini2}, by using a new notion of vector genus, several sign changing solutions are build for \eqref{eq:BEC}, and by taking the least energy nodal solution among these, one approaches the second eigenfunctions associated with the optimal partition problem $\inf_{\mathcal{P}(\Omega)} \sum_{i=1}^m\lambda_2(\omega_i)$. By putting together the previous results, one then can actually solve \eqref{eq:OPPmain} for a combination of sums of first and second eigenvalues. In order to solve the general problem with higher eigenvalues, however, it does not seem completely clear to us which variational characterization for solutions of \eqref{eq:BEC} one could take.

\smallbreak

To finally solve the general case, in \cite{RamosTavaresTerracini} (a joint work with M. Ramos and S. Terracini) we followed instead  a different strategy relying on a double approximation procedure, which I describe next. The relevant and surprising fact is that, instead of taking minimax levels for a certain energy functional, we are able to approximate the problem \eqref{eq:OPPmain} for every $k\in\N$ through a symmetric constrained energy minimization. The strategy was influential also in later papers and other contexts, see for instance \cite{KriventsovLin, MazzoleniTerraciniVelichkov, MazzoleniTreyVelichkov}.

In order to cope with the problem of not knowing the multiplicity of each set of the optimal partition \emph{a priori}, our motivation was to try to find approximate solutions of \eqref{eq:OPPmain} through the minimization process of a certain energy functional. Partially inspired by \cite{HHT} (where a different problem is treated), we let $p \in \mathbb{N}$ and consider the problem
\begin{equation}\label{eq:OPPapprox_p}
\inf _{\left(\omega_1, \ldots, \omega_m\right) \in \mathcal{P}_m(\Omega)} \sum_{i=1}^m\left(\sum_{j=1}^{k}\left(\lambda_j\left(\omega_i\right)\right)^p\right)^{1 / p} .
\end{equation}
Observe that \eqref{eq:OPPapprox_p} is a reasonably good approximation for \eqref{eq:OPPmain} for large $p$, as, given $k \in \mathbb{N}$ and any positive real numbers $a_1, \ldots, a_k$, there holds $\left(a_1^p+\cdots+a_k^p\right)^{1 / p} \rightarrow$ $\max \left\{a_1, \ldots, a_k\right\}$ as $p \rightarrow \infty$. Thus, for any given partition $\left(\omega_1, \ldots, \omega_m\right) \in \mathcal{P}(\Omega)$,
$$
\sum_{i=1}^m\left(\sum_{j=1}^{k} \lambda_j\left(\omega_i\right)^p\right)^{1 / p} \to \sum_{i=1}^m \lambda_{k}\left(\omega_i\right) \quad \text { as } \quad p \to+\infty .
$$
We proved that an optimal solution $\left(\omega_{1 p}, \ldots, \omega_{m p}\right)$ of \eqref{eq:OPPapprox_p} exists and approaches, as $p \rightarrow \infty$, a solution of our original problem \eqref{eq:OPPmain}. To show the latter, we approximated \eqref{eq:OPPapprox_p} by a system of type \eqref{eq:BEC} where the competition occurs between groups (in the spirit of Section \ref{sec4}).

The final result is the following.

\begin{theorem}[{{\cite{RamosTavaresTerracini}}}]
Given $k\geq 1$, there exists a regular optimal partition $(\omega_1,\ldots, \omega_m)$ for the problem \eqref{eq:OPPmain}. Moreover, for each $i$, at least one eigenfunction $u_{i}$ associated with $\lambda_k(\omega_i)$ is globally Lipschitz continuous, which is the optimal regularity in this case.
\end{theorem}
We observe that a free boundary condition was also proved, involving some $k$--eigenfunctions in a neighborhood of regular points of the interfaces; however, the actual statement is too technical to include in this work.

For $k=1$, finer results for the singular set are proved in the recent paper \cite{Alper}, namely that the $(N-2)$--Hausdorff dimension of the singular set is finite, together with a stratification result. In \cite{TavaresZilio}, together with A. Zilio, we characterized and proved regularity of all possible minimal partitions of problem like \eqref{eq:OPPapprox_p} (which involve combinations of eigenfunctions up to a certain order) and their eigenfunctions. On the other hand, it remains an open problem to prove the same for the original problem \eqref{eq:OPPmain} with $k\geq 2$.

\subsection{Long-range  spectral optimal partitions}

In this subsection we discuss a related class of optimal partition problems. Instead of considering classes of partitions where the sets are simply mutually disjoint, we introduce a restriction about the distance between sets. Given $r> 0$, consider the set of all $m$-partitions of $\Omega$ whose elements are at distance at least $r$:
\[
\mathcal{P}_r(\Omega)=\left\{  (\omega_1,\ldots, \omega_m)\left| \begin{array}{l} \omega_i \subset \Omega \text{ is a nonempty open set for all $i$} ,\\  \dist(\omega_i,\omega_j)\geq r\ \quad  \forall i\neq j\end{array}\right.\right\}.
\]
 It is straightforward that there exists $\bar r>0$ (which depends on $\Omega$ and on $m$) such that $\mathcal{P}_r(\Omega) \neq \emptyset$, for every $r\in [0,\bar r)$. For any such $r$, we are concerned with the following optimization problem:
\begin{equation}\label{eqn eig dist}
c_r:=\inf \left\{ \sum_{i=1}^m \lambda_1(\omega_i):\  (\omega_1,\ldots, \omega_m) \in \mathcal{P}_r(\Omega) \right\},
\end{equation}
where $\lambda_1(\cdot)$ denotes the first Dirichlet eigenvalue. In a joint paper with N. Soave, S. Terracini and A. Zilio \cite{STTZ2}, we have proved the following:
\begin{theorem}[{{\cite{STTZ2}}}]\
 \begin{enumerate}
\item \textit{Existence.} The level $c_r$ is achieved by an open optimal partition $(\Omega_{1,r},\ldots, \Omega_{m,r})$;
\item \textit{Regularity of Eigenfunctions.} If $u_{i,r}$ is a first eigenfunction associated with $\Omega_{i,r}$, then it is globally Lipschitz continuous.
\item \textit{Exterior sphere condition and exact distance between the optimal sets.} Given $x_0\in \partial \Omega_{i,r} \setminus \partial \Omega$, there exists $j\neq i$ and $y_0\in \partial \Omega_{j,r}$ such that $|x_0-y_0|=r$, and $\Omega_{i,r} \cap B_r(y_0) = \emptyset$; in particular, $\dist(\Omega_{i,r},\Omega_{j,r})=r$ and each set $\Omega_{i,r}$ satisfies an exterior sphere condition of radius $r$ at any of its boundary point.
\end{enumerate}
\end{theorem}
The latter is a statement which is specific of these long range optimal partition problems. This statement, together with \cite[Lemma 6.4]{CaffarelliPatriziQuitalo}, yields the following important information about the free boundary $\Gamma_r:=\cup_{i=1}^m\partial\Omega_{i,r}$:
\begin{enumerate}
\item[4.] \textit{Measure of the Free Boundary.} The sets $\partial \Omega_{i,r}$ have locally finite perimeter in $\Omega$.
\end{enumerate}
Under an additional regularity assumption on the free boundary $\partial \Omega_{i,r}$, we have also derived a free boundary condition satisfied by the eigenfunctions of the optimal partitions (see \cite[Theorem 1.6]{STTZ2}).  The validity of such regularity remains a crucial open problem in the general setting for optimal partition problems with a distance constraint. The difficulty is related with the construction  of admissible variations.

The approach used in \cite{STTZ2} consists of studying the following relaxed formulation of $c_r$ in terms of measurable functions rather than sets:
\begin{equation}\label{eq:weak_characterization}
\tilde c_r=\inf\left\{ \sum_{i=1}^m \int_\Omega |\nabla u_i|^2\left| \begin{array}{l} u_i \in H^1_0(\Omega), \ \int_\Omega u_i^2=1\ \forall i, \\  \dist(\text{supp}\, u_i,\text{supp}\, u_j)\geq r,\ \forall i\neq j \end{array}\right. \right\},
\end{equation}
proving the equivalence between $c_r$ and $\tilde c_{r}$. We also show a relation between \eqref{eq:weak_characterization} and an elliptic system with nonlocal competition terms
\begin{equation}\label{syst p}
-\Delta u_{i,\beta}=\lambda_{i,\beta} u_{i,\beta}+\beta u_{i,\beta}(x) \sum_{j\neq i} \int_{B_r(x)} V_r(x-y) u_{j,\beta}^2(y)\, dy
\end{equation}
where $V_r\in L^\infty(\R^N)$ satisfies $V_r>0$ a.e. in $B_r$(0), $V_r=0$ a.e. on $\R^N\setminus \overline{B_r}(0)$, and $\beta<0$. The only other results available so far regarding segregation problems driven by long-range competition are given in \cite{CaffarelliPatriziQuitalo}, where the authors analyze the spatial segregation for systems of type
\begin{equation}\label{eq:nonlocal_Caffa}
-\Delta u_{i, \beta}=\beta u_{i, \beta} \sum_{j \neq i}\left(\mathbbm{1}_{B_r} \star\left|u_j\right|^p\right) \text { in } \Omega,\qquad u_{i, \beta}=f_i \geq 0  \text { in } \Omega_r \backslash \Omega,
\end{equation}
with $1 \leqslant p \leqslant+\infty$, as $\beta\to -\infty$. In the above equation, $\mathbbm{1}_{B_r}$ denotes the characteristic function of $B_r$, the ball of center 0 and radius $r$, $\Omega_r$ is the neighborhood of radius $r$ of $\Omega$, and $\star$ stands for the convolution for $p<+\infty$, so that
$$
\left(\mathbbm{1}_{B_r} \star\left|u_j\right|^p\right)(x)=\int_{B_r(x)}\left|u_j(y)\right|^p d y \quad \forall x \in \Omega, \text { with } 1 \leqslant p<+\infty;
$$
in case $p=+\infty$, it is meant that the integral should be replaced by the supremum over $B_r(x)$ of $|u_j|$. In \cite{CaffarelliPatriziQuitalo}, the authors prove the equicontinuity of families of viscosity solutions $\left\{\mathbf{u}_\beta: \beta<0\right\}$ to \eqref{eq:nonlocal_Caffa}, the local uniform convergence to a limit configuration $\mathbf{u}$, and then study the free-boundary regularity of the positivity sets $\left\{u_i>0\right\}$ in cases $p=1$ and $p=+\infty$, mostly in dimension $N=2$.


\medbreak

The techniques adopted in the local and nonlocal cases are completely different. Powerful tools typically employed in the former ones, such as monotonicity formulas, free boundary conditions and blow-up methods, cannot be adapted in the context of optimal partitions at distance, due to the nonlocal nature of the interaction between different densities/sets. This is why the free boundary regularity for problem \eqref{eq:OPPmain} is settled, while the same problem for \eqref{eqn eig dist} is open. However, the common optimal Lipschitz regularity of $\mf{u}_r$ suggests that it should be possible to look at both problems, the local and the nonlocal ones, as a $1$-parameter family, where the parameter is the distance $r$ between the different supports. The main results of a joint paper with N. Soave and A. Zilio \cite{SoaveTavaresZilio} establish that this is possible, at least at the level of the eigenfunctions. More precisely:

\begin{theorem}[{{\cite{SoaveTavaresZilio}}}]
There exists a constant $C > 0$ such that
 \[
	\|\mathbf{u}_r\|_{\mathrm{Lip}(\overline \Omega)}:=\|\mathbf{u}_r\|_{L^\infty(\Omega)}+\| \nabla \mathbf{u}_r\|_{L^\infty(\Omega)} \leq C,
	\]
for any $0< r < \bar r$, and any minimizer $\mathbf{u}_r$ of $\tilde c_r=c_r$.
\end{theorem}
Observe that, for each fixed $r>0$, Lipschitz regularity may be proved via a barrier argument, which is possible due to the exterior sphere condition (see \cite[Theorem 3.4]{STTZsurvey}). However, the barrier used depends on the radius, and the argument breaks down as $r\to 0^+$. In \cite{SoaveTavaresZilio} we relied on different methods.

Combining this theorem with the information obtained in previous papers (and described in the previous subsection) about the local case $r=0$, we have the following (denoting the level in \eqref{eq:OPPmain} for  $k=1$ by $c_0$)
\begin{theorem}[{{\cite{SoaveTavaresZilio}}}]
There exists $C>0$ such that
\[
	c_0 \leq c_r \leq c_0 + C r \quad \text{ for sufficiently small $r>0$.}
	\]
In particular, $c_r\to c_0$ as $r\to 0$. Moreover, given any minimizer $\mathbf{u}_r$ of $c_r$ for $r>0$,  there exists $\mathbf{u}_0\in H^1_0(\Omega)\cap \mathrm{Lip}(\overline \Omega)$, solution to $c_0$, such that, up to a subsequence,
\[
\mathbf{u}_r\to \mathbf{u}_0 \quad \text{ strongly in } H^1_0(\Omega)\cap C^{0,\alpha}(\overline \Omega), \text{ for every } \alpha\in (0,1).
\]
\end{theorem}

In this way, we are establishing a relation between problems \eqref{eq:OPPmain} and \eqref{eqn eig dist}. We believe that these results may pave the way towards the development of a common free boundary regularity theory. In particular, we wonder if the very complete information known for the free boundary in the limiting problem \eqref{eq:OPPmain} can be used to deduce properties for the free boundary arising in \eqref{eqn eig dist}, at least for a small $r$.

\subsection{Spectral partition problems with volume constraint}

While the literature is full of examples of shape optimization problems with volume constraints (see for instance \cite[Chapters 2 and 3]{ShapeOptimizationBook} and references therein), up to our knowledge no one considered \emph{optimal partition problems} with volume constraints. These can be easily motivated, considering for instance the situation of a farmer that is planting several crops in a region $\Omega$, while on the other hand there is a legal limit on the amount of land that can be used for agriculture. This leads to a problem as follows:
\begin{equation}\label{eq:OPP2}
\inf\left\{\Phi(\omega_1,\ldots, \omega_m):\ \omega_i\in \mathcal{A}\ \forall i,\ \omega_i\cap \omega_j=\emptyset \ \forall i\neq j,\ \sum_{i=1}^m |\omega_i|\leq a \right\}.
\end{equation}

Given a bounded domain $\Omega \subset \mathbb{R}^N$ and $0<a<|\Omega|$, we consider the following prototypical model problem
\begin{equation}\label{eqn eig dist2}
c^a:=\inf \left\{ \sum_{i=1}^m \lambda_1(\omega_i):\  (\omega_1,\ldots, \omega_m) \in \mathcal{P}^a(\Omega) \right\},
\end{equation}
where ${\mathcal P}^a(\Omega)$ stands for the set of $m$-partitions of $\Omega$ with volume constraint $a$, \emph{i.e.,}
\begin{equation*}
{\mathcal P}^a(\Omega):= \left\{(\omega_1, \ldots, \omega_m)\;\Big|\;
\begin{array}{c}
\omega_i \subset \Omega \mbox{ are nonempty open sets for all } i, \vspace{0.05cm}\\
 \omega_i \cap \omega_j = \emptyset  \: \text{for all}\: i \not=j \mbox{ and } \sum_{i=1}^{m}|\omega_i| \leq  a\\
\end{array}
\right\}.
\end{equation*}
Here, $|\cdot|$ stands for the Lebesgue measure.

In order to investigate the problem \eqref{eqn eig dist2},  we introduce a weak formulation that involves a minimization problem where the variables are functions rather than domains, namely
\begin{equation}\label{main_func}
\tilde{c}^a = \inf_{(u_1, \ldots, u_m) \in H^a} J(u_1, \ldots, u_m),\quad \text{ where } J(u_1, \ldots, u_m) :=\sum_{i=1}^m \int_{\Omega} |\nabla u_i|^2
\end{equation}
and
\begin{equation*}
H^a : = \left\{(u_1, \ldots,  u_m)  : \begin{array}{c}\displaystyle
  u_i \in H^1_0(\Omega) \mbox{ and }\int_{\Omega} u_i^2 = 1 \mbox{ for every $i$}, \\[8pt]
\displaystyle u_i u_j \equiv 0 \mbox{ for $i\neq j$},\ \  \sum_{i =1}^m|\Omega_{u_i}|\leq  a
 \end{array}
\right \},
\end{equation*}
with $\Omega_{u_i} := \{ u_i\neq 0\}$ for all $i= 1, \ldots, m$.  This is a joint project with P. Andrade, E. Moreira dos Santos and M. Santos \cite{AndradeSantosSantosTavares}, whose main result reads as follows:
\begin{theorem}[{{\cite{AndradeSantosSantosTavares}}}]
The problem \eqref{eqn eig dist2} admits a solution. Moreover:
\begin{enumerate}
\item Given  any optimal partition $(\omega_1,\ldots, \omega_m)\in \mathcal{P}^a(\Omega)$, we have that
\[
\text{each } \Omega_i \text{ is connected and } \sum_{i=1}^k |\Omega_i|=a.
\] If $u_i$ is a first eigenfunction associated with the set $\Omega_i$, we have that $u_i$ is locally Lipschitz continuous in $\Omega$.
\item Problems \eqref{eqn eig dist2} and \eqref{main_func} are equivalent in the following sense:\begin{itemize}
\item  $c^a=\tilde{c}^a$;
\item if $(u_1,\ldots, u_m)\in H$ is an optimal solution of \eqref{main_func} and $\Omega_{u_i}:=\{u_i\neq 0\}$, then  $(\Omega_{u_1},\ldots, \Omega_{u_m})\in {\mathcal P}^a(\Omega)$ solves \eqref{eqn eig dist2};
\item  if $(\omega_1,\ldots, \omega_k)\in \mathcal{P}^a(\Omega)$ is an optimal partition for \eqref{eqn eig dist2} and $u_i$ is a first eigenfunction associated with the set $\Omega_i$, then $(u_1,\ldots, u_m)\in H^a$ is a minimizer of \eqref{main_func}.
\end{itemize}
\end{enumerate}
\end{theorem}
 An important part of the proof is based on showing the equivalence with a minimization problem for the following penalized functional:
\begin{equation*}
J_{\mu} (u_1,\ldots, u_k):= \sum_{i=1}^k \frac{\displaystyle \int_{\Omega}|\nabla u_i|^2}{\displaystyle \int_{\Omega} u_i^2}  + \mu \left[ \sum_{i=1}^k |{\Omega}_{u_i}| - a\right ]^{+}, \qquad \text{for} \ \ (u_1\ldots, u_k) \in \overline{H},
\end{equation*}
where
\begin{equation*}
\overline{H} : = \Big\{(u_1,\ldots,  u_k) \in H^1_0(\Omega;\R^k)\:  \Big|\: \quad  u_i \neq 0\ \forall i,\   u_i \cdot u_j \equiv 0\ \ \forall i\neq j \Big \}.
\end{equation*}

	The main difficulty when leading with this problem is mostly related with the production of admissible variations  which, on the other hand, give relevant information to the problem. We would also like to point out that, even though the sets $\Omega_{u_i}$ are quasi-open, we do not use this fact directly in our paper, nor the concept of $\gamma$-convergence of quasi-open set is used. Instead, the proof follows nontrivial adaptations of ideas from \cite{Brianconetal,LamboleyPierre} (shape optimization with measure constraints with one set only, no partitions) and \cite{ContiTerraciniVerziniOPP,HHT} (partition problem, no measure constraint).  The study of the regularity of the free boundary and the production of numerical simulations are the subject of current work.

We conclude by mentioning the following related (although not equivalent) problems regarding optimal partitions \cite{BogoselVelichkov,BucurVelichkov,PhilippisSpolaorVelichkov}, where the cost function is
\[
\Phi(\omega_1,\ldots, \omega_m)=\sum_{i= 1}^{k} (\lambda_{\ell}(\omega_i) + m |\omega_i|)
\] defined on partitions. This problem does not have measure constraints, although for a large $m$ there the optimal configurations will not occupy the whole $\Omega$. We refer to our introduction in \cite{AndradeSantosSantosTavares} for more details.

\subsection{Optimal partitions related with the Yamabe equation}

Several papers over the years refer to optimal partitions with nonlinear costs. This is related with the study of nodal solutions of single equations. We refer for instance  to \cite{ContiTerraciniVerziniNehari,ContiTerraciniVerziniOptimalNonlinear}, where the optimal cost in \eqref{eq:OPP} is
\begin{equation}\label{eq:nonlinear_cost}
\Phi(\omega_1,\ldots, \omega_k):=\sum_{i=1}^k c(\omega_i),
\end{equation}
$c(\omega_i)$ being the least energy solution of equations of the form $-\Delta =f(u)$ with subcritical--superlinear growth, and homogeneous Dirichlet boundary conditions in a bounded domain. As a prototypical example, the authors take
\begin{equation}\label{eq:singleeq_OPP}
-\Delta u+\lambda u= \mu |u|^{p-1}u \text{ in } \omega,\qquad u=0 \text{ on } \partial \omega
\end{equation}
with $\lambda>-\lambda_1(\Omega)$ and $1<p<2^*-1=(N+2)/(N-2)^+$ in the focusing case $\mu>0$ (the defocusing case is considered in \cite{clll}). On the other hand, the Sobolev critical case $p=2^*-1$, $\lambda\in (-\lambda_1(\Omega),0)$ is tackled by us in \cite{TavaresYou, TavaresYouZou} (recall Subsection \ref{sec:criticalcase_gradsystem}). In these papers, the existence of optimal partitions is proved, while its regularity (in the sense of Definition \ref{def:regular_partitions}) is determined when combining these references with my joint paper with S. Terracini \cite{TavaresTerracini1}. As an important application, we refer that in the special case of $m=2$ partition problems, one finds in this way a least energy nodal solution of \eqref{eq:singleeq_OPP} (recall, for instance, Theorem \ref{Sign-changing solutions}). We also refer to the recent paper \cite{ClappFernandezSaldana}, where optimal partition problems related with polyharmonic semilinear equations are considered.

Some of the things that have been mentioned  can be adapted to the context of optimal partition problems on Riemannian manifolds. One of the most relevant related problems is the study of the Yamabe equation, which has an interesting and fascinating history. It is by now classical to show that, answering the question
\begin{quote}
Given $(\mathcal{M},g)$, a closed Riemannian manifold of dimension $N\geq 3$ with metrig $g$,  is there a conformal metric with constant scalar curvature?
\end{quote}
amounts to finding positive smooth solutions to the Yamabe equation:
\begin{equation}\label{eq:Yamabe}
-\Delta_g u + \kappa_N S_gu = |u|^{2^*-2}u\qquad\text{on }\mathcal{M},
\end{equation}
where $S_g$ is the scalar curvature, $\Delta_g:=\mathrm{div}_g\nabla_g$ the Laplace-Beltrami operator, $\kappa_N:=\frac{N-2}{4(N-1)}$. In this case, the existence of a positive solution was established thanks to the combined efforts of Yamabe \cite{Yamabe}, Trudinger \cite{Trudinger}, Aubin \cite{Aubin1} and Schoen \cite{Schoen}. A detailed account is given in \cite{LeeParker}.

It is natural to consider an optimal partition problem of type \eqref{eq:OPP} with cost \eqref{eq:nonlinear_cost}, where this time $c(\omega)$ represents a least energy solution to the Yamabe equation \eqref{eq:Yamabe} in $\omega\subset M$, with homogeneous Dirichlet boundary conditions on $\partial \omega$. However, it is important to remark that optimal partitions do not always exist! In fact, there is no optimal $m$-partition for the Yamabe equation on the standard sphere $\mathbb{S}^N$ for any $m\geq 2$.  In \cite{ClappPistoiaTavares}, jointly with M. Clapp and A. Pistoia, we proved the following.

\begin{theorem}[{{\cite{ClappPistoiaTavares}}}]\label{thm:Yamabe} Assume that
$(\mathcal{M},g)$ is not locally conformally flat and $\dim \mathcal{M}\geq 10$. If $\dim \mathcal{M}=10$, assume furthermore that
\begin{equation}\label{eq:A3}
|S_g(q)|^2<\frac{5}{28}\,|W_g(q)|^2_g\qquad\forall q\in \mathcal{M},
\end{equation}
where $W_g(q)$ is the Weyl tensor of $(\mathcal{M},g)$ at $q$.

 Then, for every $m\geq 2$ there exists an optimal $m$-partition $\{\omega_1,\ldots,\omega_m\}$ for the Yamabe equation on $(\mathcal{M},g)$, such that each $\omega_i$ is connected and is regular in the sense of Definition \ref{def:regular_partitions}.

 In particular, for $m=2$, this yields the existence of a least energy nodal solution to the Yamabe equation \eqref{eq:Yamabe} having precisely two nodal domains.
\end{theorem}

We follow the approach of considering a singular perturbation, \emph{i.e.}, to approximate the problem with a system of Yamabe-equations joined by a variational competition term (in the spirit of Section \ref{sec:spectral}), proving the existence of fully nontrivial least energy solutions and studying the behavior as the competition coefficient diverges. To prove this and to prevent blowup, a new compactness criterion  is established.  To verify this criterion, we introduce a test function and perform rather delicate estimates (inspired by fine estimates established in \cite{EspositoPistoiaVetois,LeeParker}), particularly in dimension $10$ - where not only the exponents but also the coefficients of the energy expansion play a role - leading to the geometric inequality stated in assumption \eqref{eq:A3}.

In order to prove the optimal regularity of the limiting profiles $u_i$, the regularity of the free boundaries $\mathcal{M}\smallsetminus\bigcup_{i=1}^m\omega_i$ and the free boundary condition, we use local coordinates. This reduces the problem to the study of segregated profiles satisfying a system involving divergence type operators with variable coefficients. We are able to prove  \emph{a priori} bounds in H\"older spaces, by deducing an Almgren-type monotonicity formulae  and by performing a blowup analysis, combining what is known in case of the pure Laplacian \cite{CaffarelliLin,NTTV1, TavaresTerracini1,STTZsurvey} with some ideas from papers dealing with variable coefficient operators \cite{Kukavica,Mariana1,Mariana2,SoaveWeth}.  We remark that, in a recent work with M. Dias \cite{DiasTavares},  we were able to obtain uniform Lipschitz bounds, which are the optimal uniform estimates in this context.

The existence of nodal solutions to the Yamabe equation \eqref{eq:Yamabe} on an arbitrary manifold $(\mathcal{M},g)$ is largely an open problem. In \cite{AmmannHumbert}, the existence of a least energy nodal solution is established  when $(\mathcal{M},g)$ is not locally conformally flat and $\dim \mathcal{M}\geq 11$. Theorem \eqref{thm:Yamabe} recovers and extends this result. We also note that an optimal $m$-partition $\{\omega_1,\ldots,\omega_m\}$ gives rise to what in \cite{AmmannHumbert} is called \emph{a generalized metric $\bar g:=\bar u^{2^*-2}g$ conformal to $g$} by taking $\bar u:=u_1+\cdots+u_m$ with $u_i$ a positive solution to \eqref{eq:Yamabe} in $\omega_i$. So Theorem \ref{thm:Yamabe} may be seen as an extension of the main result in \cite{AmmannHumbert}.

As we mentioned before, optimal $m$-partitions on the standard sphere $\mathbb{S}^N$ do not exist. However, if one considers partitions with the additional property that every set $\omega_i$ is invariant under the action of a suitable group of isometries, then optimal $m$-partitions of this kind do exist and they give rise to sign-changing solutions to the Yamabe equation \eqref{eq:Yamabe} with precisely $m$-nodal domains for every $m\geq 2$, as shown in \cite{ClappSaldanaSzulkin_sphere}. The case of a general manifold $\mathcal{M}$ possessing some symmetries is treated in \cite{ClappPistoia_symmetries}.

\medskip

\section*{Acknowledgments}

This paper is based on my \emph{Habilitation lecture} in Mathematics, presented to Instituto Superior T\'ecnico, Universidade de Lisboa, on May 22-23, 2023. The lecture describes part of my previous work, which would not have been possible without my collaborators.  I would like to thank all of them for their lessons, ideas, support and, in many cases, for their friendship: Nicola Abatangelo, Francisco Agostinho, P\^edra Andrade, Denis Bonheure, Daniele Cassani, M\'onica Clapp, Sim\~ao Correia,  Manuel Dias, Jo\~ao Paulo Dias, Juraj F\"oldes, Massimo Grossi, Ederson Moreira dos Santos, Benedetta Noris, Gabrielle Nornberg, Filipe Oliveira, Enea Parini, Angela Pistoia, Miguel Ramos, Jos\'e Francisco Rodrigues, Alberto Salda\~na, Makson Santos, Delia Schiera, Nicola Soave, Susanna Terracini, Gianmaria Verzini, Tobias Weth, Song You, Jianjun Zhang, Alessandro Zilio, Wenming Zou.

Over the years I've also been supported in some periods by the Calouste Gulbenkian Foundation (Portugal), and supported continuously by the Portuguese government through FCT-Funda\c c\~ao para a Ci\^encia e a Tecnologia (at the time of writing, FCT is supporting my work through the projects UID/MAT/04459/2020 and PTDC/MAT-PUR/1788/2020).

I would also like to acknowledge the research centers CMAFcIO and CAMGSD.


\phantomsection

\addcontentsline{toc}{section}{References}
\makeatletter
\newcommand{\adjustmybblparameters}{\setlength{\itemsep}{2pt}\setlength{\parsep}{0pt}}
\let\ORIGINALlatex@openbib@code=\@openbib@code
\renewcommand{\@openbib@code}{\ORIGINALlatex@openbib@code\adjustmybblparameters}
\makeatother




\begin{thebibliography}{}

\end{thebibliography}


\begin{thebibliography}{100}

\bibitem{AbatangeloSaldanaTavares}
Nicola Abatangelo, Alberto Salda\~{n}a, and Hugo Tavares.
\newblock An asymptotic relationship between Lane-Emden systems and the
  1-bilaplacian equation, arXiv:2312.16696
\newblock

\bibitem{AdimurthiMancini}
Adi Adimurthi and Gabriele Mancini.
\newblock The {N}eumann problem for elliptic equations with critical
  nonlinearity.
\newblock In {\em Nonlinear analysis}, Sc. Norm. Super. di Pisa Quaderni, pages
  9--25. Scuola Norm. Sup., Pisa, 1991.

\bibitem{AftalionPacella}
Amandine Aftalion and Filomena Pacella.
\newblock Qualitative properties of nodal solutions of semilinear elliptic
  equations in radially symmetric domains.
\newblock {\em C. R. Math. Acad. Sci. Paris}, 339(5):339--344, 2004.

\bibitem{AgostinhoCorreiaTavares}
Francisco Agostinho, Sim\~{a}o Correia, and Hugo Tavares.
\newblock Classification and stability of positive solutions to the nls
  equation on the $\mathcal{T}$--metric graph.
\newblock {\em Nonlinearity}, 4(2), 2024.

\bibitem{AkAn}
Nail Akhmediev and Adrian Ankiewicz.
\newblock Partially coherent solitons on a finite background.
\newblock {\em Phys. Rev. Lett.}, 82:2661, 1999.

\bibitem{Alper}
Onur Alper.
\newblock On the singular set of free interface in an optimal partition
  problem.
\newblock {\em Comm. Pure Appl. Math.}, 73(4):855--915, 2020.

\bibitem{acf}
Hans Alt, Luis Caffarelli, and Avner Friedman.
\newblock Variational problems with two phases and their free boundaries.
\newblock {\em Trans. Amer. Math. Soc.}, 282(2):431--461, 1984.

\bibitem{AlvesSoares}
Claudianor~O. Alves and S{{\'e}}rgio H.~M. Soares.
\newblock Singularly perturbed elliptic systems.
\newblock {\em Nonlinear Anal.}, 64(1):109--129, 2006.

\bibitem{AmbrosettiColorado}
Antonio Ambrosetti and Eduardo Colorado.
\newblock Standing waves of some coupled nonlinear {S}chr{\"o}dinger equations.
\newblock {\em J. Lond. Math. Soc. (2)}, 75(1):67--82, 2007.

\bibitem{AmbrosettiProdiBook}
Antonio Ambrosetti and Giovanni Prodi.
\newblock {\em A primer of nonlinear analysis}, volume~34 of {\em Cambridge
  Studies in Advanced Mathematics}.
\newblock Cambridge University Press, Cambridge, 1993.

\bibitem{AmbrosettiRabinowitz}
Antonio Ambrosetti and Paul~H. Rabinowitz.
\newblock Dual variational methods in critical point theory and applications.
\newblock {\em J. Functional Analysis}, 14:349--381, 1973.

\bibitem{AmbrosioBraides}
Luigi Ambrosio and Andrea Braides.
\newblock Functionals defined on partitions in sets of finite perimeter. {II}.
  {S}emicontinuity, relaxation and homogenization.
\newblock {\em J. Math. Pures Appl. (9)}, 69(3):307--333, 1990.

\bibitem{AmmannHumbert}
Bernd Ammann and Emmanuel Humbert.
\newblock The second {Y}amabe invariant.
\newblock {\em J. Funct. Anal.}, 235(2):377--412, 2006.

\bibitem{vanGoethem}
Samuel Amstutz, Antonio~Andr\'{e} Novotny, and Nicolas Van~Goethem.
\newblock Minimal partitions and image classification using a gradient-free
  perimeter approximation.
\newblock {\em Inverse Probl. Imaging}, 8(2):361--387, 2014.

\bibitem{AndradeSantosSantosTavares}
P\^edra Andrade, Ederson Moreira~dos Santos, Makson Santos, and Hugo Tavares.
\newblock Optimal partition problems with volume constraint.
\newblock arXiv:2305.02870.

\bibitem{AtkinsonBrezisPeletier}
Frederik Atkinson, Haim Brezis, and Lambertus Peletier.
\newblock Nodal solutions of elliptic equations with critical {S}obolev
  exponents.
\newblock {\em J. Differential Equations}, 85(1):151--170, 1990.

\bibitem{Aubin1}
Thierry Aubin.
\newblock Probl\`emes isop\'{e}rim\'{e}triques et espaces de {S}obolev.
\newblock {\em J. Differential Geometry}, 11(4):573--598, 1976.

\bibitem{AvilaYang}
Andr{{\'e}}s~I. {{\'A}}vila and Jianfu Yang.
\newblock On the existence and shape of least energy solutions for some
  elliptic systems.
\newblock {\em J. Differential Equations}, 191(2):348--376, 2003.

\bibitem{BadialeSerra}
Marino Badiale and Enrico Serra.
\newblock {\em Semilinear elliptic equations for beginners:
Existence results via the variational approach.}.
\newblock Universitext. Springer, London, 2011.
\newblock 

\bibitem{BahriCoron}
Abbas Bahri and Jean-Michel Coron.
\newblock On a nonlinear elliptic equation involving the critical {S}obolev
  exponent: the effect of the topology of the domain.
\newblock {\em Comm. Pure Appl. Math.}, 41(3):253--294, 1988.

\bibitem{BBHU}
Ram Band, Gregory Berkolaiko, Hillel Raz, and Uzy Smilansky.
\newblock The number of nodal domains on quantum graphs as a stability index of
  graph partitions.
\newblock {\em Comm. Math. Phys.}, 311(3):815--838, 2012.

\bibitem{BartschJeanjean}
Thomas Bartsch and Louis Jeanjean.
\newblock Normalized solutions for nonlinear {S}chr\"odinger systems.
\newblock {\em Proc. Roy. Soc. Edinburgh Sect. A}, 148(2):225--242, 2018.

\bibitem{BartschJeanjeanSoave}
Thomas Bartsch, Louis Jeanjean, and Nicola Soave.
\newblock Normalized solutions for a system of coupled cubic {S}chr\"odinger
  equations on {$\mathbb{R}^3$}.
\newblock {\em J. Math. Pures Appl. (9)}, 106(4):583--614, 2016.

\bibitem{L2_5}
Thomas Bartsch, Riccardo Molle, Matteo Rizzi, and Gianmaria Verzini.
\newblock Normalized solutions of mass supercritical {S}chr\"{o}dinger
  equations with potential.
\newblock {\em Comm. Partial Differential Equations}, 46(9):1729--1756, 2021.

\bibitem{BartschSoave2}
Thomas Bartsch and Nicola Soave.
\newblock A natural constraint approach to normalized solutions of nonlinear
  {S}chr\"odinger equations and systems.
\newblock {\em J. Funct. Anal.}, 272(12):4998--5037, 2017.

\bibitem{BartschSoave}
Thomas Bartsch and Nicola Soave.
\newblock Multiple normalized solutions for a competing system of
  {S}chr\"{o}dinger equations.
\newblock {\em Calc. Var. Partial Differential Equations}, 58(1):Paper No. 22,
2019.

\bibitem{BartschWang}
Thomas Bartsch and Zhi-Qiang Wang.
\newblock Note on ground states of nonlinear {S}chr\"odinger systems.
\newblock {\em J. Partial Differential Equations}, 19(3):200--207, 2006.

\bibitem{BartschWethWillem}
Thomas Bartsch, Tobias Weth, and Michel Willem.
\newblock Partial symmetry of least energy nodal solutions to some variational
  problems.
\newblock {\em J. Anal. Math.}, 96:1--18, 2005.

\bibitem{BartschWillem95}
Thomas Bartsch and Michel Willem.
\newblock On an elliptic equation with concave and convex nonlinearities.
\newblock {\em Proc. Amer. Math. Soc.}, 123(11):3555--3561, 1995.

\bibitem{L2_6}
Thomas Bartsch, Xuexiu Zhong, and Wenming Zou.
\newblock Normalized solutions for a coupled {S}chr\"{o}dinger system.
\newblock {\em Math. Ann.}, 380(3-4):1713--1740, 2021.

\bibitem{Bellazzinietal}
Jacopo Bellazzini, Nabile Boussa\"id, Louis Jeanjean, and Nicola Visciglia.
\newblock Existence and stability of standing waves for supercritical {NLS}
  with a partial confinement.
\newblock {\em Comm. Math. Phys.}, 353(1):229--251, 2017.

\bibitem{BTWW}
Henri Berestycki, Susanna Terracini, Kelei Wang, and Juncheng Wei.
\newblock On entire solutions of an elliptic system modeling phase separations.
\newblock {\em Adv. Math.}, 243:102--126, 2013.

\bibitem{Berkolaiko}
Gregory Berkolaiko, Peter Kuchment, and Uzy Smilansky.
\newblock Critical partitions and nodal deficiency of billiard eigenfunctions.
\newblock {\em Geom. Funct. Anal.}, 22(6):1517--1540, 2012.

\bibitem{BogoselVelichkov}
Beniamin Bogosel and Bozhidar Velichkov.
\newblock A multiphase shape optimization problem for eigenvalues: qualitative
  study and numerical results.
\newblock {\em SIAM J. Numer. Anal.}, 54(1):210--241, 2016.

\bibitem{BonheureCheikhNascimento}
Denis Bonheure, Hussein Cheikh~Ali, and Robson Nascimento.
\newblock A {P}aneitz-{B}ranson type equation with {N}eumann boundary
  conditions.
\newblock {\em Adv. Calc. Var.}, 14(4):499--519, 2021.

\bibitem{BonheureSantosTavares}
Denis Bonheure, Ederson~Moreira dos Santos, and Hugo Tavares.
\newblock Hamiltonian elliptic systems: a guide to variational frameworks.
\newblock {\em Port. Math.}, 71(3-4):301--395, 2014.

\bibitem{BFMST18}
Denis Bonheure, Juraj F\"{o}ldes, Ederson Moreira~dos Santos, Alberto
  Salda\~{n}a, and Hugo Tavares.
\newblock Paths to uniqueness of critical points and applications to partial
  differential equations.
\newblock {\em Trans. Amer. Math. Soc.}, 370(10):7081--7127, 2018.

\bibitem{BSPTW}
Denis Bonheure, Ederson Moreira~dos Santos, Enea Parini, Hugo Tavares, and
  Tobias Weth.
\newblock Nodal solutions for sublinear-type problems with {D}irichlet boundary
  conditions.
\newblock {\em Int. Math. Res. Not. IMRN}, (5):3760--3804, 2022.

\bibitem{BMRT15}
Denis Bonheure, Ederson Moreira~dos Santos, Miguel Ramos, and Hugo Tavares.
\newblock Existence and symmetry of least energy nodal solutions for
  {H}amiltonian elliptic systems.
\newblock {\em J. Math. Pures Appl. (9)}, 104(6):1075--1107, 2015.

\bibitem{BonheureSerraTilli}
Denis Bonheure, Enrico Serra, and Paolo Tilli.
\newblock Radial positive solutions of elliptic systems with {N}eumann boundary
  conditions.
\newblock {\em J. Funct. Anal.}, 265(3):375--398, 2013.

\bibitem{BHV}
Virginie Bonnaillie-No{\"e}l, Bernard Helffer, and Gregory Vial.
\newblock Numerical simulations for nodal domains and spectral minimal
  partitions.
\newblock {\em ESAIM Control Optim. Calc. Var.}, 16(1):221--246, 2010.

\bibitem{BourdinBucurOudet}
Blaise Bourdin, Dorin Bucur, and {{\'E}}douard Oudet.
\newblock Optimal partitions for eigenvalues.
\newblock {\em SIAM J. Sci. Comput.}, 31(6):4100--4114, 2009/10.

\bibitem{BrezisNirenberg}
Ha{\"{\i}}m Br{{\'e}}zis and Louis Nirenberg.
\newblock Positive solutions of nonlinear elliptic equations involving critical
  {S}obolev exponents.
\newblock {\em Comm. Pure Appl. Math.}, 36(4):437--477, 1983.

\bibitem{Brianconetal}
Tanguy Brian\c{c}on, Mohammed Hayouni, and Michel Pierre.
\newblock Lipschitz continuity of state functions in some optimal shaping.
\newblock {\em Calc. Var. Partial Differential Equations}, 23(1):13--32, 2005.

\bibitem{BucurButtazzoBook}
Dorin Bucur and Giuseppe Buttazzo.
\newblock {\em Variational methods in shape optimization problems}.
\newblock Progress in Nonlinear Differential Equations and their Applications,
  65. Birkh{\"a}user Boston Inc., Boston, MA, 2005.

\bibitem{BucurButtazzoHenrot}
Dorin Bucur, Giuseppe Buttazzo, and Antoine Henrot.
\newblock Existence results for some optimal partition problems.
\newblock {\em Adv. Math. Sci. Appl.}, 8(2):571--579, 1998.

\bibitem{BucurVelichkov}
Dorin Bucur and Bozhidar Velichkov.
\newblock Multiphase shape optimization problems.
\newblock {\em SIAM J. Control Optim.}, 52(6):3556--3591, 2014.

\bibitem{ButtazzoDalMaso}
Giuseppe Buttazzo and Gianni Dal~Maso.
\newblock Shape optimization for {D}irichlet problems: relaxed formulation and
  optimality conditions.
\newblock {\em Appl. Math. Optim.}, 23(1):17--49, 1991.

\bibitem{ButtazzoTimofte}
Giuseppe Buttazzo and Claudia Timofte.
\newblock On the relaxation of some optimal partition problems.
\newblock {\em Adv. Math. Sci. Appl.}, 12(2):509--520, 2002.

\bibitem{CaffarelliLin2}
Luis Caffarelli and Fang-Hua Lin.
\newblock An optimal partition problem for eigenvalues.
\newblock {\em J. Sci. Comput.}, 31(1-2):5--18, 2007.

\bibitem{CaffarelliLin}
Luis Caffarelli and Fang-Hua Lin.
\newblock Singularly perturbed elliptic systems and multi-valued harmonic
  functions with free boundaries.
\newblock {\em J. Amer. Math. Soc.}, 21(3):847--862, 2008.

\bibitem{CaffarelliPatriziQuitalo}
Luis Caffarelli, Stefania Patrizi, and Ver\'onica Quitalo.
\newblock On a long range segregation model.
\newblock {\em J. Eur. Math. Soc. (JEMS)}, 19(12):3575--3628, 2017.

\bibitem{CassaniTavaresZhang}
Daniele Cassani, Hugo Tavares, and Jianjun Zhang.
\newblock Bose fluids and positive solutions to weakly coupled systems with
  critical growth in dimension two.
\newblock {\em J. Differential Equations}, 269(3):2328--2385, 2020.

\bibitem{CastroClapp}
Alfonso Castro and M\'{o}nica Clapp.
\newblock The effect of the domain topology on the number of minimal nodal
  solutions of an elliptic equation at critical growth in a symmetric domain.
\newblock {\em Nonlinearity}, 16(2):579--590, 2003.

\bibitem{CastroCossioNeuberger}
Alfonso Castro, Jorge Cossio, and John~M. Neuberger.
\newblock A sign-changing solution for a superlinear {D}irichlet problem.
\newblock {\em Rocky Mountain J. Math.}, 27(4):1041--1053, 1997.

\bibitem{CazenaveLions1982}
Thierry~Cazenave and Pierre-Louis Lions.
\newblock Orbital stability of standing waves for some nonlinear
  {S}chr\"odinger equations.
\newblock {\em Comm. Math. Phys.}, 85(4):549--561, 1982.

\bibitem{Cazenave2003}
Thierry Cazenave.
\newblock {\em Semilinear {S}chr\"odinger equations}, volume~10 of {\em Courant
  Lecture Notes in Mathematics}.
\newblock New York University Courant Institute of Mathematical Sciences, New
  York, 2003.

\bibitem{CeramiSoliminiStruwe}
Giovanna Cerami, Sergio Solimini, and Micha\"el Struwe.
\newblock Some existence results for superlinear elliptic boundary value
  problems involving critical exponents.
\newblock {\em J. Funct. Anal.}, 69(3):289--306, 1986.

\bibitem{C57}
Subrahmanyan Chandrasekhar.
\newblock {\em An introduction to the study of stellar structure}.
\newblock Dover Publications, Inc., New York, N. Y., 1957.

\bibitem{Chang}
Jinyong Chang.
\newblock Note on ground states of a nonlinear {S}chr{\"o}dinger system.
\newblock {\em J. Math. Anal. Appl.}, 381(2):957--962, 2011.

\bibitem{clll}
Shu-Ming Chang, Chang-Shou Lin, Tai-Chia Lin, and Wen-Wei Lin.
\newblock Segregated nodal domains of two-dimensional multispecies
  {B}ose-{E}instein condensates.
\newblock {\em Phys. D}, 196(3-4):341--361, 2004.

\bibitem{ChenLin}
Zhijie Chen and Chang-Shou Lin.
\newblock Asymptotic behavior of least energy solutions for a critical elliptic
  system.
\newblock {\em Int. Math. Res. Not. IMRN}, (21):11045--11082, 2015.

\bibitem{ChenLinZou}
Zhijie Chen, Chang-Shou Lin, and Wenming Zou.
\newblock Sign-changing solutions and phase separation for an elliptic system
  with critical exponent.
\newblock {\em Comm. Partial Differential Equations}, 39(10):1827--1859, 2014.

\bibitem{ChenZouARMA2012}
Zhijie Chen and Wenming Zou.
\newblock Positive least energy solutions and phase separation for coupled
  {S}chr{\"o}dinger equations with critical exponent.
\newblock {\em Arch. Ration. Mech. Anal.}, 205(2):515--551, 2012.

\bibitem{ChenZou}
Zhijie Chen and Wenming Zou.
\newblock An optimal constant for the existence of least energy solutions of a
  coupled {S}chr\"odinger system.
\newblock {\em Calc. Var. Partial Differential Equations}, 48(3-4):695--711,
  2013.

\bibitem{ChenZou2}
Zhijie Chen and Wenming Zou.
\newblock Positive least energy solutions and phase separation for coupled
  {S}chr\"{o}dinger equations with critical exponent: higher dimensional case.
\newblock {\em Calc. Var. Partial Differential Equations}, 52(1-2):423--467,
  2015.

\bibitem{CirantVerzini}
Marco Cirant and Gianmaria Verzini.
\newblock Bifurcation and segregation in quadratic two-populations mean field
  games systems.
\newblock {\em ESAIM Control Optim. Calc. Var.}, 23(3):1145--1177, 2017.

\bibitem{ClappFernandezSaldana}
M\'{o}nica Clapp, Juan~Carlos Fern\'{a}ndez, and Alberto Salda\~{n}a.
\newblock Critical polyharmonic systems and optimal partitions.
\newblock {\em Commun. Pure Appl. Anal.}, 20(11):4007--4023, 2021.

\bibitem{ClappMussoPistoia}
M\'{o}nica Clapp, Monica Musso, and Angela Pistoia.
\newblock Multipeak solutions to the {B}ahri-{C}oron problem in domains with a
  shrinking hole.
\newblock {\em J. Funct. Anal.}, 256(2):275--306, 2009.

\bibitem{ClappPistoia2018}
M\'{o}nica Clapp and Angela Pistoia.
\newblock Existence and phase separation of entire solutions to a pure critical
  competitive elliptic system.
\newblock {\em Calc. Var. Partial Differential Equations}, 57(1):Paper No. 23,
2018.

\bibitem{ClappPistoia_symmetries}
M\'{o}nica Clapp and Angela Pistoia.
\newblock Yamabe systems and optimal partitions on manifolds with symmetries.
\newblock {\em Electron. Res. Arch.}, 29(6):4327--4338, 2021.

\bibitem{ClappPistoia2022}
M\'{o}nica Clapp and Angela Pistoia.
\newblock Fully nontrivial solutions to elliptic systems with mixed couplings.
\newblock {\em Nonlinear Anal.}, 216:Paper No. 112694, 2022.

\bibitem{ClappPistoiaTavares}
M{\'o}nica Clapp, Angela Pistoia, and Hugo Tavares.
\newblock Yamabe systems, optimal partitions, and nodal solutions to the
  {Y}amabe equation.
\newblock {\em J. Eur. Math. Soc. (JEMS)}, to appear.

\bibitem{ClappSaldanaSzulkin_sphere}
M\'{o}nica Clapp, Alberto Salda\~{n}a, and Andrzej Szulkin.
\newblock Phase separation, optimal partitions, and nodal solutions to the
  {Y}amabe equation on the sphere.
\newblock {\em Int. Math. Res. Not. IMRN}, (5):3633--3652, 2021.

\bibitem{ClappSzulkin2019}
M\'{o}nica Clapp and Andrzej Szulkin.
\newblock A simple variational approach to weakly coupled competitive elliptic
  systems.
\newblock {\em NoDEA Nonlinear Differential Equations Appl.}, 26(4):Paper No.
  26, 2019.

\bibitem{ClappWeth2004}
M\'{o}nica Clapp and Tobias Weth.
\newblock Minimal nodal solutions of the pure critical exponent problem on a
  symmetric domain.
\newblock {\em Calc. Var. Partial Differential Equations}, 21(1):1--14, 2004.

\bibitem{ClappWeth}
M\'{o}nica Clapp and Tobias Weth.
\newblock Multiple solutions for the {B}rezis-{N}irenberg problem.
\newblock {\em Adv. Differential Equations}, 10(4):463--480, 2005.

\bibitem{ClementdeFigueiredoMitidieri}
Philippe Cl{{\'e}}ment, Djairo de~Figueiredo, and Enzo Mitidieri.
\newblock Positive solutions of semilinear elliptic systems.
\newblock {\em Comm. Partial Differential Equations}, 17(5-6):923--940, 1992.

\bibitem{ClementvanderVorst}
Philippe Cl{{\'e}}ment and Robert C. A.~M. Van~der Vorst.
\newblock On a semilinear elliptic system.
\newblock {\em Differential Integral Equations}, 8(6):1317--1329, 1995.

\bibitem{CK91}
Myriam Comte and Mariette~C. Knaap.
\newblock Existence of solutions of elliptic equations involving critical
  {S}obolev exponents with {N}eumann boundary condition in general domains.
\newblock {\em Differential Integral Equations}, 4(6):1133--1146, 1991.

\bibitem{ContiTerraciniVerziniNehari}
Monica Conti, Susanna Terracini, and Gianmaria Verzini.
\newblock Nehari's problem and competing species systems.
\newblock {\em Ann. Inst. H. Poincar{\'e} Anal. Non Lin{\'e}aire},
  19(6):871--888, 2002.

\bibitem{ContiTerraciniVerziniOptimalNonlinear}
Monica Conti, Susanna Terracini, and Gianmaria Verzini.
\newblock An optimal partition problem related to nonlinear eigenvalues.
\newblock {\em J. Funct. Anal.}, 198(1):160--196, 2003.

\bibitem{ContiTerraciniVerziniOPP}
Monica Conti, Susanna Terracini, and Gianmaria Verzini.
\newblock On a class of optimal partition problems related to the {F}u\v c\'\i
  k spectrum and to the monotonicity formulae.
\newblock {\em Calc. Var. Partial Differential Equations}, 22(1):45--72, 2005.

\bibitem{CorreiaJDE2016}
Sim\~{a}o Correia.
\newblock Characterization of ground-states for a system of {$M$} coupled
  semilinear {S}chr\"{o}dinger equations and applications.
\newblock {\em J. Differential Equations}, 260(4):3302--3326, 2016.

\bibitem{CorreiaNA2016}
Sim\~{a}o Correia.
\newblock Ground-states for systems of {$M$} coupled semilinear
  {S}chr\"{o}dinger equations with attraction-repulsion effects:
  characterization and perturbation results.
\newblock {\em Nonlinear Anal.}, 140:112--129, 2016.

\bibitem{CorreiaOliveiraTavares}
Sim\~{a}o Correia, Filipe Oliveira, and Hugo Tavares.
\newblock Semitrivial vs. fully nontrivial ground states in cooperative cubic
  {S}chr\"{o}dinger systems with {$d\geq 3$} equations.
\newblock {\em J. Funct. Anal.}, 271(8):2247--2273, 2016.

\bibitem{Dancer1988}
Edward.~N. Dancer.
\newblock The effect of domain shape on the number of positive solutions of
  certain nonlinear equations.
\newblock {\em J. Differential Equations}, 74(1):120--156, 1988.

\bibitem{deFigueiredo}
Djairo~G. de~Figueiredo.
\newblock Semilinear elliptic systems: existence, multiplicity, symmetry of
  solutions.
\newblock In {\em Handbook of differential equations: stationary partial
  differential equations. {V}ol. {V}}, Handb. Differ. Equ., pages 1--48.
  Elsevier/North-Holland, Amsterdam, 2008.

\bibitem{PhilippisSpolaorVelichkov}
Guido De~Philippis, Luca Spolaor, and Bozhidar Velichkov.
\newblock Regularity of the free boundary for the two-phase {B}ernoulli
  problem.
\newblock {\em Invent. Math.}, 225(2):347--394, 2021.

\bibitem{DiasTavares}
Manuel Dias and Hugo Tavares.
\newblock Optimal uniform bounds for competing variational elliptic systems
  with variable coefficients.
\newblock {\em Nonlinear Anal.}, 235:Paper No. 113348, 2023.

\bibitem{DovettaPistoia}
Simone Dovetta and Angela Pistoia.
\newblock Solutions to a cubic {S}chr\"{o}dinger system with mixed attractive
  and repulsive forces in a critical regime.
\newblock {\em Math. Eng.}, 4(4):Paper No. 027, 2022.

\bibitem{EspositoPistoiaVetois}
Pierpaolo Esposito, Angela Pistoia, and J\'{e}r\^{o}me V\'{e}tois.
\newblock The effect of linear perturbations on the {Y}amabe problem.
\newblock {\em Math. Ann.}, 358(1-2):511--560, 2014.

\bibitem{EvansPDE}
Lawrence~C. Evans.
\newblock {\em Partial differential equations}, volume~19 of {\em Graduate
  Studies in Mathematics}.
\newblock American Mathematical Society, Providence, RI, second edition, 2010.

\bibitem{FibichMerle2001}
Gadi Fibich and Frank Merle.
\newblock Self-focusing on bounded domains.
\newblock {\em Phys. D}, 155(1-2):132--158, 2001.

\bibitem{Fukuizumi2012}
Reika Fukuizumi, Fouad~Hadj Selem, and Hiroaki Kikuchi.
\newblock Stationary problem related to the nonlinear Schr{\"o}dinger equation
  on the unit ball.
\newblock {\em Nonlinearity}, 25(8):2271, 2012.

\bibitem{Mariana2}
Nicola Garofalo, Arshak Petrosyan, and Mariana Smit Vega~Garcia.
\newblock An epiperimetric inequality approach to the regularity of the free
  boundary in the {S}ignorini problem with variable coefficients.
\newblock {\em J. Math. Pures Appl. (9)}, 105(6):745--787, 2016.

\bibitem{Mariana1}
Nicola Garofalo and Mariana Smit Vega~Garcia.
\newblock New monotonicity formulas and the optimal regularity in the
  {S}ignorini problem with variable coefficients.
\newblock {\em Adv. Math.}, 262:682--750, 2014.

\bibitem{GidasNiNirenberg}
Basilis Gidas, Wei~Ming Ni, and Louis Nirenberg.
\newblock Symmetry and related properties via the maximum principle.
\newblock {\em Comm. Math. Phys.}, 68(3):209--243, 1979.

\bibitem{GuoJeanjean}
Tianxiang Gou and Louis Jeanjean.
\newblock Existence and orbital stability of standing waves for nonlinear
  {S}chr\"odinger systems.
\newblock {\em Nonlinear Anal.}, 144:10--22, 2016.

\bibitem{GuoJeanjean2}
Tianxiang Gou and Louis Jeanjean.
\newblock Multiple positive normalized solutions for nonlinear schr{\"o}dinger
  systems.
\newblock {\em Nonlinearity}, 31(5):2319, 2018.

\bibitem{GrillakisShatahStrauss}
Manoussos Grillakis, Jalal Shatah, and Walter Strauss.
\newblock Stability theory of solitary waves in the presence of symmetry, I.
\newblock {\em Journal of Functional Analysis}, 74(1):160--197, 1987.

\bibitem{Grossi:2021aa}
Massimo Grossi, Isabella Ianni, Peng Luo, and Shusen Yan.
\newblock Non-degeneracy and local uniqueness of positive solutions to the
  {L}ane-{E}mden problem in dimension two.
\newblock {\em J. Math. Pures Appl. (9)}, 157:145--210, 2022.

\bibitem{GrossiSaldanaTavares}
Massimo Grossi, Alberto Salda\~{n}a, and Hugo Tavares.
\newblock Sharp concentration estimates near criticality for radial
  sign-changing solutions of {D}irichlet and {N}eumann problems.
\newblock {\em Proc. Lond. Math. Soc. (3)}, 120(1):39--64, 2020.

\bibitem{Guerra2}
I.~A. Guerra.
\newblock Asymptotic behaviour of a semilinear elliptic system with a large
  exponent.
\newblock {\em J. Dynam. Differential Equations}, 19(1):243--263, 2007.

\bibitem{Guerra}
Ignacio Guerra.
\newblock Solutions of an elliptic system with a nearly critical exponent.
\newblock {\em Ann. Inst. H. Poincar{\'e} Anal. Non Lin{\'e}aire},
  25(1):181--200, 2008.

\bibitem{GuoPeng}
Qing Guo and Shuangjie Peng.
\newblock Sign-changing solutions to the slightly supercritical {L}ane-{E}mden
  system with {N}eumann boundary conditions.
\newblock arXiv: 2306.00663.

\bibitem{GuoLuoZou}
Yuxia Guo, Senping Luo, and Wenming Zou.
\newblock The existence, uniqueness and nonexistence of the ground state to the
  {$N$}-coupled {S}chr\"{o}dinger systems in {$\Bbb R^n$} {$(n\leqslant4)$}.
\newblock {\em Nonlinearity}, 31(1):314--339, 2018.

\bibitem{HanLin}
Qing Han and Fanghua Lin.
\newblock {\em Elliptic partial differential equations}, volume~1 of {\em
  Courant Lecture Notes in Mathematics}.
\newblock Courant Institute of Mathematical Sciences, New York; American
  Mathematical Society, Providence, RI, second edition, 2011.

\bibitem{HanAIHP}
Zheng-Chao Han.
\newblock Asymptotic approach to singular solutions for nonlinear elliptic
  equations involving critical {S}obolev exponent.
\newblock {\em Ann. Inst. H. Poincar{\'e} Anal. Non Lin{\'e}aire},
  8(2):159--174, 1991.

\bibitem{Hardtetal}
Robert Hardt, Maria Hoffmann-Ostenhof, Thomas Hoffmann-Ostenhof, and Nikolai
  Nadirashvili.
\newblock Critical sets of solutions to elliptic equations.
\newblock {\em J. Differential Geom.}, 51(2):359--373, 1999.

\bibitem{HeYang}
Qihan He and Jing Yang.
\newblock Quantitative properties of ground-states to an {$M$}-coupled system
  with critical exponent in {$\Bbb R^N$}.
\newblock {\em Sci. China Math.}, 61(4):709--726, 2018.

\bibitem{HeHeZhang}
Tieshan He, Lang He, and Meng Zhang.
\newblock The {B}r\'{e}zis-{N}irenberg type problem for the {$p$}-{L}aplacian:
  infinitely many sign-changing solutions.
\newblock {\em Calc. Var. Partial Differential Equations}, 59(3):Paper No. 98,
  14, 2020.

\bibitem{HeZhangWuLiang}
Tieshan He, Chaolong Zhang, Dongqing Wu, and Kaihao Liang.
\newblock An existence result for sign-changing solutions of the
  {B}r\'{e}zis-{N}irenberg problem.
\newblock {\em Appl. Math. Lett.}, 84:90--95, 2018.

\bibitem{helffer}
Bernard Helffer.
\newblock On spectral minimal partitions: a survey.
\newblock {\em Milan J. Math.}, 78(2):575--590, 2010.

\bibitem{HHT}
Bernard Helffer, Thomas Hoffmann-Ostenhof, and Susanna Terracini.
\newblock Nodal domains and spectral minimal partitions.
\newblock {\em Ann. Inst. H. Poincar{\'e} Anal. Non Lin{\'e}aire},
  26(1):101--138, 2009.

\bibitem{HHOT2}
Bernard Helffer, Thomas Hoffmann-Ostenhof, and Susanna Terracini.
\newblock Nodal minimal partitions in dimension 3.
\newblock {\em Discrete Contin. Dyn. Syst.}, 28(2):617--635, 2010.

\bibitem{HHOT3}
Bernard Helffer, Thomas Hoffmann-Ostenhof, and Susanna Terracini.
\newblock On spectral minimal partitions: the case of the sphere.
\newblock In {\em Around the research of {V}ladimir {M}az'ya. {III}}, volume~13
  of {\em Int. Math. Ser. (N. Y.)}, pages 153--178. Springer, New York, 2010.

\bibitem{ShapeOptimizationBook}
Antoine Henrot, editor.
\newblock {\em Shape optimization and spectral theory}.
\newblock De Gruyter Open, Warsaw, 2017.

\bibitem{Jeanjean}
Louis Jeanjean.
\newblock Existence of solutions with prescribed norm for semilinear elliptic
  equations.
\newblock {\em Nonlinear Anal.}, 28(10):1633--1659, 1997.

\bibitem{L2_2}
Louis Jeanjean, Jacek Jendrej, Thanh~Trung Le, and Nicola Visciglia.
\newblock Orbital stability of ground states for a {S}obolev critical
  {S}chr\"{o}dinger equation.
\newblock {\em J. Math. Pures Appl. (9)}, 164:158--179, 2022.

\bibitem{L2_1}
Louis Jeanjean and Thanh~Trung Le.
\newblock Multiple normalized solutions for a {S}obolev critical
  {S}chr\"{o}dinger equation.
\newblock {\em Math. Ann.}, 384(1-2):101--134, 2022.

\bibitem{Kawohl_uniq}
Bernhard Kawohl.
\newblock {\em Rearrangements and convexity of level sets in {PDE}}, volume
  1150 of {\em Lecture Notes in Mathematics}.
\newblock Springer-Verlag, Berlin, 1985.

\bibitem{Krasnoselskii}
Mark Krasnosel'skii.
\newblock {\em Positive solutions of operator equations}.
\newblock P. Noordhoff Ltd. Groningen, 1964.
\newblock Translated from the Russian by Richard E. Flaherty; edited by Leo F.
  Boron.

\bibitem{KriventsovLin}
Dennis Kriventsov and Fang-Hua Lin.
\newblock Regularity for shape optimizers: the nondegenerate case.
\newblock {\em Comm. Pure Appl. Math.}, 71(8):1535--1596, 2018.

\bibitem{Kukavica}
Igor Kukavica.
\newblock Quantitative uniqueness for second-order elliptic operators.
\newblock {\em Duke Math. J.}, 91(2):225--240, 1998.

\bibitem{Kwong}
Man~Kam Kwong.
\newblock Uniqueness of positive solutions of {$\Delta u-u+u^p=0$} in {${\bf
  R}^n$}.
\newblock {\em Arch. Rational Mech. Anal.}, 105(3):243--266, 1989.

\bibitem{LamboleyPierre}
Jimmy Lamboley and Michel Pierre.
\newblock Regularity of optimal spectral domains.
\newblock In {\em Shape optimization and spectral theory}, pages 29--77. De
  Gruyter Open, Warsaw, 2017.

\bibitem{LeeParker}
John~M. Lee and Thomas~H. Parker.
\newblock The {Y}amabe problem.
\newblock {\em Bull. Amer. Math. Soc. (N.S.)}, 17(1):37--91, 1987.

\bibitem{LinNiTakagi}
Chang-Shou Lin, Wei~Ming Ni, and Izumi Takagi.
\newblock Large amplitude stationary solutions to a chemotaxis system.
\newblock {\em J. Differential Equations}, 72(1):1--27, 1988.

\bibitem{Lin}
Fang-Hua Lin.
\newblock Nodal sets of solutions of elliptic and parabolic equations.
\newblock {\em Comm. Pure Appl. Math.}, 44(3):287--308, 1991.

\bibitem{LinWeiErratum}
Tai-Chia Lin and Juncheng Wei.
\newblock Erratum: ``{G}round state of {$N$} coupled nonlinear
  {S}chr{\"o}dinger equations in {${\bf R}^n$}, {$n\leq3$}'' [{C}omm. {M}ath.
  {P}hys. {\bf 255} (2005), no. 3, 629--653; mr2135447].
\newblock 277(2):573--576.

\bibitem{LinWei}
Tai-Chia Lin and Juncheng Wei.
\newblock Ground state of {$N$} coupled nonlinear {S}chr{\"o}dinger equations
  in {$\bold R^n$}, {$n\leq 3$}.
\newblock {\em Comm. Math. Phys.}, 255(3):629--653, 2005.

\bibitem{LinWei2}
Tai-Chia Lin and Juncheng Wei.
\newblock Spikes in two coupled nonlinear {S}chr\"odinger equations.
\newblock {\em Ann. Inst. H. Poincar\'e Anal. Non Lin\'eaire}, 22(4):403--439,
  2005.

\bibitem{LiuWang}
Zhaoli Liu and Zhi-Qiang Wang.
\newblock Ground states and bound states of a nonlinear {S}chr{\"o}dinger
  system.
\newblock {\em Adv. Nonlinear Stud.}, 10(1):175--193, 2010.

\bibitem{MaiaMontefuscoPellacci}
Liliane~A. Maia, Eugenio~Montefusco and Benedetta~Pellacci.
\newblock Positive solutions for a weakly coupled nonlinear {S}chr\"odinger
  system.
\newblock {\em J. Differential Equations}, 229(2):743--767, 2006.

\bibitem{Mandel}
Rainer Mandel.
\newblock Minimal energy solutions for cooperative nonlinear {S}chr{\"o}dinger
  systems.
\newblock {\em NoDEA Nonlinear Differential Equations Appl.}, 22(2):239--262,
  2015.

\bibitem{MazzoleniTerraciniVelichkov}
Dario Mazzoleni, Susanna Terracini, and Bozhidar Velichkov.
\newblock Regularity of the optimal sets for some spectral functionals.
\newblock {\em Geom. Funct. Anal.}, 27(2):373--426, 2017.

\bibitem{MazzoleniTreyVelichkov}
Dario Mazzoleni, Baptiste Trey, and Bozhidar Velichkov.
\newblock Regularity of the optimal sets for the second {D}irichlet eigenvalue.
\newblock {\em Ann. Inst. H. Poincar\'{e} C Anal. Non Lin\'{e}aire},
  39(3):529--573, 2022.

\bibitem{MollePassaseo}
Riccardo Molle and Donato Passaseo.
\newblock Infinitely many solutions for elliptic equations with non-symmetric
  nonlinearities.
\newblock {\em Calc. Var. Partial Differential Equations}, 61(3):Paper No. 115,
 2022.

\bibitem{L2_3}
Riccardo Molle, Giuseppe Riey, and Gianmaria Verzini.
\newblock Normalized solutions to mass supercritical {S}chr\"{o}dinger
  equations with negative potential.
\newblock {\em J. Differential Equations}, 333:302--331, 2022.

\bibitem{LEfullynonlinear}
Ederson Moreira~dos Santos, Gabrielle Nornberg, Delia Schiera, and Hugo
  Tavares.
\newblock Principal spectral curves for {L}ane-{E}mden fully nonlinear type
  systems and applications.
\newblock {\em Calc. Var. Partial Differential Equations}, 62(2):Paper No. 49,
  2023.

\bibitem{MussoPistoiaIndiana2002}
Monica Musso and Angela Pistoia.
\newblock Multispike solutions for a nonlinear elliptic problem involving the
  critical {S}obolev exponent.
\newblock {\em Indiana Univ. Math. J.}, 51(3):541--579, 2002.

\bibitem{MussoPistoiaJMPA2006}
Monica Musso and Angela Pistoia.
\newblock Sign changing solutions to a nonlinear elliptic problem involving the
  critical {S}obolev exponent in pierced domains.
\newblock {\em J. Math. Pures Appl. (9)}, 86(6):510--528, 2006.

\bibitem{MussoPistoia2008}
Monica Musso and Angela Pistoia.
\newblock Sign changing solutions to a {B}ahri-{C}oron's problem in pierced
  domains.
\newblock {\em Discrete Contin. Dyn. Syst.}, 21(1):295--306, 2008.

\bibitem{NTTV1}
Benedetta Noris, Hugo Tavares, Susanna Terracini, and Gianmaria Verzini.
\newblock Uniform {H}{\"o}lder bounds for nonlinear {S}chr{\"o}dinger systems
  with strong competition.
\newblock {\em Comm. Pure Appl. Math.}, 63(3):267--302, 2010.

\bibitem{MR2928850}
Benedetta Noris, Hugo Tavares, Susanna Terracini, and Gianmaria Verzini.
\newblock Convergence of minimax structures and continuation of critical points
  for singularly perturbed systems.
\newblock {\em J. Eur. Math. Soc. (JEMS)}, 14(4):1245--1273, 2012.

\bibitem{ntvAnPDE}
Benedetta Noris, Hugo Tavares, and Gianmaria Verzini.
\newblock Existence and orbital stability of the ground states with prescribed
  mass for the {$L^2$}-critical and supercritical {NLS} on bounded domains.
\newblock {\em Anal. PDE}, 7(8):1807--1838, 2014.

\bibitem{NTV2}
Benedetta Noris, Hugo Tavares, and Gianmaria Verzini.
\newblock Stable solitary waves with prescribed ${L}^2$-mass for the cubic
  {S}chr\"odinger system with trapping potentials.
\newblock {\em Discrete and Continuous Dynamical Systems - Series A},
  35(12):6085--6112, 2015.

\bibitem{NTV3}
Benedetta Noris, Hugo Tavares, and Gianmaria Verzini.
\newblock Normalized solutions for nonlinear {S}chr\"{o}dinger systems on
  bounded domains.
\newblock {\em Nonlinearity}, 32(3):1044--1072, 2019.

\bibitem{OliveiraTavares}
Filipe Oliveira and Hugo Tavares.
\newblock Ground states for a nonlinear {S}chr\"{o}dinger system with sublinear
  coupling terms.
\newblock {\em Adv. Nonlinear Stud.}, 16(2):381--387, 2016.

\bibitem{PariniWeth2015}
Enea Parini and Tobias Weth.
\newblock Existence, unique continuation and symmetry of least energy nodal
  solutions to sublinear {N}eumann problems.
\newblock {\em Math. Z.}, 280(3-4):707--732, 2015.

\bibitem{PeletiervanderVorst}
Lambertus Peletier and Robert C. A.~M. Van~der Vorst.
\newblock Existence and nonexistence of positive solutions of nonlinear
  elliptic systems and the biharmonic equation.
\newblock {\em Differential Integral Equations}, 5(4):747--767, 1992.

\bibitem{L2_7}
Benedetta Pellacci, Angela Pistoia, Giusi Vaira, and Gianmaria Verzini.
\newblock Normalized concentrating solutions to nonlinear elliptic problems.
\newblock {\em J. Differential Equations}, 275:882--919, 2021.

\bibitem{PengPengWang2016}
Shuangjie Peng, Yan-fang Peng, and Zhi-Qiang Wang.
\newblock On elliptic systems with {S}obolev critical growth.
\newblock {\em Calc. Var. Partial Differential Equations}, 55(6):Paper No. 142, 2016.

\bibitem{PengWangWang2019}
Shuangjie Peng, Qingfang Wang, and Zhi-Qiang Wang.
\newblock On coupled nonlinear {S}chr\"{o}dinger systems with mixed couplings.
\newblock {\em Trans. Amer. Math. Soc.}, 371(11):7559--7583, 2019.

\bibitem{MR3689156}
Dario Pierotti and Gianmaria Verzini.
\newblock Normalized bound states for the nonlinear {S}chr\"odinger equation in
  bounded domains.
\newblock {\em Calc. Var. Partial Differential Equations}, 56(5):Paper No. 133, 
  2017.

\bibitem{PistoiaSurvey}
Angela Pistoia.
\newblock The {L}japunov-{S}chmidt reduction for some critical problems.
\newblock In {\em Concentration Analysis and Applications to PDE - ICTS
  Workshop, Bangalore, January 2012}, Trends in Mathematics, pages 69--83.
  Springer Basel, 2013.

\bibitem{PistoiaRamosNeumann}
Angela Pistoia and Miguel Ramos.
\newblock Locating the peaks of the least energy solutions to an elliptic
  system with {N}eumann boundary conditions.
\newblock {\em J. Differential Equations}, 201(1):160--176, 2004.

\bibitem{PistoiaSaldanaTavares}
Angela Pistoia, Alberto Salda\~{n}a, and Hugo Tavares.
\newblock Existence of solutions to a slightly supercritical pure {N}eumann
  problem.
\newblock {\em SIAM J. Math. Anal.}, 55(4):3844--3887, 2023.

\bibitem{PistoiaSchieraTavares}
Angela Pistoia, Delia Schiera, and Hugo Tavares.
\newblock Existence of solutions on the critical hyperbola for a pure
  {L}ane-{E}mden system with {N}eumann boundary conditions.
\newblock {\em International Mathematics Research Notices}, 1, 745--803  2024.

\bibitem{PistoiaSoave}
Angela Pistoia and Nicola Soave.
\newblock On {C}oron's problem for weakly coupled elliptic systems.
\newblock {\em Proc. Lond. Math. Soc. (3)}, 116(1):33--67, 2018.

\bibitem{PistoiaSoaveTavares}
Angela Pistoia, Nicola Soave, and Hugo Tavares.
\newblock A fountain of positive bubbles on a {C}oron's problem for a
  competitive weakly coupled gradient system.
\newblock {\em J. Math. Pures Appl. (9)}, 135:159--198, 2020.

\bibitem{PistoiaTavares}
Angela Pistoia and Hugo Tavares.
\newblock Spiked solutions for {S}chr\"{o}dinger systems with {S}obolev
  critical exponent: the cases of competitive and weakly cooperative
  interactions.
\newblock {\em J. Fixed Point Theory Appl.}, 19(1):407--446, 2017.

\bibitem{RamosTavaresZou}
Miguel~Ramos, Hugo~Tavares, and Wenming~Zou.
\newblock A {B}ahri-{L}ions theorem revisited.
\newblock {\em Adv. Math.}, 222(6):2173--2195, 2009.

\bibitem{RamosTavaresTerracini}
Miguel Ramos, Hugo Tavares, and Susanna Terracini.
\newblock Extremality conditions and regularity of solutions to optimal
  partition problems involving {L}aplacian eigenvalues.
\newblock {\em Arch. Ration. Mech. Anal.}, 220(1):363--443, 2016.

\bibitem{RamosYang}
Miguel Ramos and Jianfu Yang.
\newblock Spike-layered solutions for an elliptic system with {N}eumann
  boundary conditions.
\newblock {\em Trans. Amer. Math. Soc.}, 357(8):3265--3284, 2005.

\bibitem{Rey1}
Olivier Rey.
\newblock Proof of two conjectures of {H}. {B}r\'{e}zis and {L}. {A}.
  {P}eletier.
\newblock {\em Manuscripta Math.}, 65(1):19--37, 1989.

\bibitem{Rogel_Salazar_2013}
Jesus Rogel-Salazar.
\newblock The {G}ross{\textendash}{P}itaevskii equation and
  {B}ose{\textendash}{E}instein condensates.
\newblock {\em European Journal of Physics}, 34(2):247--257, 2013.

\bibitem{RoseWeinstein88}
Harvey~A. Rose and Michael~I. Weinstein.
\newblock On the bound states of the nonlinear {S}chr\"odinger equation with a
  linear potential.
\newblock {\em Phys. D}, 30(1-2):207--218, 1988.

\bibitem{RoselliWillem}
Paolo Roselli and Michel Willem.
\newblock Least energy nodal solutions of the {B}rezis-{N}irenberg problem in
  dimension {$N=5$}.
\newblock {\em Commun. Contemp. Math.}, 11(1):59--69, 2009.

\bibitem{Ruf}
Bernhard Ruf.
\newblock Superlinear elliptic equations and systems.
\newblock In {\em Handbook of differential equations: stationary partial
  differential equations. {V}ol. {V}}, Handb. Differ. Equ., pages 211--276.
  Elsevier/North-Holland, Amsterdam, 2008.

\bibitem{SaldanaTavares}
Alberto Salda\~{n}a and Hugo Tavares.
\newblock Least energy nodal solutions of {H}amiltonian elliptic systems with
  {N}eumann boundary conditions.
\newblock {\em J. Differential Equations}, 265(12):6127--6165, 2018.

\bibitem{SaldanaTavaresNoDEA}
Alberto Salda\~{n}a and Hugo Tavares.
\newblock On the least-energy solutions of the pure {N}eumann {L}ane-{E}mden
  equation.
\newblock {\em NoDEA Nonlinear Differential Equations Appl.}, 29(3):Paper No.
  30, 2022.

\bibitem{SalsaPDE}
Sandro Salsa.
\newblock {\em Partial differential equations in action}, volume~99 of {\em
  Unitext}.
\newblock Springer, [Cham], third edition, 2016.
\newblock From modelling to theory, La Matematica per il 3+2.

\bibitem{SatoWang}
Yohei Sato and Zhi-Qiang Wang.
\newblock Least energy solutions for nonlinear {S}chr{\"o}dinger systems with
  mixed attractive and repulsive couplings.
\newblock {\em Adv. Nonlinear Stud.}, 15(1):1--22, 2015.

\bibitem{SchechterZou}
Martin Schechter and Wenming Zou.
\newblock On the {B}r\'{e}zis-{N}irenberg problem.
\newblock {\em Arch. Ration. Mech. Anal.}, 197(1):337--356, 2010.

\bibitem{Schoen}
Richard Schoen.
\newblock Conformal deformation of a {R}iemannian metric to constant scalar
  curvature.
\newblock {\em J. Differential Geom.}, 20(2):479--495, 1984.

\bibitem{Sirakov}
Boyan Sirakov.
\newblock On the existence of solutions of {H}amiltonian elliptic systems in
  {$\bold R^N$}.
\newblock {\em Adv. Differential Equations}, 5(10-12):1445--1464, 2000.

\bibitem{Sirakov2007}
Boyan Sirakov.
\newblock Least energy solitary waves for a system of nonlinear {S}chr\"odinger
  equations in {$\mathbb{R}^n$}.
\newblock {\em Comm. Math. Phys.}, 271(1):199--221, 2007.

\bibitem{Soave}
Nicola Soave.
\newblock On existence and phase separation of solitary waves for nonlinear
  {S}chr{\"o}dinger systems modelling simultaneous cooperation and competition.
\newblock {\em Calc. Var. Partial Differential Equations}, 53(3-4):689--718,
  2015.

\bibitem{L2_8}
Nicola Soave.
\newblock Normalized ground states for the {NLS} equation with combined
  nonlinearities.
\newblock {\em J. Differential Equations}, 269(9):6941--6987, 2020.

\bibitem{L2_9}
Nicola Soave.
\newblock Normalized ground states for the {NLS} equation with combined
  nonlinearities: the {S}obolev critical case.
\newblock {\em J. Funct. Anal.}, 279(6):108610, 2020.

\bibitem{SoaveTavares}
Nicola Soave and Hugo Tavares.
\newblock New existence and symmetry results for least energy positive
  solutions of {S}chr\"{o}dinger systems with mixed competition and cooperation
  terms.
\newblock {\em J. Differential Equations}, 261(1):505--537, 2016.

\bibitem{STTZsurvey}
Nicola Soave, Hugo Tavares, Susanna Terracini, and Alessandro Zilio.
\newblock H\"{o}lder bounds and regularity of emerging free boundaries for
  strongly competing {S}chr\"{o}dinger equations with nontrivial grouping.
\newblock {\em Nonlinear Anal.}, 138:388--427, 2016.

\bibitem{STTZ2}
Nicola Soave, Hugo Tavares, Susanna Terracini, and Alessandro Zilio.
\newblock Variational problems with long-range interaction.
\newblock {\em Arch. Ration. Mech. Anal.}, 228(3):743--772, 2018.

\bibitem{SoaveTavaresZilio}
Nicola Soave, Hugo Tavares, and Alessandro Zilio.
\newblock Free boundary problems with long-range interactions: uniform
  {L}ipschitz estimates in the radius.
\newblock {\em Math. Ann.}, 386(1-2):551--585, 2023.

\bibitem{SoaveTerracini}
Nicola Soave and Susanna Terracini.
\newblock The nodal set of solutions to some elliptic problems: singular
  nonlinearities.
\newblock {\em J. Math. Pures Appl. (9)}, 128:264--296, 2019.

\bibitem{SoaveWeth}
Nicola Soave and Tobias Weth.
\newblock The unique continuation property of sublinear equations.
\newblock {\em SIAM J. Math. Anal.}, 50(4):3919--3938, 2018.

\bibitem{SoaveZilio}
Nicola Soave and Alessandro Zilio.
\newblock Uniform bounds for strongly competing systems: the optimal
  {L}ipschitz case.
\newblock {\em Archive for Rational Mechanics and Analysis}, 218(2):647--697,
  2015.

\bibitem{Sverak}
Vladimir {\v{S}}ver{{\'a}}k.
\newblock On optimal shape design.
\newblock {\em J. Math. Pures Appl. (9)}, 72(6):537--551, 1993.

\bibitem{CIM}
Hugo Tavares.
\newblock Semilinear elliptic problems: old and new.
\newblock {\em CIM Bulletin}, 43:11--18, 2021.

\bibitem{TavaresTerracini1}
Hugo Tavares and Susanna Terracini.
\newblock Regularity of the nodal set of segregated critical configurations
  under a weak reflection law.
\newblock {\em Calc. Var. Partial Differential Equations}, 45(3-4):273--317,
  2012.

\bibitem{TavaresTerracini2}
Hugo Tavares and Susanna Terracini.
\newblock Sign-changing solutions of competition-diffusion elliptic systems and
  optimal partition problems.
\newblock {\em Ann. Inst. H. Poincar{\'e} Anal. Non Lin{\'e}aire},
  29(2):279--300, 2012.

\bibitem{TavaresWeth}
Hugo Tavares and Tobias Weth.
\newblock Existence and symmetry results for competing variational systems.
\newblock {\em NoDEA Nonlinear Differential Equations Appl.}, 20(3):715--740,
  2013.

\bibitem{TavaresYou}
Hugo Tavares and Song You.
\newblock Existence of least energy positive solutions to {S}chr\"{o}dinger
  systems with mixed competition and cooperation terms: the critical case.
\newblock {\em Calc. Var. Partial Differential Equations}, 59(1):Paper No. 26,
2020.

\bibitem{TavaresYouZou}
Hugo Tavares, Song You, and Wenming Zou.
\newblock Least energy positive solutions of critical {S}chr\"{o}dinger systems
  with mixed competition and cooperation terms: the higher dimensional case.
\newblock {\em J. Funct. Anal.}, 283(2):Paper No. 109497, 2022.

\bibitem{TavaresZilio}
Hugo Tavares and Alessandro Zilio.
\newblock Regularity of all minimizers of a class of spectral partition
  problems.
\newblock {\em Math. Eng.}, 3(1):Paper No. 2, 2021.

\bibitem{Timmermans}
Eddy Timmermans.
\newblock Phase separation of {B}ose-{E}instein condensates.
\newblock {\em Physical Review Letters}, 81(26):5718--5721, 1998.

\bibitem{Trudinger}
Neil~S. Trudinger.
\newblock Remarks concerning the conformal deformation of {R}iemannian
  structures on compact manifolds.
\newblock {\em Ann. Scuola Norm. Sup. Pisa Cl. Sci. (3)}, 22:265--274, 1968.

\bibitem{Wang1991}
Xu~Jia Wang.
\newblock Neumann problems of semilinear elliptic equations involving critical
  {S}obolev exponents.
\newblock {\em J. Differential Equations}, 93(2):283--310, 1991.

\bibitem{WangWillem}
Zhi-Qiang Wang and Michel Willem.
\newblock Partial symmetry of vector solutions for elliptic systems.
\newblock {\em J. Anal. Math.}, 122:69--85, 2014.

\bibitem{WeiWu}
Juncheng Wei and Yuanze Wu.
\newblock Ground states of nonlinear {S}chr\"{o}dinger systems with mixed
  couplings.
\newblock {\em J. Math. Pures Appl. (9)}, 141:50--88, 2020.

\bibitem{WethTMNA2006}
Tobias Weth.
\newblock Nodal solutions to superlinear biharmonic equations via decomposition
  in dual cones.
\newblock {\em Topol. Methods Nonlinear Anal.}, 28(1):33--52, 2006.

\bibitem{Wu}
Yuanze Wu.
\newblock On a {$K$}-component elliptic system with the {S}obolev critical
  exponent in high dimensions: the repulsive case.
\newblock {\em Calc. Var. Partial Differential Equations}, 56(5):Paper No. 151,
2017.

\bibitem{Yamabe}
Hidehiko Yamabe.
\newblock On a deformation of {R}iemannian structures on compact manifolds.
\newblock {\em Osaka Math. J.}, 12:21--37, 1960.

\bibitem{YinZou}
Xin Yin and Wenming Zou.
\newblock Positive least energy solutions for {$k$}-coupled {S}chr\"{o}dinger
  system with critical exponent: the higher dimension and cooperative case.
\newblock {\em J. Fixed Point Theory Appl.}, 24(1):Paper No. 5, 2022.

\bibitem{L2_4}
Zhaoyang Yun and Zhitao Zhang.
\newblock Normalized solutions to {S}chr\"{o}dinger systems with linear and
  nonlinear couplings.
\newblock {\em J. Math. Anal. Appl.}, 506(1):Paper No. 125564, 2022.

\bibitem{Zeng}
Jing Zeng.
\newblock The estimates on the energy functional of an elliptic system with
  {N}eumann boundary conditions.
\newblock {\em Bound. Value Probl.}, 2013:194, 013.

\end{thebibliography}
\EditInfo{November 18, 2023}{January 21, 2024}{Ana Cristina Moreira Freitas, Oliveira E Silva Diogo, Ivan Kaygorodov, Carlos Florentino}

\end{document}